\documentclass{article}

\usepackage[ansinew]{inputenc}
\usepackage{theorem}
\usepackage[english]{babel}
\usepackage{graphicx}
\usepackage{amsfonts}
\usepackage{amsmath}

\numberwithin{equation}{section}

\newtheorem {defn}{Definition}[section]
\newtheorem {thm}{Theorem}[section]
\newtheorem {prop}[thm]{Proposition}
\newtheorem {mthm}[thm]{Main Theorem}
\newtheorem {lemma}[thm]{Lemma}
\newtheorem {assump}[defn]{Assumption}
\newtheorem {cor}[thm]{Corollary}

\def\squarebox#1{\hbox to #1{\hfill\vbox to #1{\vfill}}}
\newcommand{\qed}{\hspace*{\fill}
\vbox{\hrule\hbox{\vrule\squarebox{.667em}\vrule}\hrule}\smallskip}

\newcommand{\bbmatrix}[1]{\left[ \begin{array}{cccccccccc} #1 \end{array} \right]}

\title{Overdetermined conservative $2D$ Systems, Invariant in One Direction and a Generalization of Potapov's theorem}
\author{Andrey Melnikov, Victor Vinnikov \\
{\small Ben-Gurion University of the Negev, Israel }}
\date{19/12/2008}
\begin{document}
\maketitle

\begin{abstract}
This work is a direct continuation of the authors work \cite{bib:NonConsTheory}. 
A special case of conservative overdetermined time invariant $2D$ systems is developed and studied.
Defining \textit{transfer function} of such systems we obtain a class $\boldsymbol{\mathcal CI}$ of inner functions $S(\lambda,t_2)$,
which are identity for $\lambda=\infty$, satisfy certain regularity assumptions and intertwines solutions of ODEs with a spectral parameter $\lambda$. Using translation model, developed in \cite{bib:TheoryNonComm} we prove that every function in the class $\boldsymbol{\mathcal CI}$ can be 
realized as a transfer function of a certain vessel.

The highlight of this theory is a generalization of Potapov's theorem \cite{bib:Potapov, bib:SpectrAnal}, which gives a very special formula for such
a function in the form of multiplication of Blacke-Potapov products, corresponding to the discrete spectrum of certain system operator $A_1(t_2)$
and of multiplicative integral, corresponding to the continuous spectrum of $A_1(t_2)$. This theorem is proved under a slightly more restrictive
assumptions, then the development of the whole theory. Namely, we suppose that the derivative of the transfer function is a continuous 
function of $t_2$ for almost all $\lambda$.

At the last part zero/pole interpolation problem \cite{bib:Inter} for matrix functions in $\boldsymbol{\mathcal CI}$ is considered and a realization
theorem of such functions appeared in \cite{bib:NonConsTheory} (theorem 8.1) is reproved. Hermitian case is also analyzed and the corresponding
realization theorem is proved.

\noindent\textbf{MSC classes:}  47N20, 47N70, 47E05, 26B30, 45D05, 34K06.

\noindent\textbf{Keywords:} system theory, vessel, colligation, inner function, Volterra operator.
\end{abstract}

\tableofcontents
\section{Introduction}
Generalizing time varying $1D$ systems to the study of $2D$ systems invariant in one direction, it turns out to be very useful
to introduce the notion of vessel. It can actually be done in different settings, for example in algebraic sense \cite{bib:Gauchman,
bib:VinnMSRI,bib:AlBallPerez} with inputs/output satisfying certain algebraic equations, or in analytic \cite{bib:Vortices, bib:Over2Dsys}
with ODEs with spectral parameter for input/output. The second approach is our main inspiration and mainly comes from the article of M. 
Liv\v sic \cite{bib:Vortices}.

In the 40-50's \cite{bib:Potapov, bib:JordanBrodskii, bib:SpectrAnal} there was developed a theory connecting non self-adjoint operator $A$
(with a small imaginary part $A-A^*$) and meromorphic functions in the upper half plane (or inside circle), called characteristic 
functions of the operator. Multiplicative structure of characteristic function was closely related to the invariant subspace 
structure of the operator $A$. See survey \cite{bib:BallCohen} on this subject. Further, this theory was developed for a pair
of commuting non self-adjoint operators \cite{bib:TheoryNonComm}. A special and important role in this research is played by conservative
systems, which are closely related to the study of the operators.

One of the strongest results in this are was obtained in \cite{bib:SpectrAnal}, where multiplicative structure of a meromorphic function
was connected to invariant subspaces of the corresponding operator $A$. Moreover, the Hilbert space was decomposed into two parts 
$\mathcal H_1 \oplus \mathcal H_2$, where $\mathcal H_1 = l^2$ - the space oh infinite sequences and $\mathcal H_2 = L^2(0,L)$. Moreover,
$A_1|_{\mathcal H_1}$ had a triangular form and $A_1|_{\mathcal H_2}$ was a Volterra and multiplication operator on a continuous from the
right function.

V.P. Potapov \cite{bib:Potapov} proved this theorem for the characteristic function using functional analysis approach only.

We are going to generalize this result to a wider class of function $\boldsymbol{\mathcal CI}$ and this is done in Part II.
Part I presents the theory of conservative vessels, based on the work \cite{bib:NonConsTheory}. We show first how such vessel arise, then
gauge equivalence of such vessels is presented. As in the non conservative case \cite{bib:NonConsTheory} one can differentiate and integrate
vessels and actually these two categories are equivalent. Main theorem of gauge equivalence is presented in theorem \ref{mthm:CGE}. 
Afterwords, transfer function and the class $\boldsymbol{\mathcal CI}$ it belongs to are defined. Translation model, which is used for the
main realization theorem of part I is presented, followed by theorem itself.

In part III we solve zero/pole interpolation problem \cite{bib:Inter} and show its applications.

\section{Background}
In our previous article \cite{bib:NonConsTheory} we have developed and studied a general theory of Vessels and corresponding overdetermined $2D$ time invariant systems. Let us recall main definitions and assumptions. An overdetermined $t_1$-invariant $2D$ system is a linear input-state-output 
(i/s/o) system, consisting of operators depending only on the variable $t_2$; in the most general case such a system is of the form \cite{bib:Vortices}
\begin{equation} \label{eq:systempre}
  I\Sigma': \left\{ \begin{array}{lll}
         x(t_1, t_2) = e^{A_1(t_2)(t_1-t_1^0)} x(t_1^0, t_2) + \int\limits_{t_1^0}^{t_1}
       e^{A_1(t_2)(t_1-y)} B_1(t_2) u(y, t_2) dy \\
    x(t_1, t_2) = F (t_2,t_2^0) x(t_1, t_2^0) + \int\limits_{t_2^0}^{t_2} F(t_2, s) B_2(s) u(t_1, s)ds \\
    y(t_1, t_2) = C(t_2) x(t_1, t_2) + D(t_2) u(t_1, t_2)
  \end{array} \right.
  \end{equation}
where for Hilbert spaces $\mathcal E, \mathcal E_*, \mathcal H_{t_2}$ there are defined
\[ \begin{array}{lll}
u(t_1,t_2) \in \mathcal E \text{ - input,} \\
y(t_1,t_2) \in \mathcal E_* \text{ - output,} \\
x(t_1,t_2) \in \mathcal H_{t_2} \text{ - state,}
\end{array}
\]
such that $u(t_1,t_2), y(t_1,t_2)$ are absolutely continuous functions of each variable when the other variable is fixed. The
transition of the system will usually be considered from $(t_1^0,t_2^0)$ to $(t_1, t_2)$. Note that $\mathcal H_{t_2}$ are a priory
different for each $t_2$, and as a result $F(t_2,t_2^0)$ has to be an evolution semi-group , i.e., satisfies the following 
\begin{defn} \label{def:EvolSemGr} Given a collection of Hilbert spaces $\{ \mathcal H_t \mid t\in I\}$ for an interval $I\subseteq\mathbb R$
and a collection of bounded invertible operators $F(s,t):\mathcal H_s \rightarrow \mathcal H_t$ for each
$s,t \in I$, we will say that $F(s,t)$ is \textbf{evolution semi-group} if the following relations
hold for all $r,t,s \in I$:
\[ \begin{array}{llll}
F(r,s) F(s,t) = F(r,t), \\
F(t,t) = Id|_{\mathcal H_t}.
\end{array}\]
\end{defn}
Demanding compatibility of transition for the system (\ref{eq:systempre}) and factorization 
\begin{equation} \label{eq:Factorization}
\begin{array}{lllll}
B_2(t_2) = \widetilde B(t_2) \sigma_2(t_2), ~~B_1(t_2) = \widetilde B(t_2) \sigma_1(t_2), \\
A_1(t_2) B_2(t_2) + F(t_2, t_2^0)\frac{\partial}{\partial s} [F (t_2^0, t_2) B_1(t_2)] = - \widetilde B(t_2) \gamma(t_2)
\end{array} \end{equation}
for some operators
\[ \widetilde B(t_2): \widetilde{\mathcal E} \rightarrow \mathcal H_{t_2}, ~~
\sigma_2(t_2), \sigma_1(t_2), \gamma(t_2): \mathcal E \rightarrow \widetilde{\mathcal E},
\]
where $ \widetilde{\mathcal E}$ is another auxiliary Hilbert space, we have reached the notion of \textbf{(integrated) vessel},
which is a collection of operator and spaces
\begin{multline*}
\mathfrak{IV} = (A_1(t_2), F(t_2,t_2^0),\widetilde B(t_2),C(t_2), D(t_2),\widetilde D(t_2); \\
\sigma_1(t_2), \sigma_2(t_2), \gamma(t_2), \sigma_{1*}(t_2), \sigma_{2*}(t_2)\gamma_*(t_2); \mathcal H_{t_2}, \mathcal{E},\mathcal E_*, \mathcal{\widetilde E},\mathcal{\widetilde E_*})
\end{multline*}
satisfying regularity assumptions 
\begin{assump} \label{assm:Regularity} 
	\begin{itemize}
	\item \textbf{Internal} regularity:
    \begin{enumerate}
    \item $A_1(t_2):\mathcal H_{t_2} \rightarrow \mathcal H_{t_2}$,
       $B_1(t_2), B_2(t_2): \mathcal E \rightarrow \mathcal H_{t_2}$, $ C(t_2): \mathcal H_{t_2} \rightarrow \mathcal E_*$
       are bounded operators (for all $t_2$) and
       $F(t_2,t_2^0): \mathcal H_{t_2^0} \rightarrow \mathcal{H}_{t_2}$ is an evolution semi-group
       (see definition \ref{def:EvolSemGr}).
    \item $F(t_2, s) B_2(s)$ and $C(s) F(s, t_2)$ are absolutely continuous as functions
       of $s$ (for almost all $t_2$) in the norm operator topology on $\mathcal L(\mathcal E, \mathcal H_{t_2})$
       and on $\mathcal L(\mathcal E_*, \mathcal H_{t_2})$, respectively.
    \end{enumerate}
	\item \textbf{Feed through} regularity: the operators $D(t_2): \mathcal E \rightarrow \mathcal E_*$ and 
		$\widetilde D(t_2):  \mathcal E_* \rightarrow \mathcal{\widetilde E_*}$ are absolutely continuous functions of $t_2$.
	\item \textbf{External input} regularity:
        \begin{enumerate}
        \item $\gamma(t_2), \sigma_2(t_2) \in L^1_{loc}( \mathcal L(\mathcal E, \widetilde{\mathcal E}))$
                in the norm operator topology.
        \item $\sigma_1(t_2) \in L(\mathcal E, \widetilde{\mathcal E})$ is absolutely continuous and invertible,
                        in the norm operator topology.
        \end{enumerate}
	\item \textbf{External output} regularity:
        \begin{enumerate}
        \item $\gamma_*(t_2), \sigma_{2*}(t_2) \in L^1_{loc}( \mathcal L(\mathcal E_*, \widetilde{\mathcal E_*}))$
                in the norm operator topology.
        \item $\sigma_{1*}(t_2) \in L(\mathcal E_*, \widetilde{\mathcal E_*})$ is absolutely continuous and invertible,
                        in the norm operator topology.
        \end{enumerate}
	\end{itemize}
\end{assump}
and the following vessel conditions: Lax condition
\renewcommand{\theequation}{\arabic{section}.\arabic{equation}.Lax}
\begin{equation} \label{eq:LaxCond} F(t_2,t_2^0) A_1(t_2^0) = A_1(t_2) F(t_2,t_2^0) \end{equation}
Input vessel condition
\renewcommand{\theequation}{\arabic{section}.\arabic{equation}.OverD}
\begin{equation} \label{eq:OverDetCondIn} 
\frac{d}{dt_2} (F (t_2^0, t_2) \widetilde B(t_2) \sigma_1(t_2)) + F(t_2^0, t_2) A_1(t_2) \widetilde B(t_2) \sigma_2(t_2)
   + F(t_2^0, t_2) \widetilde B(t_2) \gamma(t_2) = 0
\end{equation}
Output vessel condition
\begin{equation} \label{eq:OverDetCondOut} - \sigma_{1*}(t_2) \frac{d}{d t_2}[C(t_2) F(t_2,t_2^0)] + \sigma_{2*}(t_2) C(t_2) A_1(t_2) F(t_2,t_2^0)  
        \gamma_*(t_2) C(t_2) F(t_2,t_2^0) = 0 \end{equation}
Linkage condition
\renewcommand{\theequation}{\arabic{section}.\arabic{equation}.Link}
\begin{equation} \label{eq:LinkCond} \begin{array}{ll}
\sigma_{1*} D = \widetilde D \sigma_1, ~~~~~~ \sigma_{2*} D = \widetilde D \sigma_2,  \\
\widetilde D \gamma = \sigma_{2*} C \widetilde B \sigma_1  - \sigma_{1*} C \widetilde B \sigma_2 + \sigma_{1*} D' + \gamma_*D.
\end{array} \end{equation}
It is naturally associated to the system $I\Sigma$ (see (\ref{eq:systempre}))
\renewcommand{\theequation}{\arabic{section}.\arabic{equation}}
\begin{equation} \label{eq:system}
     I\Sigma: \left\{ \begin{array}{lll}
    \frac{\partial}{\partial t_1}x(t_1,t_2) = A_1(t_2) ~x(t_1,t_2) + \widetilde B(t_2) \sigma_1(t_2) ~u(t_1,t_2) \\[5pt]
    x(t_1, t_2) = F (t_2,t_2^0) x(t_1, t_2^0) + \int\limits_{t_2^0}^{t_2} F(t_2, s) \widetilde B(s) \sigma_2(s) u(t_1, s)ds \\[5pt]
    y(t_1,t_2) = C(t_2)~ x(t_1,t_2) + D(t_2) u(t_1,t_2).
  \end{array} \right.
  \end{equation}
with absolutely continuous inputs and outputs, satisfying compatibility conditions for almost all $(t_1,t_2)$:
\begin{eqnarray} \label{eq:InCC}
  \sigma_2(t_2) \frac{\partial}{\partial t_1}u(t_1, t_2) -
  \sigma_1(t_2) \frac{\partial}{\partial t_2}u(t_1,t_2) + \gamma(t_2) u(t_1,t_2) = 0, \\
  \label{eq:OutCC}
  \sigma_{2*}(t_2) \frac{\partial}{\partial t_1}y(t_1, t_2) -
  \sigma_{1*}(t_2) \frac{\partial}{\partial t_2}y(t_1,t_2) + \gamma_*(t_2) y(t_1,t_2) = 0.
\end{eqnarray}
Another important notion, which will be extensively used is the notion of \textit{adjoint vessel} \cite{bib:NonConsTheory}.
It is obtained from a simple observation that applying adjoint to the vessel conditions gives rise to a new set of conditions on
adjoint operators, which are almost vessel conditions. An \textit{adjoint systems} (for $I\Sigma$) is the sytem
\begin{equation} \label{eq:IsystemAdj}
    I\Sigma^*: \left\{ \begin{array}{lll}
    -\frac{\partial}{\partial t_1}x_*(t_1,t_2) = A_1^*(t_2) ~x_*(t_1,t_2) + C^*(t_2) \sigma^*_{1*} u_*(t_1,t_2) \\[5pt]
    -x_*(t_1,t_2) = F^*(t_2,t_2^0)~x_*(t_1^0,t_2^0) + \int\limits_{t_2^0}^{t_2} F^*(t_2, s) C^*(s) \sigma^*_{2*}(s) u_*(t_1,s) ds\\[5pt]
    y_*(t_1,t_2) = \widetilde B^*~ x_*(t_1,t_2)  + \widetilde D^*(t_2)~ u_*(t_1,t_2),
    \end{array} \right.
\end{equation}
which is associated to the vessel $\mathfrak{V^*}$ given by
\[ \mathfrak{V^*} =
(-A_1^*, -F^*(t_2,t_2^0), -C^*, \widetilde B^*, \widetilde D^*, D^*;
\sigma^*_{1*}, \sigma^*_{2*}, -\gamma_*^* - \frac{d}{dt_2}
\sigma_{1*}^*, \sigma^*_1, \sigma^*_2, -\gamma^* - \frac{d}{dt_2}
\sigma_1^*; \mathcal{H}_{t_2}, \widetilde{\mathcal E}_*,
\widetilde{\mathcal E}_*, \mathcal{E}_*, \mathcal{E}),
\]
where all the operators are functions of $t_2$ and satisfy the following axioms:
\[ \begin{array}{lllllll}
		\frac{d}{dt_2} A_1^* = A_1^* A_2^* - A_2^* A_1^* \\
		\frac{d}{dt_2} \big(F^* C^* \sigma^*_{1*} \big) - A_1^* F^* C^* \sigma^*_{2*} -   F^*C^* (\gamma_*^*+\frac{d}{dt_2} \sigma_{1*}^*) = 0 \\
		- \sigma_1^* \frac{d}{dt_2} \big(\widetilde B^* F^*\big) + \sigma_2^* \widetilde B F^* A_1^* - (\gamma^* + \frac{d}{dt_2} \sigma_1)\widetilde B^* F^*= 0 \\
		\sigma^*_1 \widetilde D^* = D^* \sigma^*_{1*}, ~~ \sigma^*_2 \widetilde D^* = D^* \sigma^*_{2*} \\
		D^* (-\gamma^*_* -\frac{d}{dt_2} \sigma_{1*}^*) = -\sigma^*_2 \widetilde B^* C^* \sigma^*_{1*} +
		\sigma^*_1 \widetilde B^* C^* \sigma^*_{2*} - \sigma^*_1 \frac{d}{dt_2} \widetilde D^* - (\gamma^* +\frac{d}{dt_2} \sigma_1^*)\widetilde D^*.
\end{array} \]
Moreover, the transfer function of the adjoint vessel $S_*(\mu,t_2)$ maps solutions of the \textit{adjoint input} ODE
\begin{equation} \label{eq:InCC*} [\sigma^*_{2*} \mu - \sigma^*_{1*} \dfrac{d}{dt_2} -
\gamma^*_* - \frac{d}{dt_2} \sigma_{1*}^*] u_*(t_2) = 0
\end{equation}
with the spectral parameter $\mu$ to solutions of the \textit{adjoint output} ODE
\begin{equation} \label{eq:OutCC*} [\sigma^*_2 \mu - \sigma^*_1 \dfrac{d}{dt_2} -
\gamma^*-\frac{d}{dt_2} \sigma_1^*] y_*(t_2) = 0
\end{equation}
with the same spectral parameter. The following relation between transfer functions has to be satisfied:
\begin{equation} \label{eq:SS*}
 S(\lambda, t_2) =  \sigma_{1*}^{-1} S^*_*(-\bar\lambda, t_2) \sigma_1
\end{equation}
One can easily verify this formula directly, using the vessel conditions and the formulas for
$S(\lambda, t_2)$, $S^*_*(-\bar\lambda, t_2)$:
\[ \begin{array}{lllllll}
  \sigma_{1*}^{-1} S^*_*(-\bar\lambda, t_2) \sigma_1 =
  \sigma_{1*}^{-1} \big[ \widetilde D^* - \widetilde B^* (-\bar\lambda I +
A_1^*)^{-1} C^* \sigma^*_{1*} \big]^* \sigma_1 = \\
 ~~~~ = \sigma_{1*}^{-1} [\widetilde D - \sigma_{1*} C(-\lambda I + A_1)^{-1} \widetilde B]\sigma_1 =
    \sigma_{1*}^{-1} \widetilde D \sigma_1 - C(-\lambda I + A_1)^{-1} \widetilde B \sigma_1 = \\
 ~~~~ = D + C (\lambda I - A_1)^{-1} \widetilde B \sigma_1 = S(\lambda, t_2).
\end{array} \]
\part{Conservative systems}
\section{Conservative systems and vessels \label{sec:Cons} }
In this part we want to talk about a more restrictive class of
vessels (and systems), satisfying a conservation law (or energy
balance). We start from the system (\ref{eq:system})
\[
     I\Sigma: \left\{ \begin{array}{lll}
    \frac{\partial}{\partial t_1}x(t_1,t_2) = A_1(t_2) ~x(t_1,t_2) + \widetilde B(t_2) \sigma_1(t_2) ~u(t_1,t_2) \\[5pt]
    x(t_1, t_2) = F (t_2,t_2^0) x(t_1, t_2^0) + \int\limits_{t_2^0}^{t_2} F(t_2, s) \widetilde B(s) \sigma_2(s) u(t_1, s)ds \\[5pt]
    y(t_1,t_2) = C(t_2)~ x(t_1,t_2) + D(t_2) u(t_1,t_2).
  \end{array} \right.
\]
and consider its adjoint (\ref{eq:IsystemAdj})
\[
    I\Sigma^*: \left\{ \begin{array}{lll}
    -\frac{\partial}{\partial t_1}x_*(t_1,t_2) = A_1^*(t_2) ~x_*(t_1,t_2) + C^*(t_2) \sigma^*_{1*}(t_2) u_*(t_1,t_2) \\[5pt]
    -x_*(t_1,t_2) = F^*(t_2,t_2^0)~x_*(t_1^0,t_2^0) + \int\limits_{t_2^0}^{t_2} F^*(t_2, s) C^*(s) \sigma^*_{2*}(s) u_*(t_1,s) ds\\[5pt]
    y_*(t_1,t_2) = \widetilde B^*~ x_*(t_1,t_2)  + \widetilde D^*(t_2)~ u_*(t_1,t_2).
    \end{array} \right.
\]
We define $I\Sigma$ to be (scattering) \textit{conservative} if the transformation
\[ (u,x,y) \rightarrow (y,x,u)
\]
is a bijection from the set of trajectories for the system (\ref{eq:system}) to the set of trajectories $(u_*,x_*,y_*)$ for the system
(\ref{eq:IsystemAdj}). Note that a preliminary necessary condition for this to be possible is that the space of inputs, satisfying
(\ref{eq:InCC}), matches with the space of outputs for (\ref{eq:IsystemAdj}), satisfying \textit{adjoint output} ODE (\ref{eq:IsystemAdj}).
Thus
\[ \mathcal E = \widetilde{\mathcal E}, ~\sigma_1(t_2) = \sigma_1^*(t_2),
~\sigma_2(t_2) = \sigma_2^*(t_2), ~\gamma^*(t_2) = -\gamma(t_2) -\frac{d}{dt_2} \sigma_1.
\]
The same conditions for the outputs of $I\Sigma$ and inputs of $I\Sigma^*$ results in
\[ \mathcal E_* = \widetilde{\mathcal E_*}, ~\sigma_{1*}(t_2) = \sigma_{1*}^*(t_2), ~\sigma_2(t_2) = \sigma_{2*}^*(t_2), ~\gamma_*^*(t_2) = -\gamma_*(t_2) -\frac{d}{dt_2} \sigma_1.
\]
An immediate consequence of such ’’adjoint pairing'' for system trajectories is the following set of energy-balance relations for
conservative systems.
\begin{thm} Suppose that (\ref{eq:system}) is a conservative system and $(u,x,y)$ is a trajectory for it. Then
\[ \frac{\partial}{\partial t_1} \langle x,x\rangle_{\mathcal H_{t_2}} + \langle\sigma_{1*} y,y\rangle_\mathcal{E} =
\langle\sigma_1 u,u\rangle_\mathcal{E}
\]
and
\[ \langle x(t_1,t_2), x(t_1,t_2) \rangle_{\mathcal{H}_{t_2}} - \langle x(t_1,t_2^0), x(t_1,t_2^0) \rangle_{\mathcal{H}_{t_2^0}} =
  \int_{t_2^0}^{t_2}[\langle\sigma_2(s) u(t_1,s),u(t_1,s)\rangle_\mathcal{E} -
     \langle\sigma_{2*}(s) y(t_1,s), y(t_1,s) \rangle_\mathcal{E_*} ds.
\]
\end{thm}
\textbf{Proof:} Let us perform the necessary calculations (notice that by the definition $x=x_*$):
\[ \begin{array}{lllllllll}
\frac{\partial}{\partial t_1} \langle x,x\rangle & = \frac{\partial}{\partial t_1} \langle x,x_*\rangle =
    \langle\frac{\partial}{\partial t_1} x,x_*\rangle + \langle x, \frac{\partial}{\partial t_1} x_* \rangle =
    \langle A_1 x + \widetilde B \sigma_1 u,x_*\rangle - \langle x, A_1^* x_* + C^* \sigma^*_{1*} u_* \rangle = \\
 & = \langle \sigma_1 u, \widetilde B^* x_*\rangle - \langle C x,\sigma^*_{1*} u_* \rangle =
    \langle \sigma_1u, y_*-\widetilde D^* u_*,\rangle - \langle y-D u,\sigma^*_{1*} u_* \rangle = \\
 & = \langle \sigma_1u, y_* \rangle - \langle \sigma_{1*} y,u_* \rangle +
    \langle (-\widetilde D \sigma_1 + \sigma_{1*} D)u,u_* \rangle =
    \langle \sigma_1u, y_* \rangle - \langle \sigma_{1*} y,u_* \rangle = \\
 & =\langle \sigma_1u, u \rangle - \langle \sigma_{1*} y,y \rangle \text{, since by definition $y_*=u, u_*=y$}
\end{array} \]
Performing the same calculations, but in the integral form with the indexes $1$ and $2$ interchanged, we shall obtain the
second equality.
\qed

So, let us consider the first equation:
\[ \frac{\partial}{\partial t_1} \langle x,x\rangle_{\mathcal H_{t_2}} + \langle\sigma_1 y,y\rangle_\mathcal{E} =
\langle\sigma_1 u,u\rangle_\mathcal{E}.\]
If one substitutes next the formulas for the derivative of $x$ according to $t_1$ (from (\ref{eq:systempre})),
it follows that the following must hold:
\[ \begin{array}{lll}
\langle A_1 x, x \rangle_\mathcal{H} + \langle x, A_1 x \rangle_\mathcal{H} +
\langle C^* \sigma_{1*} C x,x\rangle_\mathcal{H} + \\
~~~~~+ \langle \widetilde B \sigma_1 u ,x \rangle_\mathcal{H} + \langle \sigma_{1*} D u , Cx \rangle_\mathcal{H} +
\langle x, \widetilde B\sigma_1 u \rangle_\mathcal{H} + \langle C x, \sigma_{1*} D u \rangle_\mathcal{H} + \\
~~~~~~~~~~+ \langle D^* \sigma_{1*} D u , u \rangle_\mathcal{E} = \langle \sigma_1 u , u \rangle_\mathcal{E}.
\end{array} \]
We want this energy balance to hold for arbitrary input and state vectors, thus we obtain
\begin{equation*}
\begin{array}{lllll}
A_1(t_2) + A_1^*(t_2) + C^*(t_2) \sigma_{1*}(t_2) C(t_2) = 0 \\
\widetilde B\sigma_1(t_2) + C^*(t_2) \sigma_{1*}(t_2) D(t_2) = 0 \\
\sigma_1(t_2) = D^*(t_2) \sigma_{1*}(t_2) D(t_2)
\end{array}
\end{equation*}

A similar computation is correct for the second energy balance condition. We rewrite it in the integrated form, since the second system equation is such, and obtain:
\begin{align}
\label{eq:energ2}
\langle x(t_1,t_2), x(t_1,t_2) \rangle_\mathcal{H} - \langle x(t_1,t_2^0), x(t_1,t_2^0) \rangle_\mathcal{H} =
  \int_{t_2^0}^{t_2}[\langle\sigma_2(s) u(t_1,s),u(t_1,s)\rangle_\mathcal{E} -
     \langle\sigma_{2*}(s) y(t_1,s), y(t_1,s) \rangle_\mathcal{E} ds.
\end{align}
Substituting further formulas for $x(t_1,t_2)$ and for $y(t_1,s)$, we obtain
\[
\langle x(t_1,t_2),x(t_1,t_2)\rangle_{\mathcal H_{t_2}} - \langle x(\tau_1,t_2),x(\tau_1,t_2)\rangle_\mathcal{H} =
\int_{t_2^0}^{t_2}[\langle\sigma_2(s) u(t_1,s),u(t_1,s)\rangle_\mathcal{E} -
     \langle\sigma_{2*}(s) y(t_1,s), y(t_1,s) \rangle_\mathcal{E} ds,
\]
or, what is the same
\begin{multline*}
\| F(t_2, t_2^0) x(t_1, t_2^0) + \int\limits_{t_2^0}^{t_2}
       F(t_2, s) \widetilde B\sigma_2(s) u(t_1, s) ds\|_{\mathcal H_{t_2}} - \| x(t_1, t_2^0)\|_\mathcal{H} = \\
= \int_{t_2^0}^{t_2}[\langle\sigma_2(s) u(t_1,s),u(t_1,s)\rangle_\mathcal{E} -
     \langle\sigma_{2*}(s)(D(s) u(t_1, s) + C(s) x(t_1, s)),D(s) u(t_1, s) + C(s) x(t_1, s)  \rangle_\mathcal{E} ds,
\end{multline*}
which also have to hold for arbitrary input and state vectors, thus we obtain the following conditions
\renewcommand{\theequation}{\arabic{equation}.Coll}
\begin{equation*} \begin{array}{lllll}
\| F(t_2,t_2^0) x(t_1,t_2^0)\| - \| x(t_1,t_2^0)\| = \int_{t_2^0}^{t_2}
\langle \sigma_2(s) C(s) x(t_1, s), C(s) x(t_1, s) \rangle ds\\
\widetilde B\sigma_2(t_2) + C^*(t_2) \sigma_{2*}(t_2) D(t_2) = 0 \\
\sigma_2(t_2) = D^*(t_2) \sigma_{2*} (t_2)D(t_2)
\end{array} \end{equation*}
We assume next that for all $t_2$ there exist $\xi_1, \xi_2 \in \mathbb{C}$ such that $ \xi_1 \sigma_1(t_2) + \xi_2 \sigma_2(t_2) > \epsilon > 0$ (or just $\operatorname{det} \big( \xi_1 \sigma_1(t_2) + \xi_2 \sigma_2(t_2) \big)\neq 0$ in case $\operatorname{dim}\mathcal{E} < \infty$). Then from the equation
\[ \xi_1 \sigma_1(t_2) + \xi_2 \sigma_2(t_2) = D^*(t_2)
\big( \xi_1 \sigma_1(t_2) + \xi_2 \sigma_2(t_2)\big) D(t_2),
\]
it follows that the operator $D(t_2)$ is invertible, and defining a new output
\[ \bar{y} = D^{-1} y = D^{-1} (Cx + D u) = D^{-1} C x + u
\]
we obtain operators $\bar{D} = I$ and $\bar{C} = D^{-1} C$. So, we suppose that $D = \widetilde D = I$ and as a result
part of the linkage conditions become $\sigma_1 = \sigma_1 \widetilde D = D \sigma_{1*} = \sigma_{1*}$ and
$\sigma_2 = \sigma_2 \widetilde D = D \sigma_{2*} = \sigma_{2*}$, which means that $\mathcal E=\mathcal E_*$ and that
input (\ref{eq:InCC}) and output (\ref{eq:OutCC})
compatibility conditions differ only by $\gamma, \gamma_*$ and have the same $\sigma$'s. Consequently, $C(t_2)=B^*(t_2)$ if at least one
of $\sigma_1(t_2), \sigma_2(t_2)$ is invertible for almost all $t_2$. We shall usually suppose that $\sigma_1(t_2)$ is invertible,
promising a uniqueness of solutions for the input and for the output ODEs (\ref{eq:InCC}), (\ref{eq:OutCC}).

Thus, without loss of generality for conservative systems we may assume that the system
is actually of the form
\renewcommand{\theequation}{\arabic{equation}}
\begin{equation} \label{eq:systemCons}
    CI\Sigma: \left\{ \begin{array}{lll}
    \frac{\partial}{\partial t_1}x(t_1,t_2) = A_1(t_2) ~x(t_1,t_2) + B(t_2) \sigma_1(t_2) ~u(t_1,t_2) \\[5pt]
    x(t_1, t_2) = F (t_2,t_2^0) x(t_1, t_2^0) + \int\limits_{t_2^0}^{t_2} F(t_2, s) B(s) \sigma_2(s) u(t_1, s)ds \\[5pt]
    y(t_1,t_2) = u(t_1,t_2) - B^*(t_2)~ x(t_1,t_2),
    \end{array} \right.
\end{equation}
and the following two conditions are imposed:
\renewcommand{\theequation}{\arabic{section}.\arabic{equation}.Coll}
\begin{eqnarray} \label{eq:ConsCond}
A_1(t_2) + A_1^*(t_2) + B^*(t_2) \sigma_1(t_2) B(t_2) = 0 \\
\label{eq:ConsCond2}
\| F(t_2,t_2^0) x(t_1,t_2^0)\| - \| x(t_1,t_2^0)\| = \int_{t_2^0}^{t_2}
\langle \sigma_2(s) B(s) x(t_1, s), B(s) x(t_1, s) \rangle ds
\end{eqnarray}
Systems with energy balance conditions are sometimes called \textit{colligations}. So, we shall
refer the conditions (\ref{eq:ConsCond}) and (\ref{eq:ConsCond2}) as \textit{colligation conditions}.

In this manner we obtain the notion of \textit{conservative vessel}
\[ \mathfrak{CIV} = (A_1(t_2), F(t_2,t_2^0),B(t_2); \sigma_1(t_2), \sigma_2(t_2), \gamma(t_2), \gamma_*(t_2);
\mathcal H_{t_2}, \mathcal{E})
\]
satisfying regularity assumptions and the following vessel conditions:
\[ \begin{array}{lllllll}
    A_1(t_2) + A_1^*(t_2) + B^*(t_2) \sigma_1(t_2) B(t_2) = 0 & (\text{\ref{eq:ConsCond}}) \\
    \| F(t_2,t_2^0) x(t_1,t_2^0)\| - \| x(t_1,t_2^0)\| = \int_{t_2^0}^{t_2}
\langle \sigma_2(s) B^*(s) x(t_1, s), B(s) x(t_1, s) \rangle ds & (\text{\ref{eq:ConsCond2}}) \\
 F(t_2,t_2^0) A_1(t_2^0) = A_1(t_2) F(t_2,t_2^0)  & (\text{\ref{eq:LaxCond}}) \\
  \frac{d}{dt_2} (F (t_2^0, t_2) B(t_2) \sigma_1(t_2)) + F(t_2^0, t_2) A_1(t_2) B(t_2) \sigma_2(t_2)
   + F(t_2^0, t_2) B(t_2) \gamma(t_2) = 0
              & (\text{\ref{eq:OverDetCondIn}}) \\
  \sigma_1(t_2) \frac{\partial}{\partial t_2}[B^*(t_2) F(t_2,t_2^0)] - \sigma_2(t_2) B^*(t_2) A_1(t_2) F(t_2,t_2^0) -
        \gamma_*(t_2) B^*(t_2) F(t_2,t_2^0) = 0             & (\text{\ref{eq:OverDetCondOut}}) \\
 \gamma = \sigma_2 B^* B \sigma_1  - \sigma_1 B^* B \sigma_2 + \gamma_*. & (\text{\ref{eq:LinkCond}})
\end{array} \]
It is naturally associated with the system (\ref{eq:systemCons})
\begin{equation*}
     CI\Sigma: \left\{ \begin{array}{lll}
    \frac{\partial}{\partial t_1}x(t_1,t_2) = A_1(t_2) ~x(t_1,t_2) + B(t_2) \sigma_1(t_2) ~u(t_1,t_2) \\[5pt]
    x(t_1, t_2) = F (t_2,t_2^0) x(t_1, t_2^0) + \int\limits_{t_2^0}^{t_2} F(t_2, s) B(s) \sigma_2(s) u(t_1, s)ds \\[5pt]
    y(t_1,t_2) = u(t_1,t_2) - B^*(t_2)~ x(t_1,t_2).
    \end{array} \right.
  \end{equation*}
with inputs and outputs satisfying compatibility conditions (\ref{eq:InCC}), (\ref{eq:OutCC}) (with $\sigma_1=\sigma_{1*}, \sigma_2=\sigma_{2*}$):
\[ \begin{array}{llll}
  \sigma_2(t_2) \frac{\partial}{\partial t_1}u(t_1, t_2) -
  \sigma_1(t_2) \frac{\partial}{\partial t_2}u(t_1,t_2) + \gamma(t_2) u(t_1,t_2) = 0 \\
  \sigma_2(t_2) \frac{\partial}{\partial t_1}y(t_1, t_2) -
  \sigma_1(t_2) \frac{\partial}{\partial t_2}y(t_1,t_2) + \gamma_*(t_2) y(t_1,t_2) = 0
\end{array} \]

\subsection{Gauge equivalence of conservative vessels}
The notion of similarity for vessels turns out to be a notion of gauge equivalence, which generalizes the notion
of unitary equivalence for minimal systems (vessels). It turns out that minimality of a conservative system is
determined at any value of $t_2$. Moreover, observability or controllability alone, which is satisfied for any
value of $t_2$, imposes minimality of the system for all values of $t_2$.

Let $CI\Sigma$ be the i/s/o system (\ref{eq:systemCons}) associated with the vessel
$\mathfrak{CIV}$. Then the following theorem holds:
\begin{thm}
For a vessel $\mathfrak{CIV}$, satisfying regularity assumptions, if there is a uniqueness of the solution for the
output compatibility equation (\ref{eq:OutCC}) then the following conditions are equivalent:
\begin{enumerate}
\item for some $t_2^0$:
\[ \bigvee_{n\geq 0} A_1^n (t_2^0) B (t_2^0) e = H_{t_2^0}.\]
\item for all $t_2$:
\[ \bigvee_{n\geq 0} A_1^n (t_2) B(t_2) e = H_{t_2}.\]
\item Observability (for a fixed $t_2^0$):
\[ \left\{  \begin{array}{lll}
  \frac{\partial f}{dt_1} = A_1(t_2) f \\
  \frac{\partial f}{dt_2} = A_2(t_2) f \\
  v = - B^*(t_2)f
\end{array}
\right. \]
\[
v(t_1,t_2^0) = 0, \forall t_1 \Rightarrow f(0,0)=0
\]
\item Observability (for all $t_1, t_2$):
\[ \left\{  \begin{array}{lll}
  \frac{\partial f}{dt_1} = A_1(t_2)f \\
  \frac{\partial f}{dt_2} = A_2(t_2)f \\
  v = - B^*(t_2)f
\end{array}
\right. \]
\[
v(t_1,t_2) = 0, \forall t_1, t_2 \Rightarrow f(0,0)=0
\]
\end{enumerate}
\end{thm}

\textbf{Proof:} We want to show that all the conditions are equivalent.
Using $F(t_2,t_2^0)$, it is easy to verify that the solution of the system
equations is necessarily of the form
\renewcommand{\theequation}{\arabic{section}.\arabic{equation}}
\begin{equation}
\label{eq:FormvhCons} \begin{array}{lll}
f(t_1, t_2) = e^{t_1 A_1(t_2)} F(t_2,t_2^0) f(0,0) \\
v(t_1, t_2) = - B^*(t_2) e^{t_1 A_1(t_2)} F(t_2,t_2^0) f(0,0).
\end{array}
\end{equation}
Another condition is:
\[\bigvee_{n\geq 0} A_1^{*n}(t_2)B(t_2) =
\bigvee_{n\geq 0} A_1^{n}(t_2) B(t_2), \]
which comes from the colligation condition on $A_1(x)$.
Using these results we shall prove the theorem.

We shall show that $3) \Rightarrow 1) \Rightarrow 2) \Rightarrow 3)$.
\begin{itemize}
\item We want to prove $3) \Rightarrow 1)$. Suppose that the third condition
is satisfied. This means that $\forall e\in E$
\begin{equation}
\label{eq:3condCons}
<B^*(t_2^0) e^{t_1 A_1(t_2^0)} h, e> = 0,~ \forall ~t_1 \Rightarrow h = 0.
\end{equation}
Since $ e^{t_1 A_1(t_2^0)} h$ is clearly an analytic function of $t_1$,
the condition (\ref{eq:3condCons}) is equivalent to $\forall e\in E$
\[ <B^*(t_2^0)  A_1^n(t_2^0) h, e> = <h, A_1^{*n}(t_2^0) B(t_2^0) e> =0
,~ \forall ~n \Rightarrow h = 0, \]
which is true, because the only vector perpendicular to $\bigvee_{n\geq 0}A_1^{*n}(t_2^0)B(t_2^0)e$
is the zero vector and finally condition $(1)$ holds.
\item Here we want to obtain $1)\Rightarrow 2)$. Let us fix $t_2$ and
suppose by contradiction that there exists an $h \neq 0$ such that
$\forall e \in E$ $<h , \bigvee_{n\geq 0}A_1^{*n}(t_2)B(t_2)e> =0$.
Since $F(t_2,t_2^0)$ is invertible, this is equivalent to
$<F(t_2,t_2^0)h' ,\bigvee_{n\geq 0}A_1^{*n}(t_2)B(t_2)e> = 0$. And because
of analyticity of the function $e^{t_1 A_1^*(t_2)}$ the last
norm actually is:
\[ <F(t_2, t_2^0) h', e^{t_1 A_1^*(t_2)} B(t_2) e> = 0, \forall ~ t_1. \]
But it means that $\forall e \in E$
\[ <B^*(t_2) e^{t_1 A_1(t_2)} B(t_2,t_2^0) h, e> = 0, \forall ~t_1. \]
$B^*(t_2) e^{t_1 A_1(t_2)} F(t_2, t_2^0) h$ is the solution of the system equations for
initial condition $B^*(t_2^0) e^{t_1 A_1(t_2^0)}$ for $t_2 = t_2^0$, since there exist
unique solutions for the differential equation, defined by the system equations. We have found
an $h\neq 0$ such that $ <B^*(t_2^0) e^{t_1 A_1(t_2^0)}h, e>  = 0, \forall t_1$.
This means that $0\neq h \perp
\bigvee_{n\geq 0}A_1^{n}(t_2^0)B(t_2^0)e$. Contradiction.
\item Here we want to show that $2) \Rightarrow 3)$. Suppose $2)$ and we
shall show that $3)$ is satisfied too. Suppose that $v(t,t_2^0) = 0, \forall t_1$.
Then for any $t_2$ there exist a unique solution of the system equations with
initial condition $v(t,t_2^0)$. This solution is of the form
\[ v (t_1,t_2) = - B^*(t_2) e^{t_1 A_1(t_2)} F(t_2, 0) f(0, 0). \]
Suppose by contradiction that $h = f(0,0) \neq 0$. Then $\forall t$ and
$\forall e\in E$:
\[ <B(t_2) e^{t_1 A_1(t_2)} F(t_2, 0) h, e>_E = 0, \]
from where $\forall n \in \mathbb{N}$ and $\forall e\in E$:
\[ <B(t_2) A_1^n(t_2) F(t_2, 0) h, e>_E = 0, \]
which is equivalent to $\forall n \in \mathbb{N}$ and $\forall e\in E$:
\[ <F(t_2, 0) h, A_1^{*n}(t_2) B(t_2) e>_H = 0 \]
This means that $F(t_2, 0) h \perp H$ and since $F(t_2, 0)$ is invertible
$h \perp H$. Thus $h = 0$ and $f(0, 0) = 0$.
\item Finally, the conditions $4)$ and $3)$ are equivalent. This is a conclusion from the uniqueness of the
solution for a differential equation. This ends the proof of the theorem.\qed
\end{itemize}

There is also a natural notion of equivalence. Two vessels
\[ \begin{array}{lll}
\mathfrak{IV} = (A_1(t_2), F(t_2,t_2^0), B(t_2); \sigma_1(t_2), \sigma_2(t_2), \gamma(t_2), \gamma_*(t_2);
\mathcal{H}_{t_2},\mathcal{E})\\
\widetilde{\mathfrak{IV}} = (\widetilde A_1(t_2), \widetilde F(t_2,t_2^0), \widetilde B(t_2);
\sigma_1(t_2), \sigma_2(t_2), \gamma(t_2), \gamma_*(t_2); \widetilde{\mathcal H}_{t_2},\mathcal{E})
\end{array} \]
are called \textit{unitary gauge equivalent} (or \textit{unitary kinematically equivalent}),
if there exists $U(t_2): \mathcal{H}_{t_2} \rightarrow \widetilde{\mathcal{H}}_{t_2}$, unitary, such that :
\begin{equation}
\label{eq:CUconnectI} \left\{\begin{array}{ll}
\widetilde A_1(t_2)       & = U(t_2) A_{1}(t_2) U^{-1}(t_2) \\
\widetilde{F}(t_2,t_2^0) & = U(t_2) F(t_2,t_2^0) U^{-1}(t_2^0) \\
\widetilde{B}(t_2)        & = U (t_2) B(t_2)
\end{array}\right.
\end{equation}
\subsection{\label{sec:CDiffVessel}Differential form of conservative vessels}
As we have seen, one of the advantages of the $t_1$-invariant vessels is that there is always a more convenient way to work with
them. It turns out that conservative vessels are always gauge equivalent and obtain systems equations in the
form of differentiation. Of course, on first sight it seems to be impossible, because there exists a
continuum of Hilbert spaces (for each $t_2$) and the differential of the evolution semigroup is not generally
defined. In our case we can bring all the spaces to one with an isometric map. In order to do it, one needs
the following proposition
\begin{prop} There is an invertible, self-adjoint absolutely continuous operator
$\Psi(t_2) : \mathcal H_{t_2^0} \longrightarrow \mathcal H_{t_2^0}$ such that
\[ F^*(t_2,t_2^0) F(t_2,t_2^0) = \Psi(t_2)^2 \]
\end{prop}
\textbf{Proof:} First, we notice that the spectrum of $F^*(t_2,t_2^0)F(t_2,t_2^0)$
is separated from zero, since this operator is invertible, and is bounded.
Thus we can define the operator $\Psi(t_2)$ as the squareroot of this operator, using the
Riesz-Dunford calculus (where the root is taken to have the principal value for
each $\lambda \in \mathbb{C}$):
\[
\Psi(t_2) = \Psi^*(t_2) = \oint \sqrt{\lambda} (\lambda I-
F^*(t_2,t_2^0) F(t_2,t_2^0))^{-1} d\lambda
\]
Since we can take the path of integration symmetric with respect to the real axis, $\Psi(t_2)$ will
be obtained \textit{self-adjoint}. Moreover, the spectrum of this operator according to the spectral
theorem will be the square root of the spectrum of $F^*(t_2,t_2^0) F(t_2,t_2^0)$ and thus is
positive. Finally, since we have obtained self-adjoint operator $\Psi(t_2)$ with positive spectrum it would be
a \textit{positive operator}. We notice also that
\[
F^*(t_2,t_2^0) F(t_2,t_2^0) = \oint \lambda (\lambda I-
F^*(t_2,t_2^0) F(t_2,t_2^0))^{-1} d\lambda.
\]
Since $F(t_2,t_2^0)$ is an invertible operator, so will be $F^*(t_2,t_2^0)$ and consequently, $\Psi(t_2)$.

From the second colligation condition
\[ F^*(t_2,t_2^0)F(t_2,t_2^0) = I - \int_{t_2^0}^{t_2}
   F^*(t_2,t_2^0)B(t_2) \sigma_2(t_2) B^*(t_2)F(t_2,t_2^0)
\]
we obtain that
\[ \Psi(t_2) \Psi(t_2) = I - \int_{t_2^0}^{t_2}
   F^*(t_2,t_2^0)B(t_2) \sigma_2(t_2) B^*(t_2)F(t_2,t_2^0).
\]
Thus $\Psi(t_2)^2$ is an absolutely continuous function of $t_2$. In order to obtain that
$\Psi(t_2)$ has the same property, we prove the following
\begin{lemma} \label{lem:SqrtAbsCont}
 Suppose that $R(t_2)$ is an absolutely continuous operator on the interval $[a,b]$. Suppose also that it is self-adjoint and invertible. Then its square root, defined by
\[ \sqrt{R}(t_2) = \oint \sqrt{\lambda} (\lambda I - R(t_2^0))^{-1} d\lambda \]
will be an absolutely continuous operator too.
\end{lemma}
\textbf{Proof:} We shall prove it using the definition of absolute continunity. Let us
fix $\epsilon > 0$. Since $R(t_2)$ is absolutely continuous there exists $\delta > 0$ such that
for every sequence
\[ a \leq a_1 < b_1 < a_2 < b_2 < \ldots < a_k < b_k \leq b \]
so that $\sum\limits_{i=1}^{k} (b_i-a_i) < \delta_1$, it holds that
$\sum\limits_{i=1}^{k} \| R(b_i)- R(a_i) \| < \epsilon$. We shall find $\delta$ such that the same will
be correct for $\sqrt{R}(t_2)$.
\[ \begin{array}{lllll}
\sqrt{R}(b_i)- \sqrt{R}(a_i) & = \oint \sqrt{\lambda} (\lambda I -
R(b_i))^{-1} d\lambda - \oint \sqrt{\lambda} (\lambda I -
R(a_i))^{-1} d\lambda = \\
& = \oint \sqrt{\lambda} [(\lambda I - R(b_i))^{-1} - (\lambda I -
R(a_i))^{-1}] d\lambda = \\
& = \oint \sqrt{\lambda} (\lambda I -
R(b_i))^{-1}(R(a_i) - R(b_i))(\lambda I - R(a_i))^{-1} d\lambda.
\end{array} \]
Now we can evaluate the norm of this operator:
\[  \begin{array}{lllll}
\| \sqrt{R}(b_i)- \sqrt{R}(a_i) \| & =
\| \oint \sqrt{\lambda} (\lambda I -
R(b_i))^{-1}(R(a_i) - R(b_i))(\lambda I - R(a_i))^{-1} d\lambda \| \\
& \leq \oint \| \sqrt{\lambda} (\lambda I -
R(b_i))^{-1}(R(a_i) - R(b_i))(\lambda I - R(a_i))^{-1} \| d\lambda  \\
& \leq \oint \| \sqrt{\lambda} (\lambda I -
R(b_i))^{-1}\| \|(R(a_i) - R(b_i))\| \|(\lambda I - R(a_i))^{-1} \| d\lambda \\
& \leq \oint \| \sqrt{\lambda} (\lambda I - R(b_i))^{-1} \|
\|(\lambda I - R(a_i))^{-1} \| d\lambda \| R(a_i) - R(b_i) \|
\end{array} \]
So, if we succeed to prove that
\[ \oint \| \sqrt{\lambda} (\lambda I - R(b_i))^{-1} \|
\|(\lambda I - R(a_i))^{-1} \| d\lambda \leq K \in \mathbb{R}\]
we shall obtain the absolute continuity of $\sqrt{R}(t_2)$, by taking
$\delta=\frac{\delta_1}{K}$.
For this it is enough to show that
\[ \sup\limits_{s \in [a,b]} \| (\lambda I - R(s))^{-1} \| < \infty, \]
because then
\[ \begin{array}{lllll}
 \oint \| \sqrt{\lambda} (\lambda I - R(b_i))^{-1} \|
\|(\lambda I - R(a_i))^{-1} \| d\lambda \leq \\
\big( \sup\limits_{s \in [a,b]} \| (\lambda I - R(s))^{-1} \|\big)^2 L
\max\limits_{\lambda \in O} |\sqrt{\lambda}|,
\end{array} \]
where $L$ is the length of the path $O$ of the integral. As for
$\big( \sup\limits_{s \in [a,b]} \| (\lambda I - R(s))^{-1} \|\big)^2$, we just notice that
it is continuous in $\lambda$ ans $s$ function and since it is considered on a compact space
\[ (\lambda, s) \in \mathcal O \times [a,b], \]
it attains its finite maximum. \qed

Now we are ready to build a differential form of the vessel. For this let us define (by denoting here $\mathcal H = \mathcal H_{t_2^0}$)
a unitary operator $U(t_2): \mathcal H_{t_2} \longrightarrow \mathcal H$ by
\[ U(t_2) = \Psi(t_2) F(t_2^0, t_2). \]
Then we define new operators on the space $\mathcal H_{t_2^0}$:
\[ \begin{array}{lllll}
\breve F(t_2, t_2') & = U(t_2) F(t_2,t_2') U(t_2')^{-1} = \Psi(t_2) \Psi(t_2')^{-1} \\
\breve A_1(t_2)         & = U(t_2) A_1(t_2) U^{-1}(t_2) = \Psi(t_2) A_1(t_2^0) \Psi(t_2)^{-1} \\
\breve A_2(t_2)         & = \frac{d}{dt_2} \Psi(t_2) \Psi^{-1}(t_2) \\
\widetilde{\breve B}(t_2)           & = U(t_2) \widetilde B(t_2) = \Psi(t_2) F(t_2^0,t_2) \widetilde B(t_2)
\end{array} \]
Thus we obtain a \textbf{vessel in the differential form}:
\[ \mathfrak{DV} = (\breve A_1(t_2), \breve A_2(t_2),\widetilde{\breve B}(t_2);
        \sigma_1(t_2), \sigma_2(t_2), \gamma(t_2), \gamma_*(t_2);
        \mathcal H, \mathcal{E})
\]
which satisfies the following axioms:
\[ \begin{array}{lllllll}
    \widetilde A_j(t_2) + \widetilde A^*_j(t_2) + \widetilde{B}(t_2) \sigma_j(t_2) \widetilde{B}^*(t_2) = 0, ~~~~j=1,2\\
    \frac{d}{dt_2} \breve A_1(t_2) = \breve A_2(t_2) \breve A_1(t_2) - \breve A_1(t_2) \breve A_2(t_2) \\
    \frac{d}{dt_2} \big(\widetilde{\breve B}(t_2) \sigma_1(t_2)\big) - \breve A_2(t_2) \widetilde{\breve B}(t_2) \sigma_1(t_2) + \breve A_1(t_2)
            \widetilde{\breve B}(t_2) \sigma_2(t_2) + \widetilde{\breve B}(t_2) \gamma(t_2) = 0 \\
        \frac{d}{dt_2} \big( \sigma_1(t_2) \breve B^*(t_2) \big) + \sigma_1(t_2) \breve B^*(t_2) \breve A_2(t_2) +
                \sigma_2(t_2) \breve B^*(t_2)\breve A_1(t_2) + \gamma_*(t_2) \breve B^*(t_2) = 0 \\
  \gamma = \sigma_2 \breve B^* \widetilde{\breve B} \sigma_1  - \sigma_1 \breve B^* \widetilde{\breve B} \sigma_2  + \gamma_*.
\end{array} \]
The differential vessel is associated with the system
\[
    D\Sigma: \left\{ \begin{array}{lll}
    \frac{\partial}{\partial t_1}x(t_1,t_2) = \breve A_1(t_2) x(t_1,t_2) + \widetilde{\breve B}(t_2) ~\sigma_1(t_2) ~u(t_1,t_2) \\
    \frac{\partial}{\partial t_2}x(t_1,t_2) = \breve A_2(t_2) x(t_1,t_2) + \widetilde{\breve B}(t_2) ~\sigma_2(t_2) ~u(t_1,t_2) \\
    y(t_1,t_2) = u(t_1,t_2) + \breve B^*(t_2) x(t_1,t_2)
    \end{array} \right.
\]
and compatibility conditions for the input/output signals:
\[ \begin{array}{lll}
\sigma_2(t_2) \frac{\partial}{\partial t_1}u(t_1,t_2) - \sigma_1(t_2) \frac{\partial}{\partial t_2}u(t_1,t_2) +
\gamma(t_2) u(t_1,t_2) = 0, \\
\sigma_2(t_2) \frac{\partial}{\partial t_1}y(t_1,t_2) - \sigma_1(t_2) \frac{\partial}{\partial t_2}y(t_1,t_2) +
\gamma_*(t_2) y(t_1,t_2) = 0.
\end{array} \]

\textit{Conversely}, starting from a vessel in the differential form
\[ \mathfrak{V} = (\breve A_1(t_2), \breve A_2(t_2),\widetilde{\breve B}(t_2);
        \sigma_1(t_2), \sigma_2(t_2), \gamma(t_2), \gamma_*(t_2);
        \mathcal H, \mathcal{E})
\]
we can build a vessel in the integrated form
\[ \mathfrak{IV} = (A_1(t_2), F(t_2,t_2^0),\widetilde B(t_2); \sigma_1(t_2), \sigma_2(t_2), \gamma(t_2),
\gamma_*(t_2); \mathcal H_{t_2}, \mathcal{E})
\]
using the following definitions:
\[ \begin{array}{llll}
A_1(t_2)                    & = \breve A_1(t_2) \\
F(t_2,t_2^0)            & = \text{ the evolution semigroup, generated by the operator $\breve A_2(t_2)$} \\
\widetilde B(t_2) & = \widetilde{\breve B}(t_2) \\
\mathcal{H}_{t_2} & = \mathcal H \text{ - the same for all}.
\end{array} \]

In order to prove the theorem of gauge equivalence for
conservative vessels it is more convenient and readable to present
it in the differential form. This theorem is an analogue in our
framework of the standard unitary equivalence theorem for minimal
(or closely connected) conservative systems \cite{bib:BallCohen}.
\begin{mthm}\label{mthm:CGE}
Assume that we are given two minimal $t_1$-invariant vessels  $\mathfrak{V}, \widetilde{\mathfrak{V}}$ in the
differential form with transfer functions $S(\lambda,t_2)$, $\widetilde{S}(\lambda,t_2)$.
Then the vessels are gauge equivalent iff
$S(\lambda,t_2) = \widetilde{S}(\lambda,t_2) $ for all points of analyticity.
\end{mthm}
\textbf{Proof:} Notice first that from the minimality of the vessels
\[ \begin{array}{lll}
\bigvee\limits_{n\geq 0}A_{1}^{n}(t_2)B(t_2)e = H \\
\bigvee\limits_{n\geq 0}\tilde A_{1}^{n}(t_2)\tilde{B}(t_2)e = \tilde{H}
\end{array} \]
then define an isometry $U(t_2)$ on the dense sets
\begin{equation}
\label{eq:defU} U(t_2) A_1^n(t_2) B(t_2) = \tilde A_1^n(t_2) \tilde B(t_2).
\end{equation}
Define next a new operator
\begin{defn}
\label{def:dU}
\[ \frac{d U(t_2)}{dt_2} \big[ A_{1}^{n}(t_2)B(t_2) \big] =
\frac{d}{dt_2} \big[ \tilde A_{1}^{n}(t_2)\tilde B(t_2) \big] -
U(t_2) \frac{d}{dt_2} \big[ A_{1}^{n}(t_2) B(t_2) \big]. \]
\end{defn}
Then
\begin{lemma} $\frac{d U(t_2)}{dt_2}$, defined by the definition (\ref{def:dU}) is
\begin{enumerate}
\item Well defined.
\item A bounded operator from $H$ to $\tilde{H}$.
\item The derivative of $U(t_2)$ in uniform operator topology.
\end{enumerate}
\end{lemma}
\textbf{Proof:} We shall prove the lemma for each claim:
\begin{enumerate}
\item Here it is enough to show that
\[
\sum_i c_i A_1^{n_i}(t_2) B(t_2) e = 0 \Rightarrow
\frac{d U(t_2)}{dt_2} \sum_i c_i A_1^{n_i}(t_2) B(t_2) e = 0 \]
But this is a direct result of the definition
\[ \begin{array}{llll}
\frac{d U(t_2)}{dt_2} [ \sum_i c_i A_1^{n_i}(t_2) B(t_2) e ] & =
\frac{d}{dt_2} [\sum_i c_i \tilde A_1^{n_i}(t_2) \tilde B(t_2) e] -
U(t_2) \frac{d}{dt_2}[\sum_i c_i A_1^{n_i}(t_2) B(t_2) e ] = \\
& = \frac{d}{dt_2} [U(t_2) \sum_i c_i A_1^{n_i}(t_2) B(t_2) e] -
\frac{d}{dt_2}[\sum_i c_i A_1^{n_i}(t_2) B(t_2) e] \\
& = \frac{d}{dt_2} [U(t_2) 0] - U(t_2) \frac{d}{dt_2} 0 = 0.
\end{array} \]
\item A general claim that an operator defined in this way is
bounded, fails. But in the case of vessels it is true.
\begin{equation} \label{eq:dUchain}
\begin{array}{lll}
\frac{d U(t_2)}{dt_2} \big[ A_{1}^{n}(t_2)B(t_2) \big]
&  = \frac{d}{dt_2} \big[ \tilde A_{1}^{n}(t_2)\tilde B(t_2) \big] -
U(t_2) \frac{d}{dt_2} \big[ A_{1}^{n}(t_2)B(t_2) \big] \\
& = \frac{d \tilde A_{1}^{n}(t_2)}{dt_2} \tilde B(t_2) +
\tilde A_{1}^{n}(t_2) \frac{d \tilde B(t_2)}{dt_2} - U(t_2) \frac{d A_{1}^{n}(t_2)}{dt_2}B(t_2) -
U(t_2)A_{1}^{n}(t_2)\frac{d B(t_2)}{dt_2}.
\end{array} \end{equation}
Here, for a better illustration we prefer to evaluate $\frac{d B(t_2)}{dt_2}$
separately. Namely, from the vessel condition:
\begin{equation} \label{eq:phiStar}
\frac{d B(t_2)}{dt_2} = A_2(t_2) B(t_2) -
(A_1(t_2) B(t_2) \sigma_2(t_2) + B(t_2) \gamma^*(t_2)) \sigma_1^{-1}.
\end{equation}
We also have the Lax equation for the derivative of $A_1(t_2)$. From which it
is easy to obtain that:
\begin{equation} \label{eq:dA1n}
\frac{d A_1^n(t_2)}{dt_2} = A_1(t_2)^n A_2(t_2) - A_2(t_2) A_1^n(t_2).
\end{equation}
After inserting the formulas (\ref{eq:phiStar}) and (\ref{eq:dA1n}) into
(\ref{eq:dUchain}) we shall obtain:
\[ \begin{array}{lllllll}
\frac{d U(t_2)}{dt_2} \big[ A_{1}^{n}(t_2)B(t_2) \big]
& = [\tilde A_1(t_2)^n \tilde A_2(t_2) - \tilde A_2(t_2) \tilde A_1^n(t_2)] \tilde B(t_2) + \\
& + \tilde A_1^n(t_2) [A_2(t_2) B(t_2) -
(A_1(t_2) B(t_2) \sigma_2(t_2) + B(t_2) \gamma^*(t_2)) \sigma_1^{-1}] \\
& - U(t_2) [A_1^n(t_2) A_2(t_2) - A_2(t_2) A_1^n(t_2)] B(t_2) \\
& - U(t_2) A_1^n(t_2) [ A_2(t_2) B(t_2) -
(A_1(t_2) B(t_2) \sigma_2(t_2) + B(t_2) \gamma^{in*}(t_2)) \sigma_1^{-1}] \\
& = (\tilde A_2(t_2) - \tilde A_1(t_2))\tilde A_1^n(t_2)\tilde B(t_2) +
\tilde A_1^n(t_2)\tilde B(t_2) \gamma(t_2)^* \sigma_1^{-1} \\
& - U(t_2) (A_2(t_2) - A_1(t_2))  A_1^n(t_2) B(t_2)
      - U(t_2) A_1^n(t_2) B(t_2) \gamma(t_2)^* \sigma_1^{-1}.
\end{array} \]
From this equality, we may see that $\frac{d U(t_2)}{dt_2}$, applied to
$A_{1}^{n}(t_2)B(t_2)$ is actually combined
from a couple of operators. They are:
\begin{enumerate}
\item $\tilde A_2(t_2) - \tilde A_1(t_2)$, bounded and applied to
$\tilde A_{1}^{n}(t_2)\tilde B(t_2)$, which has the same norm as
$A_{1}^{n}(t_2)B(t_2)$.
\item $A_2(t_2) - A_1(t_2)$, $U(t_2)$ bounded.
\item The operator $\Gamma$ of the following form
\[ \Gamma[A_1^n(t_2) B(t_2)]
= A_1^n(t_2) B(t_2) \gamma(t_2)^*\sigma_1^{-1*}\]
\end{enumerate}
We claim that this is a bounded operator too. Suppose that
$ \gamma(t_2)^*\sigma_1^{-1*} = [G_{ki}(t_2)]$ - a continuous
matrix-function of $t_2$ from $E$ to $E$ in the standart basis $\{ e_i\}$.
Then $\forall e = \sum a_i e_i \in E$
\[ \begin{array}{lll}
A_1^n(t_2) B(t_2) \gamma(t_2)^*\sigma_1^{-1*} e =
 A_1^n(t_2) B(t_2) \sum_{i=1}^n \sum_{k=1}^n G_{ki}(t_2) a_k e_i \\
= \sum_{i,k=1}^n G_{ki}(t_2) a_k A_1^n(t_2) B(t_2) e_i,
\end{array} \]
which is sum of scalar continuous functions and thus a sum of bounded
operators. Thus we have obtained that $\frac{d U(t_2)}{dt_2}$ is a bounded
operator.
\item Let us try to evaluate the derivative of $U(t_2)$ in the strong
topology by the definition. Let us denote $h(t_2) = A_1^n(t_2) B(t_2)$ and
$\tilde{h}(t_2) = \tilde A_1^n(t_2) \tilde B(t_2)$. Then, using definition
\ref{def:dU}, we shall obtain:
\[ \begin{array}{lllclll}
\lim\limits_{\triangle t_2 \rightarrow 0}
\frac{1}{\triangle t_2}(U(t_2 + \triangle t_2) - U(t_2)) h(t_2) = \\
\lim\limits_{\triangle t_2 \rightarrow 0}
\frac{1}{\triangle t_2}(U(t_2 + \triangle t_2) h(t_2) +
                      U(t_2 + \triangle t_2) h(t_2 + \triangle t_2) -
                      U(t_2 + \triangle t_2) h(t_2 + \triangle t_2) -
                      U(t_2) h(t_2)) = \\
\lim\limits_{\triangle t_2 \rightarrow 0}
\frac{1}{\triangle t_2}(U(t_2 + \triangle t_2) h(t_2 + \triangle t_2) - U(t_2) h(t_2)) \\
- \lim\limits_{\triangle t_2 \rightarrow 0}
\frac{1}{\triangle t_2}(U(t_2 + \triangle t_2) h(t_2 + \triangle t_2) -
                      U(t_2 + \triangle t_2) h(t_2)) \stackrel{def}{=} \\
\frac{d}{dt_2} [U(t_2) h(t_2)] - \lim\limits_{\triangle t_2 \rightarrow 0}
(U(t_2 + \triangle t_2) \frac{1}{\triangle t_2}[h(t_2 + \triangle t_2) - h(t_2)]) \\
\end{array} \]
Now let us show that
$\lim\limits_{\triangle t_2 \rightarrow 0} U(t_2 + \triangle t_2) = U(t_2)$. From
the definition of $U(t_2)$ (\ref{eq:defU}) we may conclude that $U(t_2)$ is defined as an
isomorphism of continuous functions (from $E$ to $H$). We shall show that this
limit is correct on the dense set $A_1^{n}(t_2)B(t_2) e$, $e \in E$:
\[ \begin{array}{llllll}
\lim_{\triangle t_2 \rightarrow 0} [
U(t_2 + \triangle t_2) A_1^{n}(t_2)B(t_2) - U(t_2)A_1^{n}(t_2)B(t_2)] = \\
\lim_{\triangle t_2 \rightarrow 0} [
U(t_2 + \triangle t_2) A_{1}^{n}(t_2+\triangle t_2)
B(t_2+\triangle t_2) - U(t_2)A_1^{n}(t_2)B(t_2) + \\
+ U(t_2 + \triangle t_2) A_1^{n}(t_2)B(t_2) - U(t_2 + \triangle t_2)
A_1^{n}(t_2+\triangle t_2) B(t_2+\triangle t_2)] = \\
\lim_{\triangle t_2 \rightarrow 0}
\tilde A_{1}^{n}(t_2+\triangle t_2)\tilde B(t_2+\triangle t_2) -
\tilde A_{1}^{n}(t_2)\tilde B(t_2) + \\
+ \lim_{\triangle t_2 \rightarrow 0} [
U(t_2 + \triangle t_2) (A_1^{n}(t_2)B(t_2) -
A_1^{n}(t_2+\triangle t_2) B(t_2+\triangle t_2)) = 0,
\end{array} \]
since $U(t_2 + \triangle t_2)$ are uniformly bounded by $1$.
This is a continuity of $U(t_2)$ in $t_2$ on the dense set
$A_1^{n}(t_2)B(t_2)$. Finally we obtain
\[ \begin{array}{lllll}
\lim\limits_{\triangle t_2 \rightarrow 0}
\frac{1}{\triangle t_2}(U(t_2 + \triangle t_2) - U(t_2)) h(t_2) =
\frac{d}{dt_2} \tilde{h}(t_2) - U(t_2) \frac{d}{dt_2} h(t_2)
= \frac{d U(t_2)}{dt_2} h(t_2) ~~~~~~~~ \qed
\end{array} \]
\end{enumerate}
Now we are ready to prove the theorem. \newline
\noindent\textit{Proof of the main theorem \ref{mthm:CGE}:} we shall show that necessarily
\begin{equation}
\label{A2 connection}
\tilde A_2(t_2) = U(t_2) A_2(t_2) U^{-1}(t_2) + \frac{dU(t_2)}{dt_2} U^{-1}(t_2)
\end{equation}
establishing the gauge equivalence.

In the sequel, we permit ourself to omit assigning $t_2$-dependence, since its
dependence will be easily understood from the context.

First of all, using the Lax equation
\begin{equation}
\label{eq:Adot}
\dot A_1 = A_1 A_2 - A_2 A_1
\end{equation}
we find that
$A_1 A_2 =  \dot A_1 + A_2 A_1$.

Rearranging one of the vessel conditions, we obtain that:
\[
A_2 B\sigma_1=B\gamma + A_1 B\sigma_2 + \frac{d B}{dt_2} \sigma_1.
\]
Since $\sigma_1$ is invertible, we obtain:
$A_2 B = B \gamma \sigma_1^{-1}+ A_1B\sigma_2\sigma_1^{-1} + \frac{d B}{dt_2}$.
Multiplying the equation by $A_1$ and using the formula (\ref{eq:Adot}) we obtain:
\[ A_2 A_1B = - \dot A_1 B A_1B\gamma\sigma_1^{-1} + A_1^2 B \sigma_2 \sigma_1^{-1}
+ A_1\frac{d B}{dt_2}.
\]
Continuing to multiply by $A_1$ from the left and using the same formula for $\dot A_1$ twice, we obtain:
\[
A_2 A_1^2 B =- A_1 \dot A_1B - \dot A_1 A_1 B + A_1^2 B \gamma\sigma_1^{-1}
+ A_1^3 B \sigma_2 \sigma_1^{-1} + A_1^2 \frac{d B}{dt_2}.
\]

In the general case, using induction it is very easy to show that:
\begin{equation}
\label{eq:GenA2}
\begin{array}{llll}
A_2A_1^{n}B &
= - \dot A_1\underbrace {A_1\cdots A_1}_{\mbox{n-1}}B
- A_1\dot A_1\underbrace {A_1\cdots A_1}_{\mbox{n-2}}B
- \cdots- \underbrace {A_1\cdots A_1}_{\mbox{n-1}}\dot A_1B \\
& + A_1^nB\gamma\sigma_1^{-1}
+ A_1^{n+1}B\sigma_2\sigma_1^{-1}
+ A_1^n\frac{d B}{dt_2}.
\end{array}
\end{equation}
Let us define the first summation as a formula:
\begin{defn}
\label{def:BracketK}
$(A_1^n)_k = A_1^k \dot{A_1} A_1^{n-k-1}$, for $k = 0, 1,\ldots, n-1$.
\end{defn}
Then it easy to verify that
\[ \frac{d}{dt_2} A_1^n =
\dot A_1\underbrace {A_1\cdots A_1}_{\mbox{n-1}}
+ A_1\dot A_1\underbrace {A_1\cdots A_1}_{\mbox{n-2}}
+ \cdots
+ \underbrace {A_1\cdots A_1}_{\mbox{n-1}}\dot A_1 =
 \sum\limits_{k=0}^{n-1} (A_1^n)_k . \]
Thus the formula (\ref{eq:GenA2}) converts to
\[
A_2A_1^{n}B = - \sum\limits_{k=0}^{n-1} (A_1^n)_k B
+ A_1^{n+1} B \gamma \sigma_1^{-1}+ A_1^{n+1}B\sigma_2\sigma_1^{-1} + A_1^{n}\frac{d B}{dt_2}.
\]

A similar formula holds for $\tilde A_2$:
\[
\tilde A_2\tilde A_1^{n}B = -  \sum\limits_{k=0}^{n-1} (\tilde A_1^n)_k \tilde B
+ \tilde A_1^{n+1} \tilde B \gamma \sigma_1^{-1} + \tilde A_1^{n+1} \tilde B \sigma_2
  \sigma_1^{-1} + \tilde A_1^n \frac{d \tilde B}{dt_2},
\]
which can be rewritten using the action of $U$:
\begin{equation} \label{eq:A2Full} \begin{array}{llll}
\tilde A_2 U A_1^n B
& = - \sum\limits_{k=0}^{n-1}(U A_1^n U^{-1})_k U B
    + U A_1^n B \gamma \sigma_1^{-1}
    + U A_1^{n+1} B \sigma_2\sigma_1^{-1} \\
&   + U A_1^n \frac{d}{dt_2} BU.
\end{array} \end{equation}
We next simplify the first term of this formula:
\begin{lemma}
\label{lem:backForm}
The following holds:
\[ \sum\limits_{k=0}^{n-1}(U A_1^n U^{-1})_k U =
\frac{d U}{dt_2} A_1^n + U \sum\limits_{k=0}^{n-1}(A_1^n)_k
+ U A_1^n \frac{d U^{-1}}{dt_2} U. \]
\end{lemma}
\textbf{Proof:} We shall prove this formula, using definition
\ref{def:BracketK}:
\begin{equation} \label{eq:VerExpanded}\begin{array}{lll}
\sum\limits_{k=0}^{n-1}(U A_1^n U^{-1})_k U
& = \frac{d}{dt_2}(U A_1 U^{-1}) U \underbrace{A_1 \ldots A_1}_{n-1} U^{-1} \\
 & + U A_1 U^{-1} \frac{d}{dt_2}(U A_1 U^{-1}) U
        \underbrace{A_1 \ldots A_1}_{n-2} U^{-1} + \ldots \\
 & + U \underbrace{A_1 \ldots A_1}_{n-1} U^{-1} \frac{d}{dt_2}(U A_1 U^{-1})
\end{array} \end{equation}
The common factor for all the members in the sum is
$\frac{d}{dt_2}(U A_1 U^{-1})$, which may be evaluated as follows:
\[ \frac{d}{dt_2}(U A_1 U^{-1}) = \frac{d U}{dt_2} A_1 U +
U \frac{d A_1}{dt_2} U^{-1} + U A_1 \frac{d U^{-1}}{dt_2}. \]
Inserting this expression into the formula (\ref{eq:VerExpanded}) we
obtain terms of the following forms:
\begin{enumerate}
\item $U A_1^k U^{-1} \frac{d U}{dt_2} A_1^{n-k} U^{-1}$,
\item $U A_1^k U^{-1} U \frac{d A_1}{dt_2} U^{-1} A_1^{n-k}U^{-1} =
U A_1^k \frac{d A_1}{dt_2} A_1^{n-k} U^{-1},$
\item $U A_1^k \frac{d U^{-1}}{dt_2} U A_1^{n-k} U^{-1}.$
\end{enumerate}
But the sum of the first and third terms vanish, because:
\[ \begin{array}{ll}
U A_1^k U^{-1} \frac{d U}{dt_2} A_1^{n-k} U^{-1} +
U A_1^k \frac{d U^{-1}}{dt_2} U A_1^{n-k} U^{-1}
&=U A_1^k (U^{-1} \frac{d U}{dt_2} +  \frac{d U^{-1}}{dt_2} U) A_1^{n-k} U^{-1}\\
& = U A_1^k \frac{d I}{dt_2} A_1^{n-k} U^{-1} \\
& = 0.
\end{array} \]
And only the first $(k=0)$ and the last elements $(k=n)$ of the sum
remain. Thus we obtain the following formula:
\[
\sum\limits_{k=0}^{n-1}(U A_1^n U^{-1})_k U =
\frac{d U}{dt_2} A_1^n + U \sum\limits_{k=0}^{n-1}(A_1^n)_k
+ U A_1^n \frac{d U^{-1}}{dt_2} U.
\]
\qed

We may proceed with calculating $\tilde A_2$ on vectors of the form
$A_1^n B$ (spanning a dense subspace of H). From formula
(\ref{eq:A2Full}), using lemma \ref{lem:backForm} and expanding the formula
for derivative, we obtain
\[ \begin{array}{llllll}
\tilde A_2 U A_1^n B
& = - \sum\limits_{k=0}^{n-1}(U A_1^n U^{-1})_k U B + U A_1^n B \gamma \sigma_1^{-1} +
   U A_1^{n+1} B \sigma_2\sigma_1^{-1} +  U A_1^n \frac{d}{dt_2} B U \\
&   = - \frac{d U}{dt_2} A_1^n B - U \sum\limits_{k=0}^{n-1}(A_1^n)_k B - U A_1^n \frac{d U^{-1}}{dt_2} U B \\
&   + U A_1^n B \gamma \sigma_1^{-1} + U A_1^{n+1} B \sigma_2\sigma_1^{-1} + U A_1^n U^{-1} U \frac{d B}{dt_2}
    + U A_1^n U^{-1} \frac{d U}{dt_2} B \\
&   = - \frac{d U}{dt_2} A_1^n B + U A_2 A_1^n B + U A_1^n \frac{d U^{-1}}{dt_2}U B+U A_1^n U^{-1} \frac{d U}{dt_2} B \\
&   = - \frac{d U}{dt_2} A_1^n B + U A_2 A_1^n B.
\end{array}\]
Since we have been working on a dense set we may omit $A_1^n B$
to obtain:
\[ \begin{array}{ll}
\tilde A_2 U = - \frac{d U}{dt_2} + U A_2 \\
\tilde A_2 = - \frac{d U}{dt_2} U^{-1} + U A_2 U^{-1}.
\end{array} \]
And finally we conclude that on the dense set:
\[ \left\{\begin{array}{ll}
\widetilde A_1(t_2)    & = U(t_2) A_1(t_2) U^{-1}(t_2) \\
\widetilde A_2         & = U(t_2) A_2(t_2) U^{-1}(t_2) + \frac{d U(t_2)}{dt_2} U^{-1}(t_2) \\
\widetilde B(t_2)      & = U(t_2) B(t_2)
\end{array}\right. \]
which is exactly the gauge equivalence of vessels. \qed

\subsection{Transfer function of a conservative vessel}
Following the same lines as for a non conservative vessel, we perform separation of variables. Taking all the trajectory data in the form
\[ \begin{array}{lll}
u(t_1,t_2) = u_\lambda(t_2) e^{\lambda t_1}, \\
x(t_1,t_2) = x_\lambda(t_2) e^{\lambda t_1}, \\
y(t_1,t_2) = y_\lambda(t_2) e^{\lambda t_1},
\end{array} \]
we arrive at the notion of a transfer function. Note that as before $u(t_1,t_2), y(t_1,t_2)$ satisfy PDEs, but
$u_\lambda(t_2), y_\lambda(t_2)$ are solutions of ODEs with a spectral parameter $\lambda$,
\[ \begin{array}{lll}
\lambda \sigma_2(t_2)  u_\lambda(t_2) - \sigma_1(t_2) \frac{\partial}{\partial t_2}u_\lambda(t_2) +
\gamma(t_2) u_\lambda(t_2) = 0, \\
\lambda \sigma_2(t_2) y_\lambda(t_2) - \sigma_1(t_2) \frac{\partial}{\partial t_2}y_\lambda(t_2) + \gamma_*(t_2)
y_\lambda(t_2) = 0.
\end{array} \]
The corresponding i/s/o system becomes
\[ \left\{ \begin{array}{lll}
    x_\lambda(t_2) = (\lambda I - A_1(t_2))^{-1} B(t_2) \sigma_1(t_2) u_\lambda(t_2) \\
    x(t_1, t_2) = F (t_2,\tau_2) x(t_1, \tau_2) + \int\limits_{\tau_2}^{t_2} F(t_2, s) B(s) \sigma_2(s) u(t_1, s)ds \\
    y_\lambda(t_2) = u_\lambda(t_2) - B^*(t_2) x_\lambda(t_2).
    \end{array} \right.
\]
The output $y_\lambda(t_2) = u_\lambda(t_2) - B^*(t_2) x_\lambda(t_2)$ may be found from the first i/s/o equation:
\[ y_\lambda(t_2) = S(\lambda, t_2) u_\lambda(t_2), \]
using the \textit{transfer function}
\[ S(\lambda, t_2) = I - B^*(t_2) (\lambda I - A_1(t_2))^{-1} B(t_2) \sigma_1(t_2). \]
Here $\lambda$ is outside the spectrum of $A_1(t_2)$, which is independent of $t_2$ by (\ref{eq:LaxCond}). We
emphasize here that $S(\lambda,t_2)$ is a function of $t_2$ for each $\lambda$ (which is a frequency
variable corresponding to $t_1$).

\begin{prop} \label{prop:Sdef}
$S(\lambda,t_2) = I - B^*(t_2) (\lambda I - A_1(t_2))^{-1} B(t_2) \sigma_1(t_2)$ has the following properties:
    \begin{enumerate}
        \item For almost all $t_2$, $S(\lambda, t_2)$ is an analytic function of $\lambda$
          in the neighborhood of $\infty$, where it satisfies:
        \[ S(\infty, t_2) = I_{n \times n} \].
      \item For all $\lambda$, $S(\lambda, t_2)$ is an absolutely continuous function of $t_2$.
        \item The following inequalities are satisfied:
          \[ \begin{array}{llll}
            S(\lambda, t_2)^* \sigma_1(t_2) S(\lambda, t_2) = \sigma_1(t_2), & \Re{\lambda} = 0 \\
            S(\lambda, t_2)^* \sigma_1(t_2) S(\lambda, t_2) \geq \sigma_1(t_2), & \Re{\lambda} \geq 0 \\
          \end{array} \]
          for $\lambda$ in the domain of analyticity of $S(\lambda,t_2)$.
        \item For each fixed $\lambda$, multiplication by $S(\lambda, t_2)$ maps
          solutions of \linebreak $\lambda \sigma_2(t_2) u - \sigma_1(t_2) \frac{d u}{dt_2} + \gamma(t_2) u
           = 0$ to solutions of \newline
          $ \lambda \sigma_2(t_2) y - \sigma_1(t_2) \frac{d y}{dt_2} + \gamma_*(t_2) y = 0$.
    \end{enumerate}
\end{prop}
Recall \cite{bib:NonConsTheory, bib:CoddLev} that the fourth property actually means that
\begin{equation} \label{eq:InttwS}
 S(\lambda,t_2) \Phi(\lambda, t_2,\tau_2) = \Phi_*(\lambda, t_2,\tau_2) S(\lambda,\tau_2)
\end{equation}
for fundamental matrices of the corresponding equations:
\begin{equation} \label{eq:ODEPhiOUT}
 \begin{array}{lll}
\lambda \sigma_2(y) \Phi_*(\lambda,y,\tau_2) - \sigma_1(y)
\frac{\partial}{\partial y}\Phi_*(\lambda,y,\tau_2)
+ \gamma_*(y) \Phi_*(\lambda,y,\tau_2) = 0,\\
 ~~~~~~~~~~~~~~~\Phi_*(\lambda,\tau_2,\tau_2) = I
\end{array}
\end{equation}
and
\begin{equation} \label{eq:ODEPhiIN}
 \begin{array}{lll}
\lambda \sigma_2(y) \Phi(\lambda,y,\tau_2) - \sigma_1(y)
\frac{\partial}{\partial y} \Phi(\lambda,y,\tau_2)
+ \gamma(y) \Phi(\lambda,y,\tau_2) = 0,\\
 ~~~~~~~~~~~~~~~\Phi(\lambda,\tau_2,\tau_2) = I.
\end{array}
\end{equation}
and as a result $S(\lambda,t_2)$ satisfies the following differential equation
\begin{equation} \label{eq:DS}
\frac{\partial}{\partial t_2} S(\lambda,t_2) = \sigma_1^{-1} (\sigma_2 \lambda + \gamma_*) S(\lambda,t_2) -
	S(\lambda,t_2)\sigma_1^{-1} (\sigma_2 \lambda + \gamma).
\end{equation}
\textbf{Proof of Proposition \ref{prop:Sdef}:} Those are easily checked properties, following from the definition of $S(\lambda,t_2)$:
\[ S(\lambda, t_2) = I - B^*(t_2) (\lambda I - A_1(t_2))^{-1} B(t_2) \sigma_1(t_2). \]
When $\lambda \rightarrow \infty$, since all the operators are
bounded the second summand vanishes and we obtain $ S(\infty, t_2)
= I_{n \times n}$. Moreover, it will be a meromorphic function of
$\lambda$, when $\lambda > \| A_1(t_2)\|$ and we obtain the first
property.

In order to understand the second property let us rewrite $S(\lambda, t_2)$, using
the zero curvature condition in the following way:
\[ \begin{array}{lll}
    S(\lambda, t_2) & = I - B^*(t_2) (\lambda I - A_1(t_2))^{-1} B(t_2) \sigma_1(t_2) = \\
       & = I - B^*(t_2) (\lambda I - F(t_2,\tau_2)A_1(\tau_2) F(\tau_2, t_2))^{-1} B(t_2) \sigma_1(t_2) = \\
       & = I - B^*(t_2) F(t_2,\tau_2)(\lambda I - A_1(\tau_2))^{-1} F(\tau_2, t_2) B(t_2) \sigma_1(t_2).
\end{array} \]
The functions $B^*(t_2) F(t_2,\tau_2)$, $F(\tau_2, t_2) B(t_2)$,
$\sigma_1(t_2)$ are absolutely continuous in appropriate spaces,
thus their multiplication too and we obtain the second property.

The third property is a result of straightforward calculations:
\[ \begin{array}{lllll}
S(\lambda, t_2)^* \sigma_1(t_2) S(\lambda, t_2) - \sigma_1(t_2) = \\
~~~~~~ 2 \Re (\lambda) \sigma_1(t_2)
B^*(t_2) (\bar\lambda I - A_1^*(t_2))^{-1}(\lambda I - A_1(t_2))^{-1} B(t_2) \sigma_1(t_2)
\end{array} \]
Here the sign of $\Re (\lambda)$ determines the sign of
$S(\lambda, t_2)^* \sigma_1(t_2) S(\lambda, t_2) - \sigma_1(t_2)$ and thus the third property is obtained.

The fourth and the last property is a direct result of our construction. \qed

\noindent\textbf{Remark:} In the case $B(t_2)$ is a compact operator (in particular
if $\dim \mathcal E < \infty$), $S(\lambda, t_2)$ is a meromorphic function on
$\mathbb{C}\setminus i\mathbb{R}$ for all $t_2$.

\begin{defn} The class
\[ \boldsymbol{\mathcal CI} = \boldsymbol{\mathcal CI}(\sigma_1,\sigma_2, \gamma,\gamma_*)
\]
is a class of functions $S(\lambda,t_2)$ of two variables, which
\begin{enumerate}
	\item are Identity in the neighborhood of $\lambda=\infty$ for all $t_2$,
	\item are absolutely continuous as functions of $t_2$ for almost all $\lambda$,
	\item satisfy:
          \[ \begin{array}{llll}
            S(\lambda, t_2)^* \sigma_1(t_2) S(\lambda, t_2) = \sigma_1(t_2), & \Re{\lambda} = 0 \\
            S(\lambda, t_2)^* \sigma_1(t_2) S(\lambda, t_2) \geq \sigma_1(t_2), & \Re{\lambda} \geq 0 \\
          \end{array} \]
          for $\lambda$ in the domain of analyticity of $S(\lambda,t_2)$.
	\item map solutions of the input ODE (\ref{eq:InCC}) with spectral parameter $\lambda$ to the output 
		ODE (\ref{eq:OutCC}) with the same spectral parameter (i.e., they satisfy the ODE (\ref{eq:DS}))
\end{enumerate}
\end{defn}
\section{\label{sec:Translation}Translation model}
Here we would like to present a model of a vessel that is the most
convenient for solving realization problem in class $\boldsymbol{\mathcal CI}$. The same model was used in
\cite{bib:TheoryNonComm} in the constant case. Let
\[ \mathfrak{IV} = (A_1(t_2), F(t_2,\tau_2), B(t_2); \sigma_1(t_2), \sigma_2(t_2), \gamma(t_2), \gamma_*(t_2);
\mathcal H_{t_2},\mathcal E)
\]
be a vessel associated to the system
\begin{align*}
  I\Sigma: \left\{ \begin{array}{lll}
    \frac{\partial}{\partial t_1} x(t_1,t_2) = A_1(t_2) x(t_1,t_2) + B(t_2) \sigma_1(t_2) u(t_1,t_2) \\
    x(t_1, t_2) = F (t_2,\tau_2) x(t_1, \tau_2) + \int\limits_{\tau_2}^{t_2} F(t_2, s) B(s) \sigma_2(s) u(t_1, s)ds \\
    y(t_1, t_2) = u(t_1, t_2) - B^*(t_2) x(t_1, t_2)
  \end{array} \right.
\end{align*}
with inputs and outputs satisfying compatibility conditions:
\[ \begin{array}{llll}
  \sigma_2(t_2) \frac{\partial}{\partial t_1}u(t_1, t_2) -
  \sigma_1(t_2) \frac{\partial}{\partial t_2}u(t_1,t_2) + \gamma(t_2) u(t_1,t_2) = 0 \\
  \sigma_2(t_2) \frac{\partial}{\partial t_1}y(t_1, t_2) -
  \sigma_1(t_2) \frac{\partial}{\partial t_2}y(t_1,t_2) + \gamma_*(t_2) y(t_1,t_2) = 0
\end{array} \]
If the input $u(t_1,t_2) \equiv 0$ and the initial inner state $x(t_1,t_2) = h, h \in \mathcal H_{t_2}$,
then obviously the function
\[ y_{h, t_2} (t_1, s) = \int\limits_{|\lambda| = R} e^{\lambda t_1} \Phi_*(\lambda,s,t_2)
B^*(t_2) (\lambda I - A_1(t_2))^{-1}h  d\lambda
\]
is an analytic function of two variables $t_1, s$, satisfying the output PDE with initial condition
\begin{equation} \label{eq:InitTrans}
 y_{h, t_2} (t_1,t_2) = B^*(t_2) e^{A_1(t_2) t_1}  h.
\end{equation}
In other words, it is an analytic output of the system in this case. Let us denote the space of all
these outputs as $\mathcal H_{trans, t_2}$ and convert it into a Hilbert space by introducing the induced
norm, namely
\[ \| y_{h, t_2} (t_1, s) \|_{\mathcal H_{trans, t_2}} = \| h \|_{\mathcal H}. \]
The following operators can be considered on this space:
\begin{defn} \label{def:A1F}
\[ \tilde A_1 = \frac{\partial}{\partial t_1}: \mathcal H_{trans, t_2} \rightarrow \mathcal H_{trans, t_2},
~~~~~  \tilde F(t_2', t_2) = Identity = I: \mathcal H_{trans, t_2} \rightarrow \mathcal H_{trans, t_2'}. \]
\end{defn}
In the next proposition we show that all the definitions are legal.
\begin{prop} 1. The operators $\tilde A_1, \tilde F(t_2', t_2)$ are well defined \textbf{on} the space
$\mathcal H_{trans, t_2}$. 2. All the spaces $\mathcal H_{trans, t_2}$ are identical as sets.
\end{prop}
\textbf{Proof:} 1. Let us first consider the operator $\tilde A_1$.
\begin{align*}
\tilde A_1 \big( y_{h, t_2} (t_1, s) \big) & = \frac{\partial}{\partial t_1}\big( y_{h, t_2} (t_1, s)\big)=
\frac{\partial}{\partial t_1}
\int\limits_{|\lambda| = R} e^{\lambda t_1} \Phi_*(\lambda,s,t_2)
 B^*(t_2) (\lambda I - A_1(t_2))^{-1}h  d\lambda \\
& = \int\limits_{|\lambda| = R} \lambda e^{\lambda t_1} \Phi_*(\lambda,s,t_2)
      B^*(t_2) (\lambda I -  A_1(t_2))^{-1}h  d\lambda = \\
& = \int\limits_{|\lambda| = R} e^{\lambda t_1} \Phi_*(\lambda,s,t_2)
      B^*(t_2) (\lambda I -  A_1(t_2))^{-1} [ \lambda -
      A_1(t_2) +  A_1(t_2)]h  d\lambda = \\
& = \int\limits_{|\lambda| = R} e^{\lambda t_1} \Phi_*(\lambda,s,t_2)
      B^*(t_2) (\lambda I -  A_1(t_2))^{-1}  A_1(t_2) h  d\lambda
\end{align*}
and since $A_1(t_2)$ is a bounded operator, $\tilde A_1$ is well-defined. Moreover, we can also see this way that
those operators are dual for the corresponding Hilbert spaces.

Let us consider now the action of $\tilde F(t_2', t_2)$. Since this function takes an element
$y_{h, t_2} (t_1, s) \in \mathcal H_{t_2, trans}$ and considers it as an element of
$\mathcal H_{t_2', trans}$, in order to find its norm one has to evaluate the function
$y_{h, t_2'} (t_1, s)$ at $y=t_2'$ and represent it in the form (\ref{eq:InitTrans}):
\[ y_{h, t_2'} (t_1, t_2') =  B^*(t_2') e^{A_1(t_2') t_1}h'.
\]
We claim that in the case of minimal systems, it is possible. Suppose that
\[ \bigcup\limits_{n=0}^\infty A_1^n(t_2) B(t_2) \sigma_1(t_2)e = \mathcal H_{t_2}, ~~~~\forall t_2. \]
Then we can evaluate the action of $\tilde F(t_2', t_2)$ on $y_{h, t_2} (t_1, s)$ for
$h = A_1^n(t_2) B(t_2) \sigma_1(t_2) e$. In this case
\begin{align*}
y_{h, t_2'} (t_1, t_2') & = \int\limits_{|\lambda| = R} e^{\lambda t_1} \Phi_*(\lambda,t_2', t_2)
      B^*(t_2) (\lambda I -  A_1(t_2))^{-1}  A_1^n(t_2) B(t_2) \sigma_1(t_2) e  d\lambda = \\
& = \int\limits_{|\lambda| = R} \lambda^n e^{\lambda t_1} \Phi_*(\lambda,t_2', t_2)
      B^*(t_2) (\lambda I -  A_1(t_2))^{-1} B(t_2) \sigma_1(t_2) e  d\lambda = \\
& = \int\limits_{|\lambda| = R} \lambda^n e^{\lambda t_1} \Phi_*(\lambda,t_2', t_2)
      S(\lambda,t_2) e  d\lambda = \text{by intertwining property} \\
& = \int\limits_{|\lambda| = R} \lambda^n e^{\lambda t_1} S(\lambda,t_2') \Phi(\lambda,t_2', t_2)  e  d\lambda = \\
& = \int\limits_{|\lambda| = R} \lambda^n e^{\lambda t_1} B^*(t_2') (\lambda I -  A_1(t_2'))^{-1} B(t_2') \sigma_1(t_2')
     \Phi(\lambda,t_2', t_2)  e  d\lambda = \\
& = \int\limits_{|\lambda| = R} e^{\lambda t_1} B^*(t_2')
    (\lambda I -  A_1(t_2'))^{-1} A_1^n(t_2') B(t_2') \sigma_1(t_2') \Phi(\lambda,t_2', t_2)  e  d\lambda.
\end{align*}
Thus we conclude that
\[ h' = \int\limits_{|\lambda| = R} (\lambda I -  A_1(t_2'))^{-1} A_1^n(t_2') B(t_2') \sigma_1(t_2')
     \Phi(\lambda,t_2', t_2)  e  d\lambda,
\]
which is a well-defined element of $\mathcal H_{t_2'}$. Moreover, we can conclude this way that the operator
$\tilde F(t_2', t_2)$ acts from $\mathcal H_{t_2}$ to $\mathcal H_{t_2'}$ one-to-one and since it is the
inverse of itself, it is onto. Finally, in all sets the vector spaces $\mathcal H_{t_2}$ are the same for all $t_2$. \qed
\begin{prop} The vessel
\[ \mathfrak{IV} = (A_1(t_2), F(t_2,\tau_2), B(t_2); \sigma_1(t_2), \sigma_2(t_2), \gamma(t_2), \gamma_*(t_2);
\mathcal H_{t_2},\mathcal E)
\]
is gauge equivalent to the vessel
\[ \mathfrak{\widetilde{IV}} = (\frac{\partial}{\partial t_1}, I, \big(Eval|_{t_1=0,t_2}\big)^*;
 \sigma_1(t_2), \sigma_2(t_2), \gamma(t_2), \gamma_*(t_2); \mathcal H_{t_2, trans},\mathcal E).
\]
\end{prop}
\textbf{Proof:} Take the space $H_{t_2, trans}$ as above and the operators $A_1(t_2), F(t_2',t_2)$ as in
Definition \ref{def:A1F}. Define the isometric map (via induced norm)
\[ U(t_2) h = y_{h, t_2} (t_1, s).
\]
In order to show that the vessels are gauge equivalent, it is enough to show that the conditions of the
equation (\ref{eq:CUconnectI}) are satisfied. Clearly
\[ U(t_2) A_1(t_2) h = y_{A_1(t_2)h, t_2} (t_1, s) =
\dfrac{\partial}{\partial t_1} y_{h, t_2} (t_1, s) = \widetilde A_1(t_2) y_{h, t_2} (t_1, s).
\]
and
\[ \begin{array}{llll}
 \langle U(t_2) B(t_2) e,  y_{h, t_2} (t_1, p) \rangle_{\mathcal H_{trans,t_2}} & =
 \langle B(t_2)e, U^*(t_2) y_{h, t_2} (t_1, p) \rangle_{\mathcal H_{t_2}} =
 \langle B(t_2)e, h \rangle_{\mathcal H_{t_2}} = \\
 & = \langle e, B^*(t_2) h \rangle_{\mathcal H_{t_2}} =
 \langle e, Eval|_{t_1=0,t_2} y_{h, t_2} (t_1, p) \rangle_{\mathcal H_{t_2}}
\end{array} \]
since
\[ \begin{array}{llll}
Eval|_{t_1=0,t_2} y_{h, t_2} (t_1, p) =
Eval|_{t_1=0,t_2} \int\limits_{|\lambda| = R} e^{\lambda t_1} \Phi_*(\lambda,p,t_2)
B^*(t_2) (\lambda I - A_1(t_2))^{-1}h  d\lambda = \\
\int\limits_{|\lambda| = R} B^*(t_2) (\lambda I - A_1(t_2))^{-1}h  d\lambda = B^*(t_2) h .
\end{array} \]
In order to show that $\widetilde{F}(t_2,\tau_2) = U(t_2) F(t_2,\tau_2) U^{-1}(\tau_2)$, we shall
Need more calculations.
\begin{align*}
U(t_2) F(t_2,\tau_2) U^{-1}(\tau_2) y_{h, \tau_2} (t_1, s)
 & = U(t_2) F(t_2,\tau_2) h = y_{F(t_2,\tau_2) h, t_2} (t_1, s) = \\
 & =  \widetilde F(t_2,\tau_2) y_{h, \tau_2} (t_1, s),
\end{align*}
as desired. \qed

This shows that for every vessel there is a naturally built translation model, gauge equivalent
to this vessel.

\section{\label{sec:RealCons}Realization theorem}
\subsection{Definition of the realization problem}
In this section our aim will be a realization problem. In other words, we shall start from
a transfer function $S(\lambda, t_2) \in \boldsymbol{\mathcal CI}$, i.e., satisfying conditions of Proposition \ref{prop:Sdef}:
\begin{enumerate}
        \item For almost all $t_2$, $S(\lambda, t_2)$ is an analytic function of $\lambda$
          in the neighborhood of $\infty$, where it satisfies:
        \[ S(\infty, t_2) = I_{n \times n} \]
      \item For all $\lambda$, $S(\lambda, t_2)$ is an absolutely continuous function of $t_2$.
        \item The following inequalities are satisfied:
          \[ \begin{array}{llll}
            S(\lambda, t_2)^* \sigma_1(t_2) S(\lambda, t_2) = \sigma_1(t_2), & \Re{\lambda} = 0 \\
            S(\lambda, t_2)^* \sigma_1(t_2) S(\lambda, t_2) \geq \sigma_1(t_2), & \Re{\lambda} \geq 0 \\
          \end{array} \]
          for $\lambda$ in the domain of analyticity of $S(\lambda,t_2)$.
        \item For each fixed $\lambda$, multiplication by $S(\lambda, t_2)$ maps
          solutions of \linebreak $\lambda \sigma_2(t_2) u - \sigma_1(t_2) \frac{d u}{dt_2} + \gamma(t_2) u
           = 0$ to solutions of \newline
          $ \lambda \sigma_2(t_2) y - \sigma_1(t_2) \frac{d y}{dt_2} + \gamma_*(t_2) y = 0$,
\end{enumerate}
and shall build a vessel with the very transfer function, using the translation model.

The vessel that we are going to build will be denoted in the following way:
\[ \mathfrak{IV} = (A_1(t_2), F(t_2,t_2^0), B(t_2); \sigma_1(t_2), \sigma_2(t_2), \gamma(t_2), \gamma_*(t_2);
\mathcal H_{t_2},\mathcal{E}),
\]
which will have as its transfer function the given $S(\lambda, t_2)$. \\
\subsection{Construction of the Vessel}
\textbf{Building of $H_{trans,t_2}, A_1(t_2), B(t_2)$:}
For the given $S(\lambda,t_2)$ we can realize it \cite{bib:JordanBrodskii} for each $t_2$ as
\[ S(\lambda,t_2) = I - \widehat B^*(t_2) (\lambda I- \widehat A_1(t_2))^{-1} \widehat B(t_2) \sigma_1(t_2) \]
satisfying
\[
\widehat A_1(t_2) + \widehat A_1^*(t_2) = - \widehat B(t_2) \sigma_1(t_2) \widehat B^*(t_2).
\]
for an abstract Hilbert space $\mathcal H_{t_2}$. Let us for every element $h \in \mathcal H_{t_2}$ define an element
$v_{h,t_2}(t_1,s)$ of a new Hilbert space $\mathcal{H}_{trans,t_2}$ as:
\[ \begin{array}{lll}
 v_{h,t_2}(t_1,s) = \int\limits_{|\lambda| = R} e^{\lambda t_1} \Phi_*(\lambda,s,t_2) \widehat B^*(t_2)
(\lambda I - \widehat A_1(t_2))^{-1}h  d\lambda, \\
 v_{h,t_2}(t_1,t_2) = e^{\widehat A_1(t_2) t_1} \widehat B^*(t_2) h.
\end{array} \]
Here $R > \| \widehat A_1(t_2)\|$ and $\Phi_*(\lambda,s,t_2)$ was defined in (\ref{eq:ODEPhiOUT}).
The norm for this element will be the range norm
$\| v_{h, t_2}(t_1, s)\|_{\mathcal{H}_{trans,t_2}} = \|h\|_{\mathcal H_{t_2}}$.
The operators $A_1(t_2): \mathcal H_{trans,t_2} \rightarrow \mathcal H_{trans,t_2}$ and \linebreak
$B^*(t_2): \mathcal H_{trans,t_2} \rightarrow \mathcal E$ are as follows:
\[ \begin{array}{lll}
A_1(t_2) = \frac{\partial}{\partial t_1}, \\
B^*(t_2) v_h(t_1,y) = v_h(0,t_2),
\end{array} \]
which by construction will satisfy the first colligation condition (\ref{eq:ConsCond}).

\noindent
\textbf{Definition of $F(t_2', t_2)$:} This operator acts as the identity operator from $H_{trans,t_2}$
to $H_{trans,t_2'}$. It means, in particular, that we have to show that as linear spaces, the spaces $\mathcal H_{trans,t_2}$ are all the same. By the minimality of the vessel we work with the vectors
$h = \widehat A_1(t_2)^n \widehat B^*(t_2) e$ (of the dense set). So, let us start from an element
$v_{h,t_2} (t_1,s)$ of the space $\mathcal H_{trans,t_2}$ with
$h = \widehat A_1(t_2)^n \widehat B^*(t_2) \sigma_1(t_2) e$, then
$F(t_2', t_2) v_{h,t_2} (t_1,s)$ is the same function, but considered in the space $\mathcal H_{trans,t_2'}$:
\[ \begin{array}{lll}
v_{h,t_2} (t_1,s)   & = \int\limits_{|\lambda| = R} e^{\lambda t_1} \Phi_*(\lambda,s,t_2)
      \widehat B^*(t_2)(\lambda I-\widehat A_1(t_2))^{-1} h ]d\lambda \\
        & = \frac{\partial}{\partial t_1^n} [
\int\limits_{|\lambda| = R}
       e^{\lambda t_1} \Phi_*(\lambda,s,t_2) S(\lambda,t_2)  e ]d\lambda \\
   & = \text{using condition (4)!} \\
   & = \frac{\partial}{\partial t_1^n} [ \int\limits_{|\lambda| = R}
       e^{\lambda t_1} S(\lambda,s) \Phi(\lambda,s,t_2) e ]d\lambda \\
   & = \frac{\partial}{\partial t_1^n} [ \int\limits_{|\lambda| = R}
       e^{\lambda t_1} \widehat B^*(s) (\lambda I - \widehat A_1(s))^{-1} \widehat B(s)
       \sigma_1(t_2) \Phi(\lambda,s,t_2) e] d\lambda \\
   & = \int\limits_{|\lambda| = R}
       e^{\lambda t_1} \widehat B^*(s) (\lambda I - \widehat A_1(s))^{-1} \widehat A_1(s)^n \widehat B(s)
       \sigma_1(t_2) \Phi(\lambda,s,t_2) e] d\lambda \\
   & = v_{h',t_2'} (t_1,s),
\end{array} \]
where $h'$ can be evaluated from assigning $s = t_2'$ in the last expression (using the definition) and thus
\begin{equation} \label{eq:h'}
h' = \int\limits_{|\lambda| = R} (\lambda I - \widehat A_1(t_2'))^{-1}
\widehat A_1^n(t_2') \widehat B(t_2') \sigma_1(t_2') \Phi(\lambda, t_2',t_2) e d\lambda.
\end{equation}
Let us see what this formula means in the language of the reproducing kernel Hilbert spaces, introduced in
the definition of $\mathcal H_{t_2, trans}$. According to that theory,
\begin{align}
\dfrac{\sigma_1(t_2) S(\mu,t_2) - S(\lambda,t_2)\sigma_1(t_2)}{\lambda+\mu} & =
\sigma_1(t_2) B^*(t_2) (\lambda I- A_1(t_2))^{-1} (\mu I- A_1(t_2))^{-1} B(t_2) \sigma_1(t_2) = \\
\label{eq:KerDecomp}& = \sigma_1(t_2) B^*(t_2) (\lambda I- A_1(t_2))^{-1} \dfrac{1}{\mu}\sum_{k=0}^\infty \frac{A_1(t_2)}{\mu} B(t_2) \sigma_1(t_2),
\end{align}
where we actually look at this expression as a function of $\lambda$ for each $\mu$ fixed. Afterwards,
linear expressions are taken for different $\mu$'s.
Now define an isometric isomorphism $T$ from this space to $\mathcal H_{t_2, trans}$ by:
\[ T(t_2)[B^*(t_2) (\lambda I - A_1(t_2))^{-1} h ] = \int_{|\lambda| = R} e^{\lambda t_1}
\Phi_*(\lambda, s, t_2) B^*(t_2) (\lambda I - A_1(t_2))^{-1} h d\lambda
\]
with the range norm. Then we shall try with the help of this map to find how the operator $F(t_2',t_2)$ acts
on the kernels (i.e., on the functions of $\lambda$ of this specific form). We shall actually evaluate
\[ T^{-1}(t_2') F(t_2',t_2) T(t_2) [\dfrac{\sigma_1(t_2) S(\mu,t_2) - S(\lambda,t_2)\sigma_1(t_2)}{\lambda+\mu}].
\]
From the decomposition (\ref{eq:KerDecomp}) we evaluate:
\begin{align*}
T(t_2) [\dfrac{\sigma_1(t_2) S(\mu,t_2) - S(\lambda,t_2)\sigma_1(t_2)}{\lambda+\mu}] =
T(t_2) [\sigma_1(t_2) B^*(t_2) (\lambda I- A_1(t_2))^{-1} \dfrac{1}{\mu}\sum_{k=0}^\infty \frac{A_1(t_2)}{\mu} B(t_2) \sigma_1(t_2)] = \\
= \sigma_1(t_2) \int_{|\lambda| = R} e^{\lambda t_1} \Phi_*(\lambda, s, t_2)
B^*(t_2) (\lambda I- A_1(t_2))^{-1} \dfrac{1}{\mu}\sum_{k=0}^\infty \frac{A_1(t_2)}{\mu}B(t_2)\sigma_1(t_2) d\lambda =\\
= \ldots = \\
= \sigma_1(t_2) \sum_{k=0}^\infty \dfrac{1}{\mu^{k+1}} \dfrac{\partial^k}{\partial t_1^k}
\int_{|\lambda| = R} e^{\lambda t_1} \Phi_*(\lambda, s, t_2) S(\lambda,t_2) d\lambda = \\
= \sigma_1(t_2) \sum_{k=0}^\infty \dfrac{1}{\mu^{k+1}} \dfrac{\partial^k}{\partial t_1^k}
\int_{|\lambda| = R} e^{\lambda t_1} S(\lambda,s) \Phi(\lambda, s, t_2)  d\lambda = \\
= \sigma_1(t_2) \sum_{k=0}^\infty \dfrac{1}{\mu^{k+1}}\int_{|\lambda| = R} e^{\lambda t_1}
B^*(s) (\lambda I- A_1(s))^{-1} A_1(s)^k B(s) \sigma_1(s) \Phi (\lambda,s,t_2) d\lambda
\end{align*}
Suppose now that $\Phi (\lambda,s,t_2) = \sum\limits_{i=0}^\infty \Phi^i(s,t_2) \lambda^i$, then the last
expression becomes:
\begin{gather*}
\sigma_1(t_2) \sum_{k=0}^\infty \dfrac{1}{\mu^{k+1}}\int_{|\lambda| = R} e^{\lambda t_1}
B^*(s) (\lambda I- A_1(s))^{-1} A_1(s)^k B(s) \sigma_1(s) \sum\limits_{i=0}^\infty \Phi^i(s,t_2) \lambda^i d\lambda = \\
= \sigma_1(t_2) \sum_{k=0}^\infty \dfrac{1}{\mu^{k+1}} \sum\limits_{i=0}^\infty \int_{|\lambda| = R} e^{\lambda t_1}
B^*(s) (\lambda I- A_1(s))^{-1} A_1(s)^{k+i} B(s) \sigma_1(s) \Phi^i(s,t_2) d\lambda
\end{gather*}

Now we evaluate the action of $T^{-1}(t_2') F(t_2',t_2) = T^{-1}(t_2')$. For this is, we have
to apply evaluation at $(t_1,s) = (0,t_2')$ and in this way obtain the vector of $\mathcal H$. Then, we
shall insert $B^*(t_2') (\lambda I- A_1(t_2'))^{-1}$ without integral and $e^{t_1 \lambda}$,
which is the mapping back to the RKHS:
\begin{gather*}
\sigma_1(t_2) \sum_{k=0}^\infty \dfrac{1}{\mu^{k+1}}\sum\limits_{i=0}^\infty
B^*(t_2') (\lambda I- A_1(t_2'))^{-1} A_1(t_2')^{k+i} B(t_2') \sigma_1(t_2') \Phi^i(t_2',t_2) =\\
= \sigma_1(t_2) B^*(t_2') (\lambda I- A_1(t_2'))^{-1} (\mu I- A_1(t_2'))^{-1}
\sum\limits_{i=0}^\infty A_1(t_2')^i B(t_2') \sigma_1(t_2') \Phi^i(t_2',t_2).
\end{gather*}
To conclude,
\begin{multline*}
 T^{-1}(t_2') F(t_2',t_2) T(t_2) [\dfrac{\sigma_1(t_2) S(\mu,t_2) - S(\lambda,t_2)\sigma_1(t_2)}{\lambda+\mu}] = \\
\sigma_1(t_2) B^*(t_2') (\lambda I- A_1(t_2'))^{-1} (\mu I- A_1(t_2'))^{-1}
\sum\limits_{i=0}^\infty A_1(t_2')^i B(t_2') \sigma_1(t_2') \Phi^i(t_2',t_2).
\end{multline*}

\subsection{Properties of fundamental matrices $\Phi_*(\lambda,t_2,t_2^0), \Phi(\lambda,t_2,t_2^0)$}
Let us remember first the definitions of these matrices. They were defined as the fundamental matrices of the
corresponding ODEs (\ref{eq:ODEPhiOUT})
\begin{equation*}
 \begin{array}{lll}
\lambda \sigma_2(y) \Phi_*(\lambda,y,t_2^0) - \sigma_1(y)
\frac{\partial}{\partial y}\Phi_*(\lambda,y,t_2^0)
+ \gamma_*(y) \Phi_*(\lambda,y,t_2^0) = 0,\\
 ~~~~~~~~~~~~~~~\Phi_*(\lambda,t_2^0,t_2^0) = I
\end{array}
\end{equation*}
and (\ref{eq:ODEPhiIN})
\begin{equation*}
 \begin{array}{lll}
\lambda \sigma_2(y) \Phi(\lambda,y,t_2^0) - \sigma_1(y)
\frac{\partial}{\partial y} \Phi(\lambda,y,t_2^0)
+ \gamma(y) \Phi(\lambda,y,t_2^0) = 0,\\
 ~~~~~~~~~~~~~~~\Phi(\lambda,t_2^0,t_2^0) = I.
\end{array}
\end{equation*}
\begin{prop} The following formulas are correct
\[ \begin{array}{lll} \label{eq:PhiFirst}
\sigma_1(t_2) \Phi_*(\lambda,t_2,t_2^0) = (\Phi_*)^{-1*}(-\bar\lambda,t_2,t_2^0) \sigma_1(t_2^0), \\
\sigma_1(t_2) \Phi(\lambda,t_2,t_2^0) = (\Phi)^{-1*}(-\bar\lambda,t_2,t_2^0) \sigma_1(t_2^0). \\
\end{array}
\]
\end{prop}
\textbf{Proof:} we start from formula (\ref{eq:ODEPhiOUT}) and conjugate it:
\[
\bar\lambda (\Phi_*)^*(\lambda,t_2,t_2^0) \sigma_2(t_2) -
\frac{\partial}{\partial t_2}\big( (\Phi_*)^*(\lambda,t_2,t_2^0) \big) \sigma_1(t_2)
+ (\Phi_*)^*(\lambda,y,t_2^0) \gamma(t_2)_* = 0
\]
or
\[
\frac{\partial}{\partial t_2}(\Phi_*)^*(\lambda,t_2,t_2^0)  =
 (\Phi_*)^*(\lambda,t_2,t_2^0) [ \bar\lambda \sigma_2(t_2) + \gamma(t_2)_*^s - \frac{d}{dt_2} \sigma_1(t_2)]
\sigma_1(t_2)^{-1} \]
and taking the inverse:
\[
\frac{\partial}{\partial t_2} (\Phi_*)^{-1*}(\lambda,t_2,t_2^0) =
 [- \bar\lambda \sigma_2(t_2) - \gamma(t_2)_*^s + \frac{d}{dt_2} \sigma_1(t_2)] \sigma_1(t_2)^{-1}
  (\Phi_*)^{-1*}(\lambda,t_2,t_2^0).
\]
On the other hand, formula (\ref{eq:ODEPhiOUT}) may be rewritten as
\[ \begin{array}{lll}
\frac{\partial}{\partial t_2} \big( \sigma_1(t_2) \Phi_*(\lambda,t_2,t_2^0) \big) = \\
\{ \lambda \sigma_2(t_2) - \gamma_*^s(t_2) + \frac{d}{dt_2} \sigma_1(t_2)\} \Phi_*(\lambda,y,t_2^0) = \\
\{ \lambda \sigma_2(t_2) - \gamma_*^s(t_2) + \frac{d}{dt_2} \sigma_1(t_2)\}
 \sigma_1^{-1}(t_2) \sigma_1(t_2) \Phi_*(\lambda,t_2,t_2^0).
\end{array}
\]
From the last two formulas we easily conclude that $\sigma_1(t_2) \Phi_*(\lambda,t_2,t_2^0)$ and
$(\Phi_*)^{-1*}(-\bar\lambda,t_2,t_2^0)$ are linearly dependent, and using the initial condition, we find that
\[ \sigma_1(t_2) \Phi_*(\lambda,t_2,t_2^0) = (\Phi_*)^{-1*}(-\bar\lambda,t_2,t_2^0) \sigma_1(t_2^0).
\]
The formula for $\Phi(\lambda,t_2,t_2^0)$ is proved in the same way. \qed
\begin{prop}
The following formulas are correct
\[ \begin{array}{llllll} \label{eq:PhiSecond}
\frac{\partial}{\partial t_2} [\Phi_*^*(\mu,t_2,t_2^0) \sigma_1(t_2) \Phi_*(\lambda,t_2,t_2^0)] =
(\lambda +\bar\mu)\Phi_*^*(\mu,t_2,t_2^0) \sigma_2(t_2)\Phi_*(\lambda,t_2,t_2^0), \\
\frac{\partial}{\partial t_2} [\Phi^*(\mu,t_2,t_2^0) \sigma_1(t_2) \Phi(\lambda,t_2,t_2^0)] =
(\lambda +\bar\mu)\Phi^*(\mu,t_2,t_2^0) \sigma_2(t_2)\Phi(\lambda,t_2,t_2^0).
\end{array} \]
\end{prop}
\textbf{Proof:} This is a rather technical result: straightforward calculations, using differential equations
(\ref{eq:ODEPhiOUT}) and (\ref{eq:ODEPhiIN}).
\qed

\subsection{Proof of vessel conditions}
The first colligation condition we have obtained from the realization theorem for each $t_2$. Now we are going to show that all the other conditions 
holds too.

\noindent \textbf{Second colligation condition}. This condition is as follows:
\[ \frac{\partial }{\partial t_2} (\widehat F^*(t_2,t_2^0) \widehat F(t_2,t_2^0)) =
 - F^*(t_2,t_2^0) \widehat B(t_2) \sigma_2(t_2) \widehat B^*(t_2) F(t_2,t_2^0).
\]
Take an element of the dense set
$h_{t_2^0} = A_1^n(t_2^0)B(t_2^0) \sigma_1(t_2^0) e_j \in \mathcal H_{t_2^0}$,
then it is sufficient to prove that
\[ \begin{array}{llllll}
 \frac{\partial }{\partial t_2} \langle \widehat F(t_2,t_2^0) v_{h_{t_2^0}}(t_1,s), \widehat F(t_2,t_2^0)
 v_{h_{t_2^0}}(t_1,s) \rangle_{\mathcal{H}_{trans,t_2}} = \\
 \langle \sigma_2(t_2) \widehat B^*(t_2) \widehat F(t_2,t_2^0) v_{h_{t_2^0}}(t_1,s),
 \widehat B^*(t_2) \widehat F(t_2,t_2^0)v_{h_{t_2^0}}(t_1,s)  \rangle_\mathcal{E}.
\end{array} \]
So, we concentrate on the left side
\[ \begin{array}{llllll}
\frac{\partial}{\partial t_2} \| \widehat F(t_2,t_2^0)
v_h(t_1,s)\|_{\mathcal{H}_{trans,t_2^0}} =
\frac{\partial}{\partial t_2} \| v_h(t_1,s)\|_{\mathcal{H}_{trans,t_2}} = \\
\frac{\partial}{\partial t_2} \| \int\limits_{|\lambda| = R}
(\lambda I - A_1(t_2))^{-1} A_1^n(t_2) B(t_2) \sigma_1(t_2)
\Phi(\lambda,t_2,t_2^0) e_j d\lambda \|_{\mathcal{H}_{t_2}}
= \\
\frac{\partial}{\partial t_2} \| \int\limits_{|\lambda| = R}
(\lambda I - A_1(t_2))^{-1} \lambda^n B(t_2) \sigma_1(t_2)
\Phi(\lambda,t_2,t_2^0) e_j d\lambda \|_{\mathcal{H}_{t_2}}
= \\
\frac{\partial}{\partial t_2} \langle \int\limits_{|\lambda| = R_1}
(\lambda I - A_1(t_2))^{-1} \lambda^n B(t_2) \sigma_1(t_2) \Phi(\lambda,t_2,t_2^0) e_j d\lambda, \\
~~~~~~~~~~~~ \int\limits_{|\mu| = R_2}
(\mu I - A_1(t_2))^{-1} \mu^n B(t_2) \sigma_1(t_2) \Phi(\mu,t_2,t_2^0) e_j d\mu \|_{\mathcal{H}_{t_2}} =\\
\frac{\partial}{\partial t_2} \int\limits_{|\lambda| = R_1}
\int\limits_{|\mu| = R_2} \lambda^n \bar\mu^n \langle \sigma_1(t_2)
B^*(t_2) (\bar\mu I - A_1^*(t_2))^{-1}(\lambda I - A_1(t_2))^{-1} B(t_2) \sigma_1(t_2) \times \\
~~~~~~~~~~~~~~~~~~~~~~~~
\Phi(\lambda,t_2,t_2^0) e_j, \Phi(\mu,t_2,t_2^0) e_j \rangle_\mathcal{E} d\lambda d\mu = \\
\frac{\partial}{\partial t_2} \int\limits_{|\lambda| = R_1}
\int\limits_{|\mu| = R_2} \lambda^n \bar\mu^n \langle
\dfrac{S^*(\mu,t_2)\sigma_1(t_2) S(\lambda,t_2)-\sigma_1(t_2)
}{\lambda + \bar\mu}\Phi(\lambda,t_2,t_2^0) e_j,
   \Phi(\mu,t_2,t_2^0) e_j \rangle_\mathcal{E} d\lambda d\mu = \\
\frac{\partial}{\partial t_2} \int\limits_{|\lambda| = R_1}
\int\limits_{|\mu| = R_2} \lambda^n \bar\mu^n
\langle \dfrac{S^*(\mu,t_2^0)\Phi_*^*(\mu,t_2,t_2^0) \sigma_1(t_2) \Phi_*(\lambda,t_2,t_2^0)S(\lambda,t_2^0)}{\lambda + \bar\mu} e_j,e_j \rangle_\mathcal{E} d\lambda d\mu + \\
~~~~~~~~~~~~~ \frac{\partial}{\partial t_2} \int\limits_{|\lambda|
= R_1} \int\limits_{|\mu| = R_2} \lambda^n \bar\mu^n \langle
\dfrac{\Phi^{* in}(\mu,t_2,t_2^0)\sigma_1(t_2)
\Phi(\lambda,t_2,t_2^0)}{\lambda + \bar\mu} e_j, e_j
\rangle_\mathcal{E} d\lambda d\mu = \\
\frac{\partial}{\partial t_2} \int\limits_{|\lambda| = R_1}
\int\limits_{|\mu| = R_2} \lambda^n \bar\mu^n
\langle \dfrac{S^*(\mu,t_2^0)\Phi_*^*(\mu,t_2,t_2^0) \sigma_1(t_2) \Phi_*(\lambda,t_2,t_2^0)S(\lambda,t_2^0)}{\lambda + \bar\mu} e_j,e_j \rangle_\mathcal{E} d\lambda d\mu =
\end{array} \]
Now we can differentiate to obtain:
\[ \begin{array}{llllll}
= \int\limits_{|\lambda| = R_1} \int\limits_{|\mu| = R_2} \lambda^n \bar\mu^n
\langle S^*(\mu,t_2^0)\Phi_*^*(\mu,t_2,t_2^0) \sigma_2(t_2) \Phi_*(\lambda,t_2,t_2^0)S(\lambda,t_2^0) e_j,e_j \rangle_\mathcal{E} d\lambda d\mu,
\end{array} \]
where we have used the property of Proposition \ref{eq:PhiSecond}.

For the final result, reversing now the operations we have done before, we obtain that
the last expression is:
\[ \begin{array}{llllll}
= \int\limits_{|\lambda| = R_1} \int\limits_{|\mu| = R_2} \lambda^n \bar\mu^n
  \langle S^*(\mu,t_2^0)\Phi_*^*(\mu,t_2,t_2^0) \sigma_2(t_2) \Phi_*(\lambda,t_2,t_2^0)S(\lambda,t_2^0) e_j,
     e_j \rangle_\mathcal{E} d\lambda d\mu = \\
= \int\limits_{|\lambda| = R_1} \int\limits_{|\mu| = R_2}
  \lambda^n \bar\mu^n \langle \Phi^{*in}(\mu,t_2,t_2^0) S^*(\mu,t_2)
     \sigma_2(t_2) S(\lambda,t_2) \Phi(\lambda,t_2,t_2^0)e_j,
        e_j \rangle_\mathcal{E} d\lambda d\mu = \\
= \int\limits_{|\lambda| = R_1} \int\limits_{|\mu| = R_2}
  \lambda^n \bar\mu^n \langle \sigma_2(t_2) S(\lambda,t_2)
    \Phi(\lambda,t_2,t_2^0)e_j,S(\mu,t_2) \Phi(\mu,t_2,t_2^0)e_j
       \rangle_\mathcal{E} d\lambda d\mu = \\
= \int\limits_{|\lambda| = R_1} \int\limits_{|\mu| = R_2} \lambda^n \bar\mu^n
\langle \sigma_2(t_2) B^*(t_2) (\lambda I - A_1(t_2))^{-1} B(t_2) \sigma_1(t_2) \Phi(\lambda,t_2,t_2^0)e_j, \\
~~~~~~~~~~ B^*(t_2) (\mu I - A_1(t_2))^{-1}B(t_2) \sigma_1(t_2) \Phi(\mu,t_2,t_2^0)e_j
     \rangle_\mathcal{E} d\lambda d\mu = \\
= \langle \sigma_2(t_2) \int\limits_{|\lambda| = R_1} B^*(t_2) (\lambda I - A_1(t_2))^{-1}B(t_2) \sigma_1(t_2)
    \Phi(\lambda,t_2,t_2^0)e_j, \\
~~~~~~~~~~ \int\limits_{|\mu| = R_2} B^*(t_2) (\mu I -
A_1(t_2))^{-1}B(t_2) \sigma_1(t_2) \Phi(\mu,t_2,t_2^0)e_j
\rangle_\mathcal{E} d\lambda d\mu = \\
= \langle \sigma_2(t_2) \widehat B^*(t_2) v_h(t_1,s), \widehat
B(t_2)^*v_h(t_1,s)  \rangle_\mathcal{E},
\end{array} \]
since the integral of an analytic in $\lambda$ function vanishes
\[
\int\limits_{|\lambda| = R} \int\limits_{|\mu| = R} i \lambda^n
\bar\mu^n \langle \sigma_2(t_2) \Phi(\lambda,t_2,t_2^0)
e_j,\Phi(\mu,t_2,t_2^0) e_j \rangle_\mathcal{E} d\lambda d\mu
= 0.
\]

\noindent\textbf{Lax equation} is as follows:
\[ \widehat{F}(t_2, t_2^0) \widehat{A}_1(t_2^0) = \widehat{A}_1(t_2) \widehat{F}(t_2, t_2^0).
\]
In order to prove it, apply two sides of the equation to an element $v_h (t_1,s)$ to obtain
\[ \widehat{F}(t_2, t_2^0) \widehat{A}_1(t_2^0) v_h (t_1,s) = v_{A_1(t_2^0)}(t_1,s)
\]
and
\[ \begin{array}{llllll}
\widehat{A}_1(t_2) \widehat{F}(t_2, t_2^0) v_h (t_1,s)
   & = \widehat{A}_1(t_2) v_h(t_1,s) = \frac{\partial}{\partial t_1} v_h(t_1,s) = \\
   & = \frac{\partial}{\partial t_1} \int\limits_{|\lambda| = R}
         e^{\lambda t_1} \Phi_*(\lambda,t_2,t_2^0) B^*(t_2^0)
         (\lambda I - A_1(t_2^0))^{-1}h  d\lambda = \\
   & = \int\limits_{|\lambda| = R} \lambda
         e^{\lambda t_1} \Phi_*(\lambda,t_2,t_2^0) B^*(t_2^0)
         (\lambda I - A_1(t_2^0))^{-1} h  d\lambda = \\
   & = \int\limits_{|\lambda| = R} e^{\lambda t_1} \Phi_*(\lambda,t_2,t_2^0) B^*(t_2^0)
       (\lambda I - A_1(t_2^0))^{-1} A_1(t_2^0)h  d\lambda \\
   & = v_{A_1(t_2^0) h}(t_1,t_2),
\end{array}
\]
we have added and subtracted to obtain
\[ \int\limits_{|\lambda| = R} e^{\lambda t_1} \Phi_*(\lambda,t_2,t_2^0) B^*(t_2^0)
   (\lambda I - A_1(t_2^0))^{-1} A_1(t_2^0)h  d\lambda.
\]

\noindent \textbf{Linkage condition}. It is a standard result \cite{bib:HarmAnal,bib:Inter} of the theory of inner functions, that
\[ \lim\limits_{\lambda \longrightarrow \infty} \frac{\partial}{\partial t_2} S(\lambda, t_2)= 0. \]
Differentiating the formula $S(\lambda, t_2) =
\Phi_*(\lambda,t_2,t_2^0)S(\lambda, t_2^0)\Phi^{-1}(\lambda,t_2,t_2^0)$, we shall obtain:
\[ \begin{array}{lll}
 \lim\limits_{\lambda \longrightarrow \infty} \frac{\partial}{\partial t_2}S(\lambda, t_2) =
 \lim\limits_{\lambda \longrightarrow \infty}\sigma_1(t_2)^{-1} (\sigma_2(t_2) \lambda + \gamma_*^s(t_2)) S(\lambda, t_2) -\\
 ~~~~~~~~~~~~~~~~~
 \lim\limits_{\lambda \longrightarrow \infty} S(\lambda,t_2) \sigma_1(t_2)^{-1} (\sigma_2(t_2) \lambda + \gamma^s(t_2)).
\end{array} \]
Finally, inserting the realization formula for $S(\lambda, t_2)$,
we obtain:
\[ 0 = \sigma_1(t_2)^{-1} (\sigma_2(t_2) \widehat B^*(t_2) \widehat B(t_2) \sigma_1(t_2) + \gamma) -
\sigma_1(t_2)^{-1} (\sigma_2(t_2) \widehat B^*(t_2) \widehat B(t_2) \sigma_1(t_2) +
\gamma_*),
\]
from which the linkage condition is obtained:
\[ \gamma_*(t_2) = \gamma(t_2) + (\sigma_2(t_2) \widehat B^*(t_2) \widehat B(t_2) \sigma_1(t_2) -
\sigma_1(t_2) \widehat B^*(t_2) \widehat B(t_2) \sigma_2(t_2)).\] A similar
formula will be correct for the model operators too, because of
unitary equivalence. Namely,
\[ \gamma_*(t_2) = \gamma(t_2) + \sigma_2(t_2) B^*(t_2) B(t_2) \sigma_1(t_2) -
\sigma_1(t_2) B^*(t_2) B(t_2) \sigma_2(t_2)).\]

\noindent\textbf{Input vessel} condition:
\[ \frac{d}{d t_2}\big( \widehat{F}^*(t_2,t_2^0)\widehat B(t_2) \sigma_1(t_2) \big) -
 \widehat{F}^*(t_2,t_2^0)\widehat{A}_1(t_2)\widehat B(t_2) \sigma_2(t_2)
  + \widehat{F}^*(t_2,t_2^0)\widehat B(t_2) \gamma_*(t_2) = 0
\]
can be checked using bilinear forms. For this we shall take an
arbitrary element $A_1^n(t_2^0) B(t_2) \sigma_1(t_2^0) e$
\[ \begin{array}{lllllll}
 \frac{d}{d t_2} \langle F^*(t_2,t_2^0) B(t_2) \sigma_1(t_2),
       A_1^n(t_2^0) B(t_2) \sigma_1(t_2^0) e \rangle - \\
 \langle F^*(t_2,t_2^0){A}_1(t_2) B(t_2) \sigma_2(t_2)
  + F^*(t_2,t_2^0) B(t_2) \gamma_*(t_2),
    A_1^n(t_2^0) B(t_2) \sigma_1(t_2^0) e \rangle = \\
 \frac{d}{d t_2} \langle  B(t_2) \sigma_1(t_2), F(t_2,t_2^0)
       A_1^n(t_2^0) B(t_2) \sigma_1(t_2^0) e \rangle + \\
 \langle -A_1(t_2) B(t_2) \sigma_2(t_2)
  + B(t_2) \gamma_*(t_2), F(t_2,t_2^0)
    A_1^n(t_2^0) B(t_2) \sigma_1(t_2^0) e \rangle = \\
  \frac{d}{d t_2} \langle  B(t_2) \sigma_1(t_2),
     \int\limits_{|\lambda| = R} \lambda^n(\lambda I - A_1(t_2))^{-1} B(t_2) \sigma_1(t_2)
     \Phi(\lambda,t_2,t_2^0) e d\lambda \rangle + \\
  \langle - A_1^*(t_2) B(t_2) \sigma_2(t_2)
  +  B(t_2) \gamma_*(t_2),
     \int\limits_{|\lambda| = R} \lambda^n(\lambda I - A_1(t_2))^{-1} B(t_2) \sigma_1(t_2)
     \Phi(\lambda,t_2,t_2^0) e d\lambda \rangle = \\
  \frac{d}{d t_2} \langle \sigma_1(t_2),
     \int\limits_{|\lambda| = R} B^*(t_2) \lambda^n(\lambda I - A_1(t_2))^{-1} B(t_2) \sigma_1(t_2)
     \Phi(\lambda,t_2,t_2^0) e d\lambda \rangle - \\
  \langle \sigma_2(t_2), \int\limits_{|\lambda| = R} B^*(t_2) \lambda^{n+1} (\lambda I - A_1(t_2))^{-1} B(t_2) \sigma_1(t_2)
     \Phi(\lambda,t_2,t_2^0) e d\lambda \rangle + \\
  \langle \gamma_*(t_2),
     \int\limits_{|\lambda| = R}  B^*(t_2) \lambda^n(\lambda I - A_1(t_2))^{-1} B(t_2) \sigma_1(t_2)
     \Phi(\lambda,t_2,t_2^0) e d\lambda \rangle = \\
  -\frac{d}{d t_2} \langle \sigma_1(t_2), \int\limits_{|\lambda| = R}
     \lambda^n S(\lambda, t_2) \Phi(\lambda,t_2,t_2^0) e d\lambda \rangle +
     \langle \sigma_2(t_2), \int\limits_{|\lambda| = R}
     \lambda^{n+1} S(\lambda, t_2) \Phi(\lambda,t_2,t_2^0) e d\lambda \rangle - \\
  ~~~~~~~~-\langle \gamma_*(t_2), \int\limits_{|\lambda| = R}
     \lambda^n S(\lambda, t_2) \Phi(\lambda,t_2,t_2^0) e d\lambda \rangle = \\
  -\frac{d}{d t_2} \langle \sigma_1(t_2), \int\limits_{|\lambda| = R}
     \lambda^n \Phi_*(\lambda,t_2,t_2^0) S(\lambda, t_2^0) e d\lambda \rangle +
     \langle \sigma_2(t_2), \int\limits_{|\lambda| = R}
     \lambda^{n+1} \Phi_*(\lambda,t_2,t_2^0) S(\lambda, t_2^0) e d\lambda \rangle - \\
  ~~~~~~~~- \langle \gamma_*(t_2), \int\limits_{|\lambda| = R}
     \lambda^n \Phi_*(\lambda,t_2,t_2^0) S(\lambda, t_2^0) e d\lambda \rangle = \\
  -\langle I, \int\limits_{|\lambda| = R} \lambda^n [
     \frac{d}{dt_2} \big( \sigma_1(t_2) \Phi_*(\lambda,t_2,t_2^0)
     \big) -( \lambda \sigma_2 - \gamma_*^*(t_2)) \Phi_*(\lambda,t_2,t_2^0) ]
     S(\lambda, t_2^0) e d\lambda \rangle = \\
  \langle I, \int\limits_{|\lambda| = R} \lambda^n [
     - \sigma_1(t_2) \frac{d}{dt_2} \Phi_*(\lambda,t_2,t_2^0)
     + (\lambda \sigma_2 + \gamma_*(t_2))
     \Phi_*(\lambda,t_2,t_2^0)]
     S(\lambda, t_2^0) e d\lambda \rangle =0,
\end{array}  \]
because it is the differential equation for $\Phi_*(\lambda,t_2,t_2^0)$
(\ref{eq:ODEPhiIN}).

\noindent \textbf{Output vessel condition} in the integrated form is the following
\[ \begin{array}{lll}
 \frac{d}{dt_2} \big( \widehat{F}(t_2^0, t_2) \widehat B(t_2)
 \sigma_1(t_2) \big) +  \widehat F(t_2^0, t_2) \widehat{A}_1(t_2) \widehat B(t_2) \sigma_2(t_2) +
 \widehat F(t_2^0, t_2) \widehat B(t_2) \gamma(t_2) = 0
\end{array} \]
Let us use the formula for $\widehat{F}(t_2^0, t_2)$. Then the
last formula becomes:
\[ \begin{array}{lll}
 \frac{d}{dt_2} \big(
 \int\limits_{|\lambda| = R} (\lambda I - A_1(t_2^0))^{-1} B(t_2^0) \sigma_1(t_2^0)
 \Phi(\lambda,t_2^0,t_2) e d\lambda \big) + \\
  \int\limits_{|\lambda| = R} (\lambda I - A_1(t_2^0))^{-1} B(t_2^0) \sigma_1(t_2^0)
 \Phi(\lambda,t_2^0,t_2) \sigma_1(t_2)^{-1} [\lambda \sigma_2(t_2) + \gamma(t_2)] e d\lambda =  \\
 \int\limits_{|\lambda| = R} (\lambda I - A_1(t_2^0))^{-1} B(t_2^0) \sigma_1(t_2^0)
  \times \\
 ~~~~~~~~~~~~~~~
 [\frac{d}{dt_2} \Phi(\lambda,t_2^0,t_2) +
 \Phi(\lambda,t_2^0,t_2) \sigma_1(t_2)^{-1}(\lambda \sigma_2(t_2) + \gamma(t_2))] e d\lambda = 0,
\end{array} \]
which is true because
\[ \Phi(\lambda,t_2^0,t_2) = \Phi(\lambda,t_2,t_2^0)^{-1} \]
and from (\ref{eq:ODEPhiOUT})
\[ \frac{d}{dt_2} \Phi(\lambda,t_2^0,t_2) +
 \Phi(\lambda,t_2^0,t_2) \sigma_1(t_2)^{-1}(\lambda \sigma_2(t_2) + \gamma(t_2)) = 0. \]
We have thus showed that all the vessel conditions are fullfiled.

\subsection{Regularity}
We shall show first that \textbf{$F(t_2^0,t_2) B(t_2), F^*(t_2,t_2^0) B(t_2)$ are absolutely continuous.}
It is necessary to show that the obtained operators satisfy the conditions of assumption
\ref{assm:Regularity}. So let us check the norms:
\[
F(t_2^0,t_2) B(t_2) =
   \int\limits_{|\lambda| = R} (\lambda I - A_1(t_2^0))^{-1}
   B(t_2^0) \sigma_1(t_2^0) \Phi(\lambda,t_2^0,t_2) d\lambda.
\]
This must be an absolutely continuous function of $t_2$ for all
$t_2^0$. But
\[ \begin{array}{llllll}
\| F(t_2^0,b_i) B(b_i) - F(t_2^0,a_i) B(a_i) \| = \\
~~~~~~~
  \| \int\limits_{|\lambda| = R} (\lambda I - A_1(t_2^0))^{-1}
  B(t_2^0) \sigma_1(t_2^0) [\Phi(\lambda,t_2^0,b_i) - \Phi(\lambda,t_2^0,a_i)] d\lambda \| \leq \\
~~~~~~~
  \int\limits_{|\lambda| = R} \| (\lambda I - A_1(t_2^0))^{-1}
  B(t_2^0) \sigma_1(t_2^0) [\Phi(\lambda,t_2^0,b_i) - \Phi(\lambda,t_2^0,a_i)] \| d\lambda \leq \\
~~~~~~~
  \int\limits_{|\lambda| = R} \| (\lambda I - A_1(t_2^0))^{-1}
  B(t_2^0) \sigma_1(t_2^0)\| \| [\Phi(\lambda,t_2^0,b_i) - \Phi(\lambda,t_2^0,a_i)] \| d\lambda \leq \\
~~~~~~~
  \int\limits_{|\lambda| = R} \| (\lambda I - A_1(t_2^0))^{-1} B(t_2^0) \sigma_1(t_2^0)\|  d\lambda
  \sup\limits_{\lambda} \| \Phi(\lambda,t_2^0,b_i) - \Phi(\lambda,t_2^0,a_i) \|
\end{array} \]
Since $\Phi(\lambda,t_2^0,t_2)$ is an entire function of
$\lambda$ for a fixed $t_2$, and is an absolutely continuous
function of $t_2$ for a fixed $\lambda$, we claim that
\[ \sup\limits_{\lambda} \| \Phi(\lambda,t_2^0,b_i) - \Phi(\lambda,t_2^0,a_i) \| < \epsilon \]
will be as small as desired, provided $|b_i - a_i| < \delta$.
Indeed,
\[ \begin{array}{llllll}
\sup\limits_{\lambda} \| \Phi(\lambda,t_2^0,b_i) - \Phi(\lambda,t_2^0,a_i) \| =
\sup\limits_{\lambda} \| \int\limits_{a_i}^{b_i}
\Phi(\lambda,t_2^0,s) \sigma_1^{-1}(s)\big( \sigma_2(s) \lambda + \gamma(s) \big) ds \| \leq \\
~~~~~~ \leq \sup\limits_{\lambda} \int\limits_{a_i}^{b_i}
\| \Phi(\lambda,t_2^0,s) \sigma_1^{-1}(s)\big( \sigma_2(s) \lambda + \gamma(s) \big) \| ds \leq \\
~~~~~~ \leq \sup\limits_{\lambda} \int\limits_{a_i}^{b_i}
\| \lambda \Phi(\lambda,t_2^0,s) \sigma_1^{-1}(s) \sigma_2(s) \| ds +
\sup\limits_{\lambda} \int\limits_{a_i}^{b_i}
\| \Phi(\lambda,t_2^0,s) \sigma_1^{-1}(s) \gamma(s) \| ds \leq \\
~~~~~~ \leq \sup\limits_{\lambda, s} \| \lambda \Phi(\lambda,t_2^0,s) \sigma_1^{-1}(s)\|
\int\limits_{a_i}^{b_i} \| \sigma_2(s) \| ds + \sup\limits_{\lambda, s}
\| \Phi(\lambda,t_2^0,s) \sigma_1^{-1}(s)\| \int\limits_{a_i}^{b_i} \| \gamma(s) \| ds \leq \\
~~~~~~ \leq K (b_i-a_i) < \epsilon,
\end{array} \]
since $\Phi(\lambda,t_2^0,s) \sigma_1^{-1}(s)$ and
$\lambda\Phi^{in}(\lambda,t_2^0,s) \sigma_1^{-1}(s)$ are
continuous in two variables and are considered on the compact
space $\mathcal O \times[t_2^0, t_2]$ and
$\int\limits_{a_i}^{b_i} \| \sigma_2(s) \|$ and
$\int\limits_{a_i}^{b_i} \| \gamma(s) \| ds $ exists and are as
small as $|b_i -a_i|$.

The other multiplier is constant:
\[ \int\limits_{|\lambda| = R} \| (\lambda I - A_1(t_2^0))^{-1} B(t_2^0) \sigma_1(t_2^0)\|  d\lambda
< \infty \]
and thus the absolute continuity is obtained.

\noindent\textbf{Boundedness of} $\widehat F(t_2,t_2^0)$. We shall use here the Gronwall formula, saying that
\[ \frac{d}{d t_2} \| h_{t_2}\| \leq c(t_2) \| h_{t_2}\| \Longrightarrow
\| h(t_2)\| \leq \exp[\int_{t_2^0}^{t_2} c(y) dy] \|h(t_2^0)\|.
\]
Using the second colligation condition for $h_{t_2} = A_1^n(t_2)
B(t_2) \sigma_1(t_2) e_j \in \mathcal H_{t_2}$: we obtain that
\[ \begin{array}{llllll}
\frac{d}{dt_2} \| h_{t_2} \|_{H_{t_2}} = \| \sigma_2(t_2) B^*(t_2)
h_{t_2}, B^*(t_2) h_{t_2} \|_{\mathcal H_{t_2}} = \langle B(t_2)
\sigma_2(t_2) B^*(t_2) h_{t_2}, h_{t_2} \rangle^{1/2}_{\mathcal
H_{t_2}}
\end{array} \]
and in order to use the Gronwall formula it is enough to prove that
\[ \| B(t_2) \sigma_2(t_2) B^*(t_2) \|_{\mathcal H_{t_2}} \leq c(t_2)
\]
for a continuous $c(t_2)$. We shall do it by the definition:
\[
\| B(t_2) \sigma_2(t_2) B^*(t_2) \|_{\mathcal H_{t_2}}
 = \sup\limits_{h_{t_2} \in \mathcal H_{t_2}} \dfrac{\| B(t_2) \sigma_2(t_2) B^*(t_2) h_{t_2}\|_{\mathcal H_{t_2}}}
       {\| h_{t_2}\|_{\mathcal H_{t_2}}},
\]
but
\[ \begin{array}{llllll}
\| B(t_2) \sigma_2(t_2) B^*(t_2) h_{t_2} \|^2_{\mathcal H_{t_2}} =
\langle B(t_2) \sigma_2(t_2) B^*(t_2) h_{t_2},B(t_2) \sigma_2(t_2) B^*(t_2) h_{t_2} \rangle_{\mathcal H_{t_2}} = \\
\langle B^*(t_2) B(t_2) \sigma_2(t_2) B^*(t_2) h_{t_2}, \sigma_2(t_2) B^*(t_2) h_{t_2} \rangle_{\mathcal H_{t_2}} \leq \\
\| B^*(t_2) B(t_2)\| \langle \sigma_2(t_2) B^*(t_2) h_{t_2}, \sigma_2(t_2) B^*(t_2) h_{t_2} \rangle_{\mathcal H_{t_2}} \leq \\
\| B^*(t_2) B(t_2)\| \| \sigma_2(t_2) \|^2 \| B^*(t_2) h_{t_2}
\|^2_{\mathcal H_{t_2}}.
\end{array} \]
Thus
\[ \begin{array}{llllll}
\sup\limits_{h_{t_2}\in \mathcal H_{t_2}} \dfrac{\| B(t_2)
\sigma_2(t_2) B^*(t_2) h_{t_2} \|_{\mathcal H_{t_2}}}
       {\| h_{t_2}\|_{\mathcal H_{t_2}}} \leq \\
~~~~~ \leq \| B^*(t_2) B(t_2)\|^{1/2} \| \sigma_2(t_2) \|
\sup\limits_{h_{t_2}\in \mathcal H_{t_2}}
 \dfrac{ \langle B^*(t_2) h_{t_2}, B^*(t_2) h_{t_2} \rangle^{1/2}_{\mathcal H_{t_2}}}
 {\| h_{t_2}\|_{\mathcal H_{t_2}}} =\\
~~~~~~~~~~~~~ = \| B^*(t_2) B(t_2)\|^{1/2} \| \sigma_2(t_2) \| \| B(t_2) \|.
\end{array} \]
But we can use here that $\| B(t_2)\| = \| B^*(t_2)\|$. On the
other hand,
\[ \| B(t_2) e\|^2 = \langle B(t_2) e, B(t_2) e \rangle_{\mathcal H_{t_2}} =
\langle B^*(t_2) B(t_2) e, e \rangle_{\mathcal E} \leq \| B^*(t_2)
B(t_2) \| \| e\|^2_{\mathcal E}.
\]
Finally, we conclude that
\[ \sup\limits_{h_{t_2}\in \mathcal H_{t_2}} \dfrac{\| B(t_2) \sigma_2(t_2) B^*(t_2) h_{t_2} \|_{\mathcal H_{t_2}}}
       {\| h_{t_2} \|_{\mathcal H_{t_2}}} \leq
  \| B^*(t_2) B(t_2)\| \| \sigma_2(t_2) \|
\]
but each element $B^*(t_2) B(t_2)$ is actually $S_0(t_2)
\sigma_1(t_2)^{-1}$, where $S_0(t_2)$ is the first coefficient of
the infinity expansion of $S(\lambda,t_2)- I$. Thus $S_0(t_2)$ is a
continuous function. Thus $S_0(t_2) \sigma_1(t_2)^{-1}$ is a
continuous function and finally $B^*(t_2) B(t_2)$ is too. Thus we
have proved that
\[ \begin{array}{llllll}
\frac{d}{d t_2} \| h_{t_2}\|_{\mathcal H_{t_2}} & = \langle
\sigma_2(t_2) B^*(t_2) h_{t_2}, B^*(t_2) h_{t_2} \rangle_{\mathcal
H_{t_2}}
= \langle B(t_2) \sigma_2(t_2) B^*(t_2) h_{t_2},  h_{t_2} \rangle_{\mathcal H_{t_2}} \leq \\
& \leq \| B^*(t_2) B(t_2)\| \| \sigma_2(t_2) \| \| h_{t_2} \|_{\mathcal H_{t_2}}
\end{array} \]
and taking
\[
c(t_2) = \| B^*(t_2) B(t_2)\| \| \sigma_2(t_2) \|
\]
we obtain the desired result.

\part{Triangular forms and multiplicative integrals}
In the fundamental paper of M.S. Brodskii, M. Liv\v sic \cite{bib:SpectrAnal} it is established a very special (and useful) form
for an operator $A\in B(\mathcal H)$ having ''small,, imaginary part, i.e., $\dim(A-A_*)<\infty$. Denoting $\mathcal E = \frac{A-A^*}{i} \mathcal H$
and representing the operator
\[ \frac{A-A^*}{i} f  = J(f,e) e ~~~(f\in\mathcal H, e\in\mathcal E, J = \pm 1),
\]
there is defined characteristic function
\[ w(\lambda) = I - i \langle(A - \lambda I)^{-1}e,e\rangle J.
\]
For example, a simple operator 
\[ A_0 f = \lambda_0 f, ~~(f\in\mathcal H_0, \dim \mathcal H_0 =1, \Im \lambda_0 \neq 0)
\]
has characteristic function $w_0(\lambda) = \frac{\lambda - \bar\lambda_0}{\lambda - \lambda_0}$. This function $w_0(\lambda)$ maps upper half plane
(of $\mathbb C$) conformally to $D = \{z\in\mathbb C\mid |z| \leq 1\}$ or to $\mathbb C\backslash D$, depending on the sign of 
$\operatorname{Im} \lambda_0$.

Second important example is the operator
\[ A_1 f = i \int_0^x f(t) dt, ~~(0\leq x\leq L, f\in L^2(0,L)),
\]
whose spectrum is $\lambda = 0$. This operator has characteristic function
\[ w_1(\lambda) = e^{\frac{iL}{\lambda}}
\]
and has essential singularity at the point $\lambda = 0$. These two operators (or their analogues to our setting) are the building ''blocks,,
for constructing a model for any operator $A_1(t_2)$. Performing cascade connections of operators, one will obtain \cite{bib:NonConsTheory}
an operator, whose characteristic function is multiplication of the original ones. Also the converse holds, finding an invariant subspace
for $A$ and for $A^*$, one can cascadly decompose the operator and represent its characteristic function as a multiplication of characteristic
functions of its cascade ingredients. The highlight of this theory is that finding a maximal chain of invariant subspaces one can prove that
\cite{bib:JordanBrodskii} any operator in this class has a triangular model. In this model $\mathcal H = \mathcal H_1 \oplus \mathcal H_2$,
where
\[ H_1 = l^2 = \{ (c_i) = (c_0,c_2,\ldots) \mid \sum |c_i|^2 < \infty\},~~H_2 = L^2[0,L]
\]
and $A$ is cascade connection of $A_1\in B(\mathcal H_1), A_2\in\mathcal H_2$, defined by
\[ \begin{array}{lll}
 (A_1 c)_i = (\lambda_i c_i + i \sum\limits_{k=i+1}^\infty c_i \Pi_k J \Pi_k^*, \\
 A_2 f(x) = f(x) c(x) + i \int_0^x f(t) \Pi(t) J \Pi^*(x) dt.
\end{array} \]
for some explicitly defined $\Pi_k$, $\Pi(t)$, satisfying
\[ \Pi_k J \Pi^*_k = 2 \Im \lambda_i, ~~\operatorname{tr}(\Pi^*(x)\Pi(x)) =1, 0\leq x \leq L.
\]
and continuous from the right $c(t)$ function. The result of this construction is that the corresponding characteristic function has
the form
\[ w(\lambda) = \prod\limits_{i=0}^n \frac{\lambda-\bar\lambda_i}{\lambda-\lambda_i} \operatorname{exp}(i\stackrel{\leftarrow}\int_0^L \frac{\Pi(t) J \Pi^*(x)}{\lambda - c(x)} dx).
\]
This means that its integral part $w_2(\lambda) = \operatorname{exp}(i\int_0^L \frac{\Pi(t) J \Pi^*(x)}{\lambda - c(x)} dx)$ can be represented
as a function of two variables 
\begin{equation} \label{eq:w2Int}
w_2(\lambda,s) = \operatorname{exp}(i\stackrel{\leftarrow}\int_0^s \frac{\Pi(t) J \Pi^*(x)}{\lambda - c(x)} dx),
\end{equation}
which satisfies
\[ \frac{\partial}{\partial s} w_2(\lambda,s) = i \frac{\Pi(s) J \Pi^*(s)}{\lambda - c(s)} w_2(\lambda,s).
\]

Following these ideas we present first finite cascade connections of vessels (see \cite{bib:NonConsTheory} for
its definition). Let $S(\lambda, t_2) \in \boldsymbol{\mathcal CI}$ be an $n\times n$ matrix-function that maps solutions of the input differential equation with the
spectral parameter $\lambda$
\begin{equation} \label{eq:CompInput}
\sigma_2 \frac{\partial}{\partial t_1}u - \sigma_1 \frac{\partial}{\partial t_2} u + \gamma u = 0,
\end{equation}
to the output differential equation with the same spectral parameter $\lambda$
\begin{equation} \label{eq:CompOutput}
 \sigma_2 \frac{\partial}{\partial t_1} y - \sigma_1 \frac{\partial}{\partial t_2}y +
\gamma_* y = 0.
\end{equation}
This can be written \cite{bib:NonConsTheory} by means of fundamental matrices of the input $\Phi(\lambda,t_2,t_2^0)$ and of the output
$\Phi_*(\lambda,t_2,t_2^0)$ ODEs. Namely:
\begin{equation} \label{eq:SS0}
 S(\lambda, t_2) \Phi(\lambda,t_2,t_2^0) = \Phi_*(\lambda,t_2,t_2^0) S(\lambda, t_2^0).
\end{equation}
The corresponding conservative vessel is a collection of operators and spaces
\[ \mathfrak{V} = (A_1, A_2, B; \sigma_1, \sigma_2, \gamma, \gamma_*; \mathcal{H}, \mathcal{E}), \]
such that $A_1,A_2 \in B(\mathcal{H})\hspace{2mm} \text{and} \hspace{2mm} B \in B(\mathcal{E}, \mathcal H)$,
which are all functions of $t_2$ satisfying the following axioms:
\[ \begin{array}{lllllll}
   \frac{d}{dt_2} A_1 = A_2 A_1 - A_1 A_2 \\
   A_1 + A_1^* + B \sigma_1 B^* = 0 \\
   A_2 + A_2^* + B \sigma_2 B^* = 0 \\
   \frac{d}{dt_2} \big(B \sigma_1\big) - A_2 B \sigma_1 + A_1  B \sigma_2 + B \gamma = 0 \\
   \frac{d}{dt_2} \big(B \sigma_1 \big) + A_2^* B \sigma_1  - A_1 B \sigma_2 - B (\gamma_*^* + \frac{d}{dt_2}\sigma_1) = 0 \\
   \gamma = \sigma_2 B^* B \sigma_1 -\sigma_1 B^* B \sigma_2 + \gamma_*, \\
   \sigma_1 = \sigma_1^*, ~~~\sigma_2 = \sigma_2^*,~~~ \Re\gamma =\Re \gamma_* = \frac{\partial}{\partial t_2} \sigma_1.
   \end{array}
\]
Performing coupling of such vessels, one can obtain a new vessel with transfer function,
which is the product of all the previous \cite{bib:TheoryNonComm}.
The new inner Hilbert space is the direct sum of all inner spaces.
In this manner one obtains the matrix function as a product of simpler ones.

The converse of this process is our main interest. 

\textit{We suppose that the regularity assumptions \ref{assm:Regularity} are slightly more restrictive, namely, that
all the absolutely continuous conditions are replaced by continuous derivatives conditions}.

given a function $S(\lambda,t_2)$
satisfying (\ref{eq:SS0}), can it be written as a finite product
of ``elementary'' vessels (=''bulding,, blocks in Liv\v sic-Brodskii setting)? We specify exactly what elementary means.
\begin{defn} A vessel $\mathfrak{V}$ is \textbf{elementary} if the spectrum of $A_1$ is a single point.
\end{defn}
Suppose that we are given a vessel $\mathfrak{V}$, such that
\[ A_1 = z_1,
\]
where $z_1$ is a point of the spectrum. (We allow repetitions of spectrum points.)

The first simple vessel will be
\[ \mathfrak{V_1} =
(A_1^1, A_2^1, B_1; \sigma_1, \sigma_2, \gamma_1, \gamma_{1*};
\mathcal{H} = \mathbb C, \mathcal{E}) \]
where $B_1$ is a solution of the following differential equation with the spectral parameter $z_1$:
\[ z_1 \sigma_2 u - \sigma_1 \frac{\partial}{\partial t_2} u + \gamma u = 0
\]
and (for $\theta_1(t_2)$ an arbitrary differentiable real valued function)
\[ \begin{array}{llll}
A_1^1 = z_1, A_2^1 = -\frac{B_1\sigma_2B_1^*}{2 B_1\sigma_1B_1^*} + i \theta_1'(t_2), \\
\gamma_1 = \gamma, \\
\gamma_{1*} = \gamma_1 + \sigma_2 B_1^* B_1 \sigma_1 - \sigma_1 B_1^* B_1 \sigma_2.
\end{array}
\]
Notice that the characteristic function of the vessel $\mathfrak{V_1}$ is
\[ S_1(\lambda,t_2) = I + B_1^* (\lambda I + A_1)^{-1} B_1 \sigma_1 =
I + B_1^* (\lambda I + z_1)^{-1} B_1 \sigma_1.
\]
Define next vessel
\[ \mathfrak{V_2} =
(A_1^2, A_2^2, B_2; \sigma_1, \sigma_2, \gamma_2, \gamma_{2*};
\mathcal{H} = \mathbb C, \mathcal{E}), \]
where $B_2$ is a solution of the differential equation with the spectral parameter
$z_2$
\[ z_2 \sigma_2 u - \sigma_1 \frac{\partial}{\partial t_2} u + \gamma_2 u = 0,
\]
and for any real valued differentiable function $\theta_2$,
\[ \begin{array}{llll}
A_1^2 = z_2, A_2^2 = -\frac{B_2\sigma_2B_2* }{2 B_2\sigma_1B_2^* } + i \theta_2'(t_2), \\
\gamma_2 = \gamma_{1*}, \\
\gamma_{2*} = \gamma_2 + \sigma_2 B_2^* B_2 \sigma_1 - \sigma_1 B_2^* B_2 \sigma_2.
\end{array}
\]
Notice that coupling the vessels $\mathfrak{V_1}$ and $\mathfrak{V_2}$
one obtains a vessel $\mathfrak{V_2}\vee \mathfrak{V_1}$ with
the corresponding characteristic function
\[ \begin{array}{llll}
 S_{\mathfrak{V_2}\vee \mathfrak{V_1}}(\lambda,t_2) & = I + \bbmatrix{B_1^* & B_2^*}
\big( \lambda I + \bbmatrix{z_1 & 0 \\ -B_2 \sigma_1 B_1^* & z_2} \big)^{-1}
\bbmatrix{B_1 \\ B_2} \sigma_1 = \\
& = S_{\mathfrak{V_2}} S_{\mathfrak{V_1}}.
\end{array} \]
\begin{defn} \label{def:SimpVes}
Define recursively the following vessels (for $2 \leq i \leq n$):
\[ \mathfrak{V_i} =
(A_1^i, A_2^i, B_i, C_i, I; \sigma_1, \sigma_2, \gamma_i, \gamma_{i*};
\mathcal{H} = \mathbb C, \mathcal{E}) \]
where $B_i$ is a solution of the differential equation with the spectral parameter $z_i$
\[ z_i \sigma_2 u - \sigma_1 \frac{\partial}{\partial t_2} u + \gamma_i u = 0
\]
and
\[ \begin{array}{llll}
A_1^i = z_i, A_2^i = -\frac{B_i\sigma_2B_i* }{2 B_i\sigma_1B_i^* } + i \theta_i'(t_2), \\
\gamma_i = \gamma_{(i-1)*}, \\
\gamma_{i*} = \gamma_i + \sigma_2 B_i^* B_i \sigma_1 - \sigma_1 B_i^* B_i \sigma_2.
\end{array}
\]\end{defn}
\textbf{Remarks}

\noindent\textbf{1. }Notice that there is no additional condition on $B_i$'s,
but in order to perform infinite couplings
at the next stage, we demand that
\[ \operatorname{tr} \big(\sigma_1(t_2) B(t_2,s) B(t_2,s)^* \big) = 1
\]
ensuring control of convergences.

\noindent\textbf{2. }The operator $A_2$ can be defined on the dense set $A_1^n B$ in a very specific way
\begin{equation} \label{eq:A2Dense}
\begin{array}{lll}
A_2 A_1^n B = -\frac{d}{dt_2} A_1 \underbrace{A_1 \ldots A_1}_{n-1}B - A_1 \frac{d}{dt_2} A_1 \underbrace{A_1 \ldots A_1}_{n-2}B-
\underbrace{A_1 \ldots A_1}_{n-1}\frac{d}{dt_2} A_1 B + \\
+ A_1^n B \gamma \sigma_1^{-1} + A_1^{n+1}B\sigma_2\sigma_1^{-1}+A_1^n \frac{d}{dt_2}B.
\end{array}
\end{equation}

\noindent\textbf{3.} The characteristic function of the final vessel is
\[ S_{{\mathfrak{V}_N\vee \ldots\vee\mathfrak{V_1}}(\lambda,t_2)} =
S_{\mathfrak{V}_N} \ldots S_{\mathfrak{V_1}}.
\]

\vspace{5mm}

\noindent Another conclusion from the coupling construction is the following.
\begin{thm} Coupling $N$ simple vessels, defined by definition \ref{def:SimpVes}, the vessel
\[ \mathfrak{V} = (A_1, A_2, B; \sigma_1, \sigma_2, \gamma_0, \gamma_{N*};
\mathcal{H} = \oplus_{i=1}^N\mathbb C, \mathcal{E})
\]
is obtained, whose characteristic function is the multiplication of all simple characteristic functions. The corresponding operators are
\[ \begin{array}{llllllll}
B = \bbmatrix{B_1 \\ B_2\\ \vdots \\ B_n}, \\
A_1 = \bbmatrix{z_1 & 0 & 0 & \cdots & 0 \\
          -B_2 \sigma_1 B_1^* & z_2 & 0 & \cdots & 0 \\
          -B_3 \sigma_1 B_1^* & -B_3 \sigma_1 B_2^* & z_3 & \cdots & 0 \\
          \vdots & \vdots & \vdots & \ddots & \vdots \\
          -B_n \sigma_1 B_1^* & -B_n \sigma_1 B_2^* & -B_n \sigma_1 B_3^* & \cdots & z_n},
\end{array} \]
and $A_2$ can be defined in the following way:
\begin{equation} \label{eq:A2}
A_2 = \bbmatrix{-\frac{B_1\sigma_2B_1* }{2 B_1\sigma_1B_1* } + i \theta_1'(t_2) &  0 & \cdots & 0 \\
       -B_2 \sigma_2 B_1^* & -\frac{B_2\sigma_2B_2* }{2 B_2\sigma_1B_2^* } + i \theta_2'(t_2) &  \cdots & 0 \\
       \vdots & \vdots & \ddots & \vdots \\
       -B_n \sigma_2 B_1^* & -B_n \sigma_2 B_2^* & \cdots &
 -\frac{B_n\sigma_2B_n* }{2 B_n\sigma_1 B_n^* } + i \theta_n'(t_2)}.
\end{equation}
\end{thm}
Notice that formula (\ref{eq:A2Dense}) is equivalent to (\ref{eq:A2}), which is an interesting fact by itself.

\section{Construction of a vessel with the discrete spectral data of $A_1$}
At the first stage, we would like to build a vessel which maps solutions of the input ODE to a certain output ODE with the same
spectral parameter $\lambda$, having operator $\widetilde A_1(t_2)$ with the discrete part of the original $A_1(t_2)$.
Let the initial ODE with the spectral parameter $\lambda$ be of the form
\[ \lambda \sigma_1 u - \sigma_2 \dfrac{\partial}{\partial t_2} u - \gamma_0(t_2) u = 0,
\]
where $\sigma_1, \sigma_2, \gamma_0(t_2) \in B(\mathcal E)$. Suppose that $\lambda^{(h)}$ is the discrete spectrum
that we want to
obtain for the operator $A_1$ in the final construction.
This will be achieved by the coupling of elementary vessels in the following way
\begin{defn} \label{def:SpecDisc}
The \textbf{discrete spectral data} are $\{ \lambda^{(h)} \}$ and the \textbf{discrete auxiliary data} are $\{ b^{(h)} \}$,
where $b^{(h)}$ is a $\dim(\mathcal E)\times 1$ matrix function, satisfying:
\[ \begin{array}{lllll}
        \sigma_2 b^{(h)} \lambda^{(h)} - \sigma_1 \frac{\partial}{\partial t_2} b^{(h)} + \gamma_h(t_2) b^{(h)} = 0, \\
        \sum_{h=1}^\infty b^{(h)*} \sigma_1 b^{(h)} < \infty \text{ for $N = \infty$}.
\end{array} \]
\end{defn}
Notice that each step in the finite coupling is constructed from the previous one using the following recurrence formula
\[ \gamma_{h+1} = \gamma_h + \sigma_2 b^{(h)} b^{(h)*} \sigma_1 - \sigma_1 b^{(h)} b^{(h)*} \sigma_2.
\]
In order to build the final vessel (see theorem \ref{thm:DisVes}) we define the inner Hilbert space $\mathcal H = l^2$
with the usual inner product of $l^2$. Then, for an arbitrary vector $(v_h)_{h=1}^\infty \in l^2$, we define the operators of the final vessel
$A_1 \in B(\mathcal H), B^*\in B(\mathcal H,\mathcal E)$ as follows
\[  \begin{array}{lllll}
(y_h) = A_1 (v_h), \text{ where } y_h = - \sum_{j=1}^{h-1} b^{(h)*} \sigma_1 b^{(j)} v_j + \lambda^{(h)} v_h, \\
B^* (v_h) = \sum_{h=1}^\infty b^{(h)} v_h, \\
\end{array} \]
Notice that from the definition of $B^*$ we can find $B$ using the uniqueness of the adjoint of an operator.
\begin{thm} \label{thm:DisVes} The following set
\[ \mathfrak{V} =
(A_1, A_2, B; \sigma_1, \sigma_2, \gamma_0, \gamma^d;
\mathcal{H} = l^2 , \mathcal{E})
\]
with $A_1,A_2,B$ defined above for $\gamma^d = \lim\limits_{h\rightarrow\infty} \gamma_h$ is a conservative vessel.
\end{thm}
\textbf{Proof:} It is possible to build a vessel in the integrated form with the transfer function $S(\lambda,t_2)$.
This vessel is
\[ \mathcal V = (\widetilde A_1(t_2), \widetilde F, \widetilde B(t_2); \sigma_1,\sigma_2,\gamma; \mathcal H_{t_2}; \mathcal E),
\]
where $\mathcal H_{t_2}$ is a family of Hilbert spaces (for each $t_2$) and $\widetilde F$ is an involution semi group. The
construction of this vessel provided an operator $\Psi = \sqrt{F^* F}$, which acted on a chosen Hilbet space $\mathcal H_{t_2^0}$.
In this manner a unitary operator $U(t_2): \mathcal H_{t_2} \longrightarrow \mathcal H_{\tau_2}$, which is defined by
\[ U(t_2) = \Psi(t_2) F(\tau_2, t_2), \]
enables identification of all the Hilbert spaces $\mathcal H_{t_2}$ with the chosen one $\mathcal H_{t_2^0}$.
This was a crucial step in construction of an equivalent vessel in a differential form with the same transfer function in section
\ref{sec:CDiffVessel}.

Suppose now that we have built for $t_2^0$ a colligation $\mathcal C$ with $A_1=A_1(t_2^0), \mathcal H=\mathcal H_{t_2^0}$,
then from the vessel condition
\[ A_1(t_2) = F(t_2,t_2^0) A_1 F(t_2^0,t_2)
\]
we obtain that the chain $\mathcal G_\alpha$ is mapped by $F(t_2,t_2^0)$ to a maximal chain of $A_1(t_2)$-invariant subspaces in $\mathcal H_{t_2}$.
Applying the unitary operator $U(t_2)$, we obtain another maximal chain $G_\alpha(t_2)$ of $\Psi(t_2)A_1\Psi(t_2)^{-1}$-invariant subspaces in
$\mathcal H$
\[ G_\alpha(t_2) = U(t_2) F(t_2,t_2^0) G_\alpha = \Psi(t_2) G_\alpha.
\]
Now one immediately sees that the one-dimensional space $G_1$ consists of an invariant vector $v_1$
for one value of the spectrum $\lambda_1$. Moreover, $\Psi(t_2) v_1$ is an invariant vector for
$A_1(t_2)$ for the same spectral value. Projecting on the space $G_1(t_2)$ for each $t_2$,
we obtain a vessel with one point spectrum and transfer function $S_1(\lambda,t_2)$.
Defining further $S^1(\lambda,t_2) = S(\lambda,t_2) S_1(\lambda,t_2)^{-1}$, we obtain a vessel
with one discrete point spectrum less by one comparatively to the original vessel. 
I we continue in this manner while all the discrete points disappear, we shall obtain a transfer function
\[
S^c(\lambda,t_2) = S(\lambda,t_2)\prod_{i=1}^N S_i(\lambda,t_2)^{-1}, ~~\text{$N$ is finite of $\infty$}
\]
with a purely continuous spectrum. \qed
\section{Further construction of vessels for the continuous spectral data of $A_1$
from the initial vessel with associated $\gamma^d$}
Suppose now that we have a vessel with purely continuous spectrum. It means that for $t_2^0$ there exists $c(s)$ - continuous from the left 
scalar function and $\beta(t_2^0,s)$ such that
\[ S(\lambda,t_2^0) = \stackrel{\leftarrow}{\int_0^l} \exp (\frac{\beta(t_2^0,t) \beta^*(t_2^0,s)\sigma_1 }{\lambda + c(s)}) ds
\]
and there corresponds a triangular model with $\mathcal H = L_{\sigma_1}^2(0,L;\mathbb C^p)$ and (for each \linebreak $f(t) \in \mathcal H$)
\[ \begin{array}{llll}
A_1 (f(t)) = -\int_0^t \beta(t_2^0,t) \sigma_1 \beta^*(t_2^0,s) f(s) ds - c(t) f(t), \\
B^* (f(t)) = \int_0^L \beta(t_2^0,s) f(s) ds.
\end{array}
\]
Suppose also that we have constructed a conservative vessel 
\[  \mathcal V = (A_1(t_2), F(t_2,t_2^0), B(t_2); \sigma_1,\sigma_2,\gamma, \gamma_*; \mathcal H_{t_2}; \mathcal E), 
\]
with the prescribed transfer function $S(\lambda,t_2)$ and which for $t_2^0$ has the triangular
structure above. It can be done using the general theorem on reconstruction of vessels from their transfer functions.
Then its transfer function may be written in the following form
\[ \begin{array}{llll}
S(\lambda,t_2) = I + B^*(t_2) (\lambda I - A_1(t_2^0)^{-1} B(t_2) \sigma_1 = \\
~~~~ = I + \big( F^{-1}(t_2,t_2^0) B(t_2)\big)^* F^*(t_2,t_2^0)F(t_2,t_2^0) (\lambda I - A_1(t_2^0)^{-1}
F^{-1}(t_2,t_2^0) B(t_2) \sigma_1 
\end{array} \]
Using the vessel conditions, if one denotes $b_0(t_2) = F^{-1}(t_2,t_2^0) B(t_2)$, 
$X^{-1}(t_2) = F^*(t_2,t_2^0)F(t_2,t_2^0)$ then these functions satisfy
\[ \begin{array}{llll}
\frac{\partial}{\partial t_2} \big( b_0(t_2) \sigma_1 \big) + A_1(t_2^0) b_0(t_2) \sigma_2 + b_0(t_2) \gamma = 0, \\
\frac{\partial}{\partial t_2} \big( X^{-1}(t_2) \big) = - X^{-1}(t_2) b_0(t_2) \sigma_2 b^*_0(t_2) X^{-1}(t_2) ,~~~~~ X^{-1}(t_2^0) = I, \\
A_1(t_2^0) X(t_2) + X(t_2) A_1^*(t_2^0) = - b_0(t_2) \sigma_1 b^*_0(t_2)
\end{array} \]
and
\[ S(\lambda,t_2) = I + b_0^*(t_2) X(t_2)^{-1} (\lambda I - A_1(t_2^0)^{-1} b_0(t_2) \sigma_1.
\]
Notice that the third condition is a straightforward conclusion of the first and the second.
The following diagram illustrates the operators and the spaces involved
\[ \begin{array} {lllllll}
\mathcal E & \stackrel{\times b_0(t_2,t)}{\longrightarrow} & \mathcal H_{t_2^0} & \stackrel{X(t_2)}{\longrightarrow} & 
\mathcal H_{t_2^0} & \stackrel{\int_0^L b_0^*(t_2,t) \cdot dt}{\longrightarrow} & \mathcal E\\
\end{array} \]

Performing projection of this vessel for each $s\in[0,L]$ on the invariant subspace $\mathcal G_{s} \subseteq \mathcal H_{t_2^0}$, we shall obtain
from the linkage condition
\[ \gamma(s) = \gamma + \sigma_1 B P_{t_2}(s) B \sigma_2 - \sigma_2 B P_{t_2}(s) B \sigma_1,
\]
where $P_{t_2}(s)$ is the orthogonal projection on the space $F(t_2,t_2^0) \mathcal G_{s}$. Let us show first that $\gamma(s)$ obtained in this way
is a differentiable with respect to $s$ function:
\[ \gamma(s+\Delta s) - \gamma(s) = \sigma_1 B [P_{t_2}(s+\Delta s) - P_{t_2}(s)] B \sigma_2 - 
	\sigma_2 B [P_{t_2}(s+\Delta s) - P_{t_2}(s)] B \sigma_1
\]
\begin{lemma} The following limit exists
\[ \lim\limits_{\Delta s\rightarrow 0} \frac{B^*(t_2)(P_{t_2}(s+\Delta s) - P_{t_2}(s)) B(t_2) \sigma_1}{\Delta s},
\]
for almost all $s$.
\end{lemma}
\textbf{proof:}
Notice that for $t_2^0$ can be evaluated explicitly
\begin{equation} \label{eq:t20beta*beta}
\begin{array}{lll}
\lim\limits_{\Delta s\rightarrow 0} \frac{B^*(t_2^0)(P_0(s+\Delta s) - P_0(s)) B(t_2) \sigma_1}{\Delta s} =
\lim\limits_{\Delta s\rightarrow 0} \frac{\int_s^{s+\Delta s} \beta(t_2^0,y) \beta^*(t_2^0,y) \sigma_1}{\Delta s} =
\beta(t_2^0,s) \beta^*(t_2^0,s) \sigma_1.
\end{array} \end{equation}
Using Dunford Shwartz calculus, we shall obtain that
\[ P_{t_2}(s) = \oint f(\lambda) (\lambda I - F(t_2,t_2^0) P_{t_0}(s) F(t_2^0,t_2) F^*(t_2^0,t_2) P_{t_0}(s) F^*(t_2,t_2^0))^{-1} d\lambda,
\]
where the integral is taken around the spectrum of the idempotent $F P_{t_0}(s) F^{-1}$ and $f(\lambda)$ is an analytic function,
which obtains $1$ for real values of $\lambda\in[1,\| F P_{t_0}(s) F^{-1}\|]$ and $f(0) = 0$. We also obtain that
\begin{multline*}
 P_{t_2}(s) = \oint f(\lambda) (\lambda I - F(t_2,t_2^0) P_{t_0}(s) F(t_2^0,t_2) F^*(t_2^0,t_2) P_{t_0}(s) F^*(t_2,t_2^0))^{-1} d\lambda = \\
 =  \oint f(\lambda) (\lambda I - J(t_2,s))^{-1} d\lambda,
\end{multline*}
where we denoted
\[ J(t_2,s) = F(t_2,t_2^0) P_{t_0}(s) F(t_2^0,t_2) F^*(t_2^0,t_2) P_{t_0}(s) F^*(t_2,t_2^0). \]
Thus we obtain
\begin{multline*}
B^*(t_2)(P_{t_2}(s+\Delta s) - P_{t_2}(s)) B(t_2) \sigma_1 = \\
 = B^*(t_2) \oint f(\lambda) (\lambda I - J(t_2,s+\Delta s) )^{-1} d\lambda B(t_2) \sigma_1 - B^*(t_2) \oint f(\lambda) (\lambda I - J(t_2,s))^{-1}
  d\lambda B(t_2) \sigma_1= \\
 = B^*(t_2) \oint f(\lambda) [(\lambda I - J(t_2,s+\Delta s))^{-1} - (\lambda I - J(t_2,s))^{-1}] B(t_2) \sigma_1 d\lambda =
\end{multline*}
Using here the resolvent formula of Hilbert, we shall obtain
\begin{multline*}
B^*(t_2)(P_{t_2}(s+\Delta s) - P_{t_2}(s)) B(t_2) \sigma_1 = \\
 = B^*(t_2) \oint f(\lambda) (\lambda I - J(t_2,s+\Delta s))^{-1} [J(t_2,s+\Delta s) - J(t_2,s)] (\lambda I - J(t_2,s))^{-1} B(t_2) \sigma_1 d\lambda.
\end{multline*}
Notice further that
\begin{multline*}
F(t_2^0,t_2)[J(t_2,s+\Delta s) - J(t_2,s)]F^*(t_2^0,t_2) = \\
= P_{t_0}(s+\Delta s) F(t_2^0,t_2) F^*(t_2^0,t_2) P_{t_0}(s+\Delta s) - P_{t_0}(s) F(t_2^0,t_2) F^*(t_2^0,t_2) P_{t_0}(s) = \\
= P_{t_0}(s+\Delta s) F(t_2^0,t_2) F^*(t_2^0,t_2) P_{t_0}(s+\Delta s) - P_{t_0}(s+\Delta s) F(t_2^0,t_2) F^*(t_2^0,t_2) P_{t_0}(s) + \\
+ P_{t_0}(s+\Delta s) F(t_2^0,t_2) F^*(t_2^0,t_2) P_{t_0}(s) - P_{t_0}(s) F(t_2^0,t_2) F^*(t_2^0,t_2) P_{t_0}(s) = \\
= P_{t_0}(s+\Delta s) F(t_2^0,t_2) F^*(t_2^0,t_2) [P_{t_0}(s+\Delta s) - P_{t_0}(s)] + \\
+ [P_{t_0}(s+\Delta s)-P_{t_0}(s)] F(t_2^0,t_2) F^*(t_2^0,t_2) P_{t_0}(s).
\end{multline*}
Inserting this expression back we obtain
\begin{multline*}
B^*(t_2)(P_{t_2}(s+\Delta s) - P_{t_2}(s)) B(t_2) \sigma_1 = \\
= B^*(t_2) \oint f(\lambda) [(\lambda I - J(t_2,s+\Delta s))^{-1} - (\lambda I - J(t_2,s))^{-1}] B(t_2) \sigma_1 d\lambda  = \\
= B^*(t_2) \oint f(\lambda) (\lambda I - J(t_2,s+\Delta s))^{-1} F(t_2^0,t_2)[P_{t_0}(s+\Delta s) F(t_2^0,t_2) F^*(t_2^0,t_2) [P_{t_0}(s+\Delta s) - P_{t_0}(s)] + \\
+ [P_{t_0}(s+\Delta s)-P_{t_0}(s)] F(t_2^0,t_2) F^*(t_2^0,t_2) P_{t_0}(s)] (\lambda I - J(t_2,s))^{-1} B(t_2) \sigma_1 d\lambda  = \\
= B^*(t_2) \oint f(\lambda) (\lambda I - J(t_2,s+\Delta s))^{-1} F(t_2^0,t_2) P_{t_0}(s+\Delta s) F(t_2^0,t_2) F^*(t_2^0,t_2)\times \\
		[P_{t_0}(s+\Delta s) - P_{t_0}(s)](\lambda I - J(t_2,s))^{-1} B(t_2) \sigma_1 d\lambda + \\
+ B^*(t_2) \oint f(\lambda) (\lambda I - J(t_2,s+\Delta s))^{-1} F(t_2^0,t_2) [P_{t_0}(s+\Delta s)-P_{t_0}(s)] \times \\
		F(t_2^0,t_2) F^*(t_2^0,t_2) P_{t_0}(s)] (\lambda I - J(t_2,s))^{-1} B(t_2) \sigma_1 d\lambda 
\end{multline*}
Let us concentrate on the first expression
\begin{multline*}
B^*(t_2) \oint f(\lambda) (\lambda I - J(t_2,s+\Delta s))^{-1} F(t_2^0,t_2) P_{t_0}(s+\Delta s) F(t_2^0,t_2) F^*(t_2^0,t_2)\times \\
		[P_{t_0}(s+\Delta s) - P_{t_0}(s)](\lambda I - J(t_2,s))^{-1} B(t_2) \sigma_1 d\lambda = \\
= B^*(t_2) F(t_2,t_2^0) P_{t_0}(s) \oint f(\lambda) (\lambda I - J'(t_2,s+\Delta s))^{-1} \times \\
		[P_{t_0}(s+\Delta s) - P_{t_0}(s)](\lambda I - J'^*(t_2,s))^{-1} F^*(t_2^0,t_2) B(t_2) \sigma_1 d\lambda,
\end{multline*}
where 
\[ J'(t_2,s) =  F(t_2^0,t_2) F^*(t_2^0,t_2) P_{t_0}(s) F^*(t_2,t_2^0) F(t_2,t_2^0) P_{t_0}(s) \in B(\mathcal H_{t_2^0}).
\]
Using local observability of the vessel we conclude that 
\[ B^*(t_2) F(t_2,t_2^0) = \sum s_i B^*(t_2^0) A_1^n(t_2^0) F(t_2,t_2^0) \sum s_i B^*(t_2^0) V_i
\]
and plugging it into the last expression, we shall obtain
\[ \sum s_i s_j^* B^*(t_2^0)  V_i \oint f(\lambda) (\lambda I - J(t_2,s+\Delta s))^{-1} [P_{t_0}(s+\Delta s) - P_{t_0}(s)]
(\lambda I - J'^*(t_2,s))^{-1} V_j^* B(t_2^0) \sigma_1
\]
Notice that for fixed $i,j$ similarly to (\ref{eq:t20beta*beta}) the following limit exists
\begin{multline*} 
\lim\limits_{\Delta s\rightarrow 0} \frac{ B^*(t_2^0)  V_i \oint f(\lambda) (\lambda I - J(t_2,s+\Delta s))^{-1} [P_{t_0}(s+\Delta s) - P_{t_0}(s)]
(\lambda I - J'^*(t_2,s))^{-1} V_j^* B(t_2^0) \sigma_1}{\Delta s} = \\
= \beta^*(t_2^0,s)  V_i \oint f(\lambda) (\lambda I - J(t_2,s))^{-1} (\lambda I - J'^*(t_2,s))^{-1} V_j^* \beta(t_2^0,s) \sigma_1
\end{multline*}
But since $\sum s_i s_j^* B^*(t_2^0)  V_i \oint f(\lambda) (\lambda I - J(t_2,s))^{-1} P_{t_0}(s)
(\lambda I - J'^*(t_2,s))^{-1} V_j^* B(t_2^0) \sigma_1$ and \linebreak
$\sum s_i s_j \beta^*(t_2^0,s)  V_i \oint f(\lambda) (\lambda I - J(t_2,s))^{-1} (\lambda I - J'^*(t_2,s))^{-1} V_j^* \beta(t_2^0,s) \sigma_1$ converge 
uniformly in $s$ to their limits we conclude that the second expression is the derivative of the first and this is what we wanted.
\qed

Let us consider next the output vessel condition for various values of $s$:
\[
\sigma_1 \frac{\partial}{\partial t_2}[B^*(t_2) P_{t_2}(s) F(t_2,t_2^0)] - \sigma_2 B^*(t_2) P_{t_2}(s) F(t_2,t_2^0) A_1(t_2^0)  -
        \gamma_*(s) B^*(t_2) P_{t_2}(s) F(t_2,t_2^0) = 0.
\]
This is an ODE for the function $Y(t_2,s) = B^*(t_2) P_{t_2}(s) F(t_2,t_2^0)$, since its derivative is a continuous function of $t_2$ 
\cite{bib:CoddLev}
\[
\sigma_1 \frac{\partial}{\partial t_2}[Y(t_2,s)] - \sigma_2 Y(t_2,s) A_1(t_2^0)  - \gamma_*(s) Y(t_2,s) = 0
\]
Form the theory of ODE with a parameter $s$ we conclude that since $\gamma_*(s)$ is differentiable in $s$ function, so will be
$Y(t_2,s)$ on an interval $[0,L]$. Moreover, the mixed partial derivatives of the second degree are equal:
\[ \frac{\partial^2}{\partial t_2\partial s} Y(t_2,s) = \frac{\partial^2}{\partial s\partial t_2} Y(t_2,s).
\]
Then for the transfer function of the compressed vessel, we obtain
\[ \begin{array}{lll}
S(\lambda,s,t_2)  = \\
= I - B^*(t_2) P_{t_2}(s) F(t_2,t_2^0) (\lambda I - A_1(t_2^0))^{-1} F(t_2^0,t_2) F^*(t_2^0,t_2) F^*(t_2,t_2^0) P_{t_2}(s) B(t_2) \sigma_1 = \\
= I - Y(t_2,s) (\lambda I - A_1(t_2^0))^{-1} F(t_2^0,t_2) F^*(t_2^0,t_2) Y^*(t_2,s) \sigma_1
\end{array} \]
and equality of mixed partial derivatives is immediate.

Let us write down the formula for the logarithmic derivative (w.r.t. to $s$) of $S(\lambda,s,t_2)$ for a fixed $t_2$:
\begin{multline*}
\frac{\partial S(\lambda,s,t_2)}{\partial s}S(\lambda,s,t_2) = 
\lim\limits_{\Delta s\rightarrow 0} \frac{S(\lambda,s+\Delta s,t_2) - S(\lambda,s,t_2)}{\Delta s} S^{-1}(\lambda,s,t_2) = \\
= \lim\limits_{\Delta s\rightarrow 0} \frac{S(\lambda,s+\Delta s,t_2) S^{-1}(\lambda,s,t_2) - I}{\Delta s}.
\end{multline*}
Let us concentrate on the expression $\frac{S(\lambda,s+\Delta s,t_2) S^{-1}(\lambda,s,t_2) - I}{\Delta s}$. From the vessel definitions
it follows that $S(\lambda,s+\Delta s,t_2) S^{-1}(\lambda,s,t_2)$ is the transfer function of the vessel, which is obtained as the projection of
the initial vessel on the invariant subspaces $P_{t_2}(s+\Delta s) - P_{t_2}(s)$, i.e., it is of the form
\[ I - B^*(t_2)(P_{t_2}(s+\Delta s) - P_{t_2}(s)) (\lambda I - A_1'(t_2))^{-1}
(P_{t_2}(s+\Delta s) - P_{t_2}(s)) B(t_2) \sigma_1,
\]
for $A_1'(t_2) = F(t_2,t_2^0)(P_0(s+\Delta s) - P_0(s)) A_1(t_2^0) (P_0(s+\Delta s) - P_0(s)) F(t_2^0,t_2)$. We want to show that 
$\frac{\partial S(\lambda,s,t_2)}{\partial s}S(\lambda,s,t_2)$
is of the form $\frac{\beta(t_2,s) \alpha(t_2,s)}{\lambda + c(s)}$ for $\Delta s \rightarrow 0$, so multiplying the limit expression by
$\lambda + c(s)$, we shall obtain
\[ \begin{array}{llll}
- (\lambda + c(s)) \frac{S(\lambda,s+\Delta s,t_2) S^{-1}(\lambda,s,t_2) - I}{\Delta s} = \\
= (\lambda + c(s)) \frac{B^*(t_2)(P_{t_2}(s+\Delta s) - P_{t_2}(s)) 
		(\lambda I - A_1'(t_2))^{-1} (P_{t_2}(s+\Delta s) - P_{t_2}(s)) B(t_2) \sigma_1}{\Delta s} = \\
= \frac{B^*(t_2)(P_{t_2}(s+\Delta s) - P_{t_2}(s)) (\lambda - A_1'(t_2) + A_1'(t_2) + c(s))
		(\lambda I - A_1'(t_2))^{-1} (P_{t_2}(s+\Delta s) - P_{t_2}(s)) B(t_2) \sigma_1}{\Delta s} = \\
= \frac{B^*(t_2)(P_{t_2}(s+\Delta s) - P_{t_2}(s)) B(t_2) \sigma_1}{\Delta s} + \\
~~~~~~~~~~~~~ + \frac{B^*(t_2)(P_{t_2}(s+\Delta s) - P_{t_2}(s)) (A_1'(t_2) + c(s))
		(\lambda I - A_1'(t_2))^{-1} (P_{t_2}(s+\Delta s) - P_{t_2}(s)) B(t_2) \sigma_1}{\Delta s} = \\
\end{array} \]
and consequently, we need the following lemma
\begin{lemma} For a continuous function $f \in L_{\sigma_1}^2(0,L;\mathbb C^p)$ the following holds
\[ \lim\limits_{\Delta s\rightarrow 0}
\frac{B^*(t_2^0) (A_1(t_2^0) + c(s))(P_0(s+\Delta s) - P_0(s)) f }{\Delta s} = 0.
\]
\end{lemma}
\textbf{Proof:} Let us evaluate
\begin{multline*}
B^*(t_2^0) (A_1(t_2^0) + c(s))(P_0(s+\Delta s) - P_0(s)) f = \\
= \int_s^{s+\Delta s} \beta(t_2^0,t) \big( [c(s)-c(t)] f(t) + \int_s^t \beta(t_2^0,t) \sigma_1 \beta^*(t_2^0,y) f(y) dy \big) dt
\end{multline*}
and using the fact that for any continuous function $f(t)$, $\lim\limits_{\Delta s \rightarrow 0} \frac{\int_s^{s+\Delta s} f(t) dt}{\Delta s} = f(s)$
we obtain that
\begin{multline*}
\lim\limits_{\Delta s \rightarrow 0}\frac{B^*(t_2^0) (A_1(t_2^0) + c(s))(P_0(s+\Delta s) - P_0(s)) f}{\Delta s} = \\
= \lim\limits_{\Delta s \rightarrow 0} \frac{\int_s^{s+\Delta s} \beta(t_2^0,t) \big( [c(s)-c(t)] f(t) + \int_s^t \beta(t_2^0,t) \sigma_1 \beta^*(t_2^0,y) f(y) dy \big) dt}{\Delta s} = \\
= \beta(t_2^0,s) \big( [c(s)-c(s)] f(s) + \int_s^s \beta(t_2^0,s) \sigma_1 \beta^*(t_2^0,y) f(y) dy \big) = 0,
\end{multline*}
as desired.
\qed

Since we have local observability 
\[ \bigvee_{n\geq 0,e \in \mathcal E} A_1^n(t_2^0)\widetilde B(t_2^0)e = 
\bigvee_{n\geq 0,e \in \mathcal E} F(t_2^0,t_2) A_1^n(t_2)\widetilde B(t_2)e = 
\mathcal H_{t_2^0},
\]
in order to have
\[ \lim\limits_{\Delta s\rightarrow 0}
\frac{B^*(t_2) F(t_2,t_2^0) \big(\lambda I - A_1(t_2^0))^{-1}(A_1(t_2^0) + c(s))(P_0(s+\Delta s) - P_0(s)) F(t_2^0,t_2)B(t_2) \sigma_1}{\Delta s} = 0
\]
it is enough to show that for each $f \in L_{\sigma_1}^2(0,L;\mathbb C^p)$
\[ \lim\limits_{\Delta s\rightarrow 0}
\frac{B^*(t_2^0) (A_1(t_2^0) + c(s))(P_0(s+\Delta s) - P_0(s)) f}{\Delta s} = 0,
\]
which follows from the lemma.

So, $S(\lambda,t_2)$ has a realization in the form of multiplicative integral:
\[ S(\lambda,t_2) = \stackrel{\leftarrow}{\int_0^L} exp(\frac{K(t_2,y)}{\lambda+c(y)}) dy, 
\]
which means that it is possible to define
\[ 
S(\lambda,t_2,s) = \stackrel{\leftarrow}{\int_0^s} exp(\frac{K(t_2,y)}{\lambda+c(y)}) dy
\]
such that
\[ \frac{\partial}{\partial s} S(\lambda,t_2,s) = \frac{K(t_2,s)}{\lambda + c(s)} S(\lambda,t_2,s).
\]
Then using the equality 
$\frac{\partial^2}{\partial s\partial t_2} S(\lambda,t_2,s) = \frac{\partial^2}{\partial t_2\partial s} S(\lambda,t_2,s)$, we shall obtain
\begin{multline*} 
\frac{\partial^2}{\partial s\partial t_2} S(\lambda,t_2,s) =
\frac{\partial}{\partial s} [\sigma_1^{-1}(\lambda \sigma_2 + \gamma(s)) S(\lambda,t_2,s) - S(\lambda,t_2,s) \sigma_1^{-1}(\lambda \sigma_2 + \gamma)]=\\
= \sigma_1^{-1} \frac{d}{ds}\gamma(s) S(\lambda,t_2,s) + 
\sigma_1^{-1}(\lambda \sigma_2 + \gamma(s))  \frac{K(t_2,s)}{\lambda + c(s)}S(\lambda,t_2,s) -  \frac{K(t_2,s)}{\lambda + c(s)}S(\lambda,t_2,s) \sigma_1^{-1}(\lambda \sigma_2 + \gamma), \\
\frac{\partial^2}{\partial t_2\partial s} S(\lambda,t_2,s) = \frac{\partial}{\partial t_2} \frac{K(t_2,s)}{\lambda + c(s)} S(\lambda,t_2,s) =
\frac{\frac{\partial}{\partial t_2}K(t_2,s)}{\lambda + c(s)} S(\lambda,t_2,s) + \frac{K(t_2,s)}{\lambda + c(s)} \frac{\partial}{\partial t_2}S(\lambda,t_2,s) = \\
= \frac{\frac{\partial}{\partial t_2}K(t_2,s)}{\lambda + c(s)} S(\lambda,t_2,s) + \frac{K(t_2,s)}{\lambda + c(s)} [\sigma_1^{-1}(\lambda \sigma_2 + \gamma(s)) S(\lambda,t_2,s) - S(\lambda,t_2,s) \sigma_1^{-1}(\lambda \sigma_2 + \gamma)],
\end{multline*}
or after cancellations and multiplying on $S^{-1}(\lambda,t_2,s)$ on the right ($K=K(t_2,s)$)
\[
\sigma_1^{-1} \frac{d}{ds}\gamma(s) + \sigma_1^{-1}(\lambda \sigma_2 + \gamma(s))  \frac{K}{\lambda + c(s)}
= \frac{\frac{\partial}{\partial t_2}K}{\lambda + c(s)} + \frac{K}{\lambda + c(s)} \sigma_1^{-1}(\lambda \sigma_2 + \gamma(s)).
\]
Let us multiply this expression by $\lambda + c(s)$ to obtain
\[
(\lambda+c(s))\sigma_1^{-1} \frac{d}{ds}\gamma(s) + \sigma_1^{-1}(\lambda \sigma_2 + \gamma(s))  K
= \frac{\partial}{\partial t_2}K + K \sigma_1^{-1}(\lambda \sigma_2 + \gamma(s)).
\]
Substituting here $\lambda = \lambda + c(s) - c(s)$ and rearranging, we shall obtain
\begin{multline*} \frac{\partial}{\partial t_2} K + K \sigma_1^{-1}(-c(s) \sigma_2 + \gamma(s)) = \\
\sigma_1^{-1}(-c(s) \sigma_2 + \gamma(s)) K + (\lambda+c(s))(\sigma_1^{-1} \frac{d}{ds} \gamma(s) + \sigma_1^{-1}\sigma_2 K - K \sigma_1^{-1}\sigma_2)).
\end{multline*}
Consequently,
\begin{eqnarray}
\label{eq:DGamma}\sigma_1^{-1} \frac{d}{ds} \gamma(s) + \sigma_1^{-1}\sigma_2 K - K \sigma_1^{-1}\sigma_2 = 0, \\
\label{eq:DK} 
\begin{array}{lll} \frac{\partial}{\partial t_2} K = \sigma_1^{-1}(-c(s) \sigma_2 + \gamma(s)) K - K \sigma_1^{-1}(-c(s) \sigma_2 + \gamma(s)).
\end{array}
\end{eqnarray}
Form here we conclude that
\[ K(t_2,s) = \beta(t_2,s) \alpha(t_2,s),
\]
where $\beta(t_2,s), \alpha(t_2,s)$ satisfy the following differential equations:
\[ \begin{array}{lll}
\frac{\partial}{\partial t_2} \beta(t_2,s) = \sigma_1^{-1}[-c(s)\sigma_2 + \gamma(s)] \beta(t_2,s),   & \beta(t_2^0,s)\text{ from above} ,\\
\frac{\partial}{\partial t_2} \alpha(t_2,s) = \alpha(t_2,s) \sigma_1^{-1} [c(s)\sigma_2 - \gamma(s)], & \alpha(t_2^0,s) = \beta^*(t_2^0,s)\\
\end{array} \]
and
\begin{equation} \label{eq:GamaBetAlfa}
 \frac{d}{ds} \gamma(s) = \sigma_1 \beta(t_2,s) \alpha(t_2,s) \sigma_1^{-1} \sigma_2 - \sigma_2 \beta(t_2,s) \alpha(t_2,s).
\end{equation}
Notice also that from the properties of the fundamental matrices, for 
\[ \beta(t_2,s) = \Phi(t_2,c(s)) \beta(t_2^0,s)
\]
we obtain that
\[ \begin{array}{lll}
\alpha(t_2,s) & = \alpha(t_2^0,s) \Phi^{-1}(t_2,c(s)) = \alpha(t_2^0,s) \sigma_1^{-1} \Phi^*(t_2,-c(s)) \sigma_1 = \\
		& = \beta^*(t_2^0,s)\sigma_1 \sigma_1^{-1} \Phi^*(t_2,-c(s)) \sigma_1 = \beta^*(t_2,s) \sigma_1,
\end{array} \]
which shows that
\[ \frac{\partial}{\partial s} S(\lambda,t_2,s) = \frac{\beta(t_2,s)\beta^*(t_2,s) \sigma_1}{\lambda + c(s)} S(\lambda,t_2,s).
\]
We conclude that $S(\lambda,t_2,s)$ has the multiplicative structure similar to (\ref{eq:w2Int}), which was our goal.
\part{Zero/pole interpolation}
Recall that in \cite{bib:NonConsTheory} (theorem 8.1) we have proved an important realization:
\begin{thm} \label{tm:NonConsRealiz}
Suppose that $S(\lambda, t_2)$ intertwins sollutions of ODEs with a spectral parameter $\lambda$. 
Then there exists vessel $\mathfrak{DV}$ in the differential form
\[ \mathfrak{DV} = (A_1(t_2), A_2(t_2), B(t_2), C(t_2), D(t_2), \widetilde D;
        \sigma_1, \sigma_2, \gamma, \sigma_{1*}, \sigma_{2*}, \gamma_*;
                \mathcal{H}, \mathcal{E}, \mathcal{E}_*, \widetilde{\mathcal E}, \widetilde{\mathcal E}_*),
\]
with this transfer function and for which
\begin{eqnarray}
C(t_2) = \oint\limits_{Spec A_1} \Phi_*(\lambda,t_2,t_2^0) C_0(\lambda I - A_1)^{-1} d\lambda \\
B(t_2) = \oint\limits_{Spec A_1} (\lambda I - A_1)^{-1} B_0 \sigma_1(t_2) \Phi^* (-\bar\lambda,t_2,t_2^0) d\lambda \\
D(t_2) = S(\infty, t_2)
\end{eqnarray}
and
\[ S(\lambda, t_2) = D(t_2) + C(t_2) (\lambda I - A_1)^{-1} B(t_2) \sigma_1. \]
\end{thm}
We are going to present here a different proof of this theorem, using the technique of zero/pole interpolation.
\section{Zero-Pole Interpolation of $S(\lambda, t_2)$}
Suppose that we want to solve a zero-pole interpolation problem.
Following \cite{bib:Inter}, in order to have (for each $t_2$)
zero $(C(t_2), A_\pi)$ - pole $(A_\xi, B(t_2)\sigma_1)$ data,
sufficient for reconstructing $S(\lambda, t_2)$ up to similarity,
there must exist a solution $X(t_2)$ of Sylvester equation
\begin{equation} \label{eq:Sylvester}
 X(t_2) A_\pi - A_\xi X(t_2) = B(t_2) \sigma_1 C(t_2).
\end{equation}
From this equation we obtain that
\[ X^{-1} A_\xi = A_\pi X^{-1} - X^{-1} B(t_2) \sigma_1 C(t_2) X^{-1}
\]
In this case the matrix $X(t_2)$ is called the \textit{null-pole
coupling matrix} for $S(\lambda,t_2)$.

\noindent \textbf{Remarks:} If one starts from a matrix
$S(\lambda,t_2)$ with the usual properties, then one can
explicitly evaluate zero-pole data $(C(t_2), A_\pi), (A_\xi,
B(t_2)\sigma_{1*}), X(t_2)$ by definition. Then differentiating the
Sylvester equation (\ref{eq:Sylvester}), one obtains:
\[ X' A_\pi - A_\xi X' = (-A_\xi B(t_2) \sigma_{2*} - B(t_2) \gamma_*)
C(t_2) + B(t_2) \sigma_{1*} \sigma_{1*}^{-1}(\sigma_{2*} C(t_2)
A_\pi + \gamma_* C(t_2))
\]
or
\begin{equation}\label{eq:ConSylv} (X' - B(t_2) \sigma_{2*} C(t_2)) A_\pi - A_\xi (X' - B(t_2) \sigma_{2*}
C(t_2))= 0.
\end{equation}
\begin{enumerate}
    \item Notice that if additionally, the spectrum of $A_\pi$ is disjoint from the
    spectrum of $A_\xi$, then from the uniqueness of solution of (\ref{eq:ConSylv}) we obtain that
    \begin{equation} \label{eq:SylvD} X' = B(t_2) \sigma_{2*} C(t_2).
    \end{equation}
    In other words, $X(t_2)$ will additionally satisfy
    differential equation (\ref{eq:SylvD}).
    \item Suppose that $X(t_2)$ satisfies differential
    equation (\ref{eq:SylvD}) and for $t_2^0$ the algebraic
    Sylvester equation (\ref{eq:Sylvester}) is satisfied:
    \[ X(t_2^0) A_\pi - A_\xi X(t_2^0) = B(t_2^0) \sigma_1 C(t_2^0),  \]
    then solving differential equation (\ref{eq:SylvD}) with
    this initial condition $X(t_2^0)$, we obtain that $X(t_2)$
    satisfies the algebraic differential equation
    (\ref{eq:Sylvester}) for each $t_2$.
    \item Notice that invertibility of the matrix $X(t_2)$ is not
    globally promised. If $\det X(t_2^0) \neq 0$, then there could be values of $t_2$ for them
    $\det X(t_2) = 0$.
\end{enumerate}

\begin{lemma} Suppose that $S(\lambda, t_2)$ maps solutions as above. Suppose also
$(C(t_2), A_\pi)$ is the right pole data. Then there exists a unitary matrix $U$, such that
$\tilde C(t_2) = C(t_2) U$ satisfies the output differential equation with the spectral matrix-parameter $A_1 = U^* A_\pi U$:
\[ \tilde C(t_2)' \sigma_1 - A_1 \tilde C(t_2)\sigma_2 + \tilde C(t_2) \gamma_* = 0 \]
\end{lemma}
\textbf{Proof:} In order to prove it, one has to use theorem \ref{tm:NonConsRealiz}, where
one builds a minimal realization of $S(\lambda, t_2)$
\[ S(\lambda,t_2) = D(t_2) + \tilde C(t_2) (\lambda I - A_1)^{-1} \tilde B(t_2) \sigma_1
\]
with
\begin{eqnarray}
\tilde C(t_2) = \oint\limits_{Spec A_1} \Phi_*(\lambda,t_2,t_2^0) C_0(\lambda I - A_1)^{-1} d\lambda \\
\tilde B(t_2) = \oint\limits_{Spec A_1} (\lambda I - A_1)^{-1} B_0 \Phi^* (\lambda,t_2,t_2^0) d\lambda \\
D(t_2) = S(\infty, t_2)
\end{eqnarray}
for a fixed realization at $t_2^0$, $S(\lambda,t_2^0) = D_0 + C_0 (\lambda I - A_1)^{-1} B_0 \sigma_1$.
The left pole data, which can be easily read from a minimal realization, is $\tilde C(t_2), A_1$. On the other hand,
we are given the left pole data $C(t_2), A_\pi$. Since all minimal realizations are unitary equivalent, there is
a unitary matrix such that $\tilde C(t_2) = C(t_2)  U(t_2), A_1 = U(t_2) A_\pi U^*(t_2)$. It is immediate from
lemma 8.2 \cite{bib:NonConsTheory}, that $\tilde C(t_2)$ satisfies the output differential equation with the spectral matrix-parameter $A_1$.
\qed

In this manner we see that it is possible to obtain the right pole data in a "convenient" form (namely satisfying the differential
equation with the spectral matrix-parameter). Notice that the same considerations hold for the left null data, too (considering the inverse
matrix). In this case the corresponding differential equation will be adjoint output, because of the formula (\ref{eq:SS*})
connecting $S(\lambda,t_2)$ and $S_*(\lambda,t_2)$. Then
\begin{lemma} Suppose that $S(\lambda,t_2)$ is as above.
If we are given null-pole data $(C(t_2),A_\pi; B(t_2), A_\xi)$ with $C(t_2)$ and $B(t_2)$, satisfying
output and the adjoint output equations with the spectral matrix-parameters $A_\pi$ and $A_\xi$, respectively, then the null-pole
coupling matrix $X(t_2)$ satisfies the differential equation (\ref{eq:SylvD}).
\end{lemma}
\textbf{Proof:} Recall that the null-pole coupling matrix $X(t_2)$ is defined as one satisfying the Sylvester equation
(\ref{eq:Sylvester}) and such that
\[ S(\lambda,t_2) = I + C(t_2) (\lambda I - A_\pi)^{-1} X^{-1} B(t_2) \sigma_1.
\]
Let us remember now equation (\ref{eq:DS}):
\[ \frac{\partial}{\partial t_2} S(\lambda, t_2) = \sigma_{1*}^{-1} (\sigma_{2*} \lambda + \gamma_*) S(\lambda,t_2) -
S(\lambda,t_2) \sigma_1^{-1} (\sigma_2 \lambda + \gamma).
\]
Substituting here the differential equations:
\[ \begin{array} {lllll}
\sigma_1 C(t_2)' = \sigma_2 C(t_2) A_\pi + \gamma_* C(t_2) \\
B(t_2)' \sigma_1 = -A_\xi B(t_2) \sigma_2 - B(t_2) \gamma_*
\end{array} \]
We obtain that
\[ \begin{array} {lllll}
\frac{\partial}{\partial t_2} S(\lambda, t_2) = \frac{\partial}{\partial t_2} C(t_2) (\lambda I - A_\pi)^{-1} X^{-1} B(t_2) \sigma_1 = \\
\sigma_1^{-1} (\sigma_2 C(t_2) A_\pi + \gamma_* C(t_2)) (\lambda I - A_\pi)^{-1} X^{-1} B(t_2) \sigma_1 +
C(t_2) (\lambda I - A_\pi)^{-1} (X^{-1})' B(t_2) \sigma_1 + \\
~~~~~~~~+ C(t_2) (\lambda I - A_\pi)^{-1} X^{-1} (-A_\xi B(t_2) \sigma_2 - B(t_2) \gamma_*).
\end{array} \]
Then equation (\ref{eq:DS}) becomes
\[ \begin{array} {lllll}
\sigma_1^{-1} (\sigma_2 C(t_2) A_\pi + \gamma_* C(t_2)) (\lambda I - A_\pi)^{-1} X^{-1} B(t_2) \sigma_1 +
C(t_2) (\lambda I - A_\pi)^{-1} (X^{-1})' B(t_2) \sigma_1 + \\
~~~~~~~~+ C(t_2) (\lambda I - A_\pi)^{-1} X^{-1} (-A_\xi B(t_2) \sigma_2 - B(t_2) \gamma_*) = \\
\sigma_1^{-1} (\lambda \sigma_2 + \gamma_* ) \big(I + C(t_2) (\lambda I - A_\pi)^{-1} X^{-1} B(t_2) \sigma_1\big) -
\big( I + C(t_2) (\lambda I - A_\pi)^{-1} X^{-1} B(t_2) \sigma_1 \big)  \sigma_1^{-1} (\sigma_2 \lambda + \gamma)
\end{array} \]
After some cancellations (using Sylvester equation \ref{eq:Sylvester}) it means that
\[ \begin{array} {lllll}
-\sigma_1^{-1} \sigma_2 C(t_2)  X^{-1} B(t_2) \sigma_1 + C(t_2) (\lambda I - A_\pi)^{-1} (X^{-1})' B(t_2) \sigma_1 +
C(t_2) X(t_2)^{-1} B(t_2) \sigma_2 + \\
+ C(t_2) (\lambda I - A_\pi)^{-1} X(t_2)^{-1} B(t_2) (\gamma - \gamma_*)C(t_2) +
 (\lambda I - A_\pi)^{-1} X^{-1} B(t_2) \sigma_1 C(t_2) X^{-1} B(t_2) \sigma_2 = \\
 ~~~~~~~~~= \sigma_1^{-1} (\gamma_* - \gamma)
\end{array} \]
or using the linkage condition (\ref{eq:LinkCond})
\[ \begin{array} {lllll}
C(t_2) (\lambda I - A_\pi)^{-1} (X^{-1})' B(t_2) \sigma_1 +
C(t_2) (\lambda I - A_\pi)^{-1} X^{-1} B(t_2) \sigma_1 (\gamma-\gamma_* + C(t_2) X^{-1} B(t_2) \sigma_2 )= 0.
\end{array} \]
Finally, using linkage condition (\ref{eq:OutCC}) again,
\[ \begin{array} {lllll}
C(t_2) (\lambda I - A_\pi)^{-1} [ (X^{-1})' + X^{-1} B(t_2) \sigma_2 C(t_2) X^{-1}] B(t_2) \sigma_1 = 0.
\end{array} \]
Because the realization is minimal, it means that
\[ [ (X^{-1})' + X^{-1} B(t_2) \sigma_2 C(t_2) X^{-1}] B(t_2) \sigma_1 = 0. \]
Rewrite it as
\[ [ X' - B(t_2) \sigma_2 C(t_2)] X^{-1}(t_2) B(t_2) \sigma_1 = 0 \]
and multiply by $(\lambda I - A_\pi)^{-1}$ on the left
\[ (\lambda I - A_\pi)^{-1} [ X' - B(t_2) \sigma_2 C(t_2)] X^{-1}(t_2) B(t_2) \sigma_1 = 0. \]
Using equation (\ref{eq:ConSylv}) we obtain that
\begin{eqnarray*} (\lambda I - A_\pi)^{-1} [ X' - B(t_2) \sigma_2 C(t_2)] X^{-1}(t_2) B(t_2) \sigma_1 =   \\
 = [ X' - B(t_2) \sigma_2 C(t_2)] (\lambda I - A_\xi)^{-1} X^{-1}(t_2) B(t_2) \sigma_1 = 0.
\end{eqnarray*}
Again, since the realization is minimal it means that
\[ (X^{-1})' + X^{-1} B(t_2) \sigma_2 C(t_2) X^{-1} = 0,
\]
or in other words (\ref{eq:SylvD}) holds.
\qed
\begin{thm} \label{thm:uniqueRealiz}
Suppose that $(C(t_2),A_\pi), (A_\xi, B(t_2)), X(t_2)$
are the zero-pole data, with invertible coupling matrix $X(t_2)$ on
the interval $\mathcal I$. Suppose also that additionally to
(\ref{eq:Sylvester}), $X(t_2)$ satisfies (\ref{eq:SylvD}). Then
there exists a unique matrix function $S(\lambda, t_2)$, which maps
solutions of (\ref{eq:InCC}) with spectral parameter
$\lambda$ to solutions of (\ref{eq:OutCC}) with the same
spectral parameter and which is identity at infinity.
\end{thm}
\textbf{Proof:} Notice first that since we have normalized the desired function as $D(t_2) = I$ at
infinity, in order to have linkage conditions, we have to consider the case, when
\[ \sigma_{1*} = \sigma_1, ~~~\sigma_{2*} = \sigma_2.
\]
For each $t_2$ as in \cite{bib:Inter} one defines
\[ \widetilde B(t_2) = X(t_2)^{-1} B(t_2).
\]
Then let us check that the resulting matrix function
\[ S(\lambda, t_2) = I + C(t_2) (\lambda I - A_\pi)^{-1} \widetilde B(t_2) \sigma_1
\]
satisfies all the requirements. In order to do it, we first evaluate the differentiation of
$\widetilde B(t_2)$:
\[ \begin{array}{llll}
\widetilde B(t_2)' & = \big(X(t_2)^{-1}\big)' B(t_2) + X(t_2)^{-1} B(t_2)'= \\
& = -\widetilde B(t_2) \sigma_2 C(t_2) \widetilde B(t_2) - A_\pi \widetilde B(t_2) \sigma_2 \sigma_1^{-1} +
\widetilde B(t_2) \sigma_1 C(t_2) \widetilde B(t_2) \sigma_2 \sigma_1^{-1} - \widetilde B(t_2) \gamma_* \sigma_1^{-1}.
\end{array} \]
Then
\[ \begin{array}{lllllll}
\frac{\partial}{\partial t_2} S(\lambda, t_2) = \sigma_1^{-1} (\sigma_2 C(t_2) A_\pi + \gamma_* C(t_2))
    (\lambda I - A_\pi)^{-1} \widetilde B(t_2) \sigma_1 + \\
~~~~~ + C(t_2) (\lambda I - A_\pi)^{-1} \\
~~~~~~~~~~ \big( -\widetilde B(t_2) \sigma_2 C(t_2) \widetilde B(t_2) - A_\pi \widetilde B(t_2) \sigma_2 \sigma_1^{-1} +
\widetilde B(t_2) \sigma_1 C(t_2) \widetilde B(t_2) \sigma_2 \sigma_1^{-1} - \widetilde B(t_2) \gamma_* \sigma_1^{-1} \big)
\sigma_1  =  \\
= \sigma_1^{-1}(\sigma_2 \lambda + \gamma_*)C(t_2) (\lambda I - A_\pi)^{-1} \widetilde B(t_2) \sigma_1 -
    \sigma_1^{-1} \sigma_2 C(t_2) \widetilde B(t_2) \sigma_1 - C(t_2) \widetilde B(t_2) \sigma_2 + \\
+ C(t_2) (\lambda I - A_\pi)^{-1} \widetilde B(t_2)
\big( - \sigma_2 C(t_2) \widetilde B(t_2) \sigma_1 +
   \sigma_1 C(t_2) \widetilde B(t_2) \sigma_2 - (\lambda \sigma_2 + \gamma_*) \big).
\end{array} \]
Consequently, defining $\gamma$ from the linkage condition (\ref{eq:LinkCond})
\[ \gamma = \sigma_2 C(t_2) \widetilde B(t_2) \sigma_1 - \sigma_1 C(t_2) \widetilde B(t_2) \sigma_2 + \gamma_*
\]
one obtains that
\[ \begin{array}{lllll}
\frac{\partial}{\partial t_2} S(\lambda, t_2) - \sigma_1^{-1} (\sigma_2\lambda + \gamma_*) S(\lambda,t_2) +
    S(\lambda,t_2) \sigma_1^{-1} (\sigma_2\lambda + \gamma) = \\
= \sigma_1^{-1}(\sigma_2 \lambda + \gamma_*)C(t_2) (\lambda I - A_\pi)^{-1} \widetilde B(t_2) \sigma_1 -
    \sigma_1^{-1} \sigma_2 C(t_2) \widetilde B(t_2)\sigma_1 + C(t_2) \widetilde B(t_2) \sigma_2 + \\
+ C(t_2) (\lambda I - A_\pi)^{-1} \widetilde B(t_2)
\big( - \sigma_1^{-1} \sigma_2 C(t_2) \widetilde B(t_2) \sigma_1 +
   C(t_2) \widetilde B(t_2) \sigma_2 - (\lambda \sigma_2 + \gamma_*) \big) - \\
- \sigma_1^{-1} (\sigma_2\lambda + \gamma_*) \big( I + C(t_2) (\lambda I - A_\pi)^{-1} \widetilde B \sigma_1\big) + \\
+ \big( I + C(t_2) (\lambda I - A_\pi)^{-1} \widetilde B \big) \sigma_1^{-1} (\sigma_2\lambda + \gamma) = \\
= -\sigma_1^{-1} \gamma_* + \sigma_1^{-1} \gamma - \sigma_1^{-1} \sigma_2 C(t_2) \widetilde B(t_2) \sigma_1 +
    C(t_2) \widetilde B(t_2) \sigma_2 + \\
+ C(t_2) (\lambda I - A_\pi)^{-1} \widetilde B(t_2) \big( - \sigma_2 C(t_2) \widetilde B(t_2) \sigma_1 +
    \sigma_1 C(t_2) \widetilde B(t_2)\sigma_2 - \gamma_* + \gamma \big) = 0
\end{array} \]
and the theorem follows. \qed

Using this theorem we obtain an alternative proof of realization of theorem \ref{tm:NonConsRealiz}:
\begin{cor}
Let $S(\lambda,t_2)$ be a function as above. Suppose that at $t_2^0$ there is a realization
\[ S(\lambda, t_2^0) = I + C(t_2^0) (\lambda I - A_1)^{-1} B(t_2^0) \sigma_1.
\]
Then there exists a realization of $S(\lambda,t_2)$
\[ S(\lambda,t_2) = I + C(t_2) (\lambda I + A_1)^{-1} B(t_2) \sigma_1
\]
with $C(t_2), B(t_2)$, satisfying the output and input adjoint differential equations with matrix parameter $A_1$.
\end{cor}
\textbf{Proof:} Take $(C(t_2),A_1)$ and $(B(t_2), A_\xi)$, the null-pole triple for the matrix $S(\lambda,t_2)$. Then they can be
chosen so that the output differential equation (with $A_1$) and the adjoint output differential equations (with $A_\xi$) are satisfied.
Define now
\[ X(t_2) = B(t_2) \sigma_1 C(t_2).
\]
Then $X(t_2)$ is a coupling matrix (by simple calculations) for this null-pole data and the corresponding realization
\[ S(\lambda,t_2) = I + C(t_2) (\lambda I + A_1)^{-1} X^{-1} B(t_2) \sigma_1
\]
is unique by theorem \ref{thm:uniqueRealiz}. Moreover, $X(t_2)^{-1} B(t_2)$ satisfies the input adjoint differential equation with
the spectral matrix parameter $A_1$. This finishes the proof. \qed

\section{Hermitian realization theorem}
As a special case of realization theorem discussed in \cite{bib:NonConsTheory}, one can
consider the Hermitian case
\begin{thm} Suppose that $S(\lambda,t_2)$ is as above and additionally satisfies
\begin{enumerate}
    \item $S(\infty,t_2) = \infty$,
    \item $S(\lambda,t_2) \sigma_1^{-1} S^*(\lambda,t_2) \geq \sigma_1^{-1}$ for $\Re \lambda > 0$,
    \item $S(\lambda,t_2) \sigma_1^{-1} S^*(-\bar\lambda,t_2) = \sigma_1^{-1}$ for all $t_2$,
\end{enumerate}
then there exists a conservative vessel with $S(\lambda,t_2)$ as the transfer function. In this case
\[ S(\lambda,t_2) = I + C(t_2) (\lambda I + A_1)^{-1} C^*(t_2) \sigma_1,
\]
where $C(t_2)$ satisfies the output differential equation with the matrix parameter $A_1$.
\end{thm}
\textbf{Proof:} Notice that theorem 6.1.1 of \cite{bib:Inter} means that for each $t_2$, there exists an invertible, Hermitian
matrix $X(t_2)$, satisfying
\begin{equation} \label{eq:ForX} X(t_2) A_1 + A_1^* X(t_2) = -C^*(t_2) \sigma_1 C(t_2)
\end{equation}
and a minimal realization
\[ S(\lambda,t_2) = I + C(t_2) (\lambda I + A_1)^{-1} X^{-1}(t_2) C(t_2)^* \sigma_1, \]
when $X(t_2)$ is the null-pole coupling matrix associated with the right pole pair $(C(t_2),A_1)$ and the left
null pair $(-A_1^*, -C(t_2)\sigma_1)$ for $S(\lambda,t_2)$. Here the adjoint is taken with respect to the standard inner
product on the finite dimensional auxiliary Hilbert space.

Notice that from the theorem it also follows that in order to obtain that
$S(\lambda,t_2) \sigma_1^{-1} S^*(\lambda,t_2) \geq \sigma_1^{-1}$ for $\Re \lambda > 0$, we have to demand
positive-definiteness of $X(t_2)$, because in using this realization and formula (\ref{eq:ForX})
(see theorem 6.2.2 from \cite{bib:Inter})
\[
S(\lambda,t_2) \sigma_1^{-1} S^*(-\bar\lambda,t_2) - \sigma_1^{-1} =
C(t_2) (\lambda I - A_1)^{-1} X^{-1}(t_2) (\lambda I - A_1^*)^{-1}C(t_2)^*
\]
and the positive definiteness of $X(t_2)$ is immediate from minimality of the realization.

Now we would like to find a unitary equivalence to obtain the so called \textit{colligation condition}.
Define $Y(t_2) = \sqrt{X(t_2)}$ and define
\begin{equation} \label{eq:KinEq}
\begin{array}{lll}
\tilde C(t_2) = C(t_2) Y^{-1}(t_2) \\
\tilde A_1(t_2) = Y^{-1}(t_2) A_1 Y(t_2)
\end{array}
\end{equation}
Define the inner product on the Hilbert space, on which $\tilde A_1(t_2) $ acts as (for $<,>$ -  the standard inner product)
\[ \langle v,u \rangle = <X(t_2) v, u>
\]
Then $Y(t_2)$ is a self adjoint operator, because
\[ \langle Y(t_2) v, u \rangle = <X(t_2) Y(t_2) v, u> = <v, Y(t_2) u> = \langle v, Y(t_2) u \rangle,
\]
and the adjoint of $\tilde A_1$ is evaluated from
\[ \langle \tilde A_1(t_2) v,u \rangle = <X(t_2) Y^{-1}(t_2) A_1 Y(t_2) v,u > =
<X(t_2)^{-1} v, Y(t_2) A_1^* Y^{-1}(t_2) u> = \langle  v, \tilde A_1(t_2)u \rangle
\]
and $\tilde A_1^* = Y(t_2) A_1^* Y^{-1}(t_2)$. Then
\[ \begin{array}{llll}
\langle \tilde A_1(t_2) v,u \rangle + \langle \tilde A_1^*(t_2) v, u \rangle =
    <X(t_2)Y^{-1}(t_2) A_1 Y(t_2) v,u > + <X(t_2) Y(t_2) A_1^* Y^{-1}(t_2) v, u> = \\
<Y(t_2)[A_1 X(t_2) + A_1^* X(t_2)] Y^{-1}(t_2) v,u> = -<Y(t_2) C^*(t_2) \sigma_1 C(t_2)Y^{-1}(t_2) v,u> = \\
-<X(t_2) \tilde C^*(t_2) \sigma_1 \tilde C(t_2) v,u> = -\langle \tilde C^*(t_2) \sigma_1 \tilde C(t_2) v,u \rangle.
\end{array}\]
In other words
\[ \tilde A_1(t_2) + \tilde A_1^*(t_2) = - \tilde C(t_2) \sigma_1 \tilde C(t_2),
\]
which is precisely the colligation condition. Finally,
\[ S(\lambda,t_2) = I + C(t_2) (\lambda I + A_1)^{-1} X(t_2)^{-1} C^*(t_2) \sigma_1 =
                  I + \tilde C(t_2) (\lambda I + \tilde A_1(t_2) )^{-1} \tilde C^*(t_2) \sigma_1.
\]
\qed

\noindent\textbf{Remarks:}
\begin{enumerate}
    \item It is also possible to use kinematic equivalence in order to obtain the same inner Hilbert space
        with the standard inner product on it. This kinematic equivalence is made by means of matrix $Y(t_2)$,
        using formulas (\ref{eq:KinEq}). Notice that in this case the operator $\tilde A_1(t_2)$ is non-constant.
    \item On the other hand, it is possible to keep the operator $A_1$ constant, varying the Hilbert spaces
        with $t_2$. In this case, it is enough to define the inner product as (for $v \in \mathcal H$, $u \in \mathcal E$)
        \[ \langle v , u\rangle = <X(t_2)^{-1} v,u>
        \]
        Then $C^*(t_2) u = X(t_2)^{-1}u$ and the formula
        \[ S(\lambda,t_2) = I + C(t_2) (\lambda I + A_1)^{-1} C^*(t_2) \sigma_1
        \]
        is obtained.
\end{enumerate}



\begin{thebibliography}{99}
\bibitem[ABP]{bib:AlBallPerez} D. Alpay, J. A. Ball, Y. Peretz, \textit{System theory, operator models and
scattering: the time-varying case}, J. Operator Theory 47(2002), pp. 245-286.
\bibitem[BC]{bib:BallCohen} J.A. Ball, N.Cohen, \textit{De Branges-Rovnyak operator models and systems theory: a survey},
Operator Theory: Advances and Applications, Vol. 50, Birkhauser-Verlag, Basel 1991.
\bibitem[BV]{bib:Over2Dsys} J.A. Ball and V. Vinnikov, \textit{Overdetermined Multidimensional
Systems: State Space and Frequency Domain Methods}, Mathematical Systems Theory, D.~Gilliam and J. Rosenthal, eds.,
Inst. Math. and its Appl. Volume Series, Vol. 134,
Springer-Verlag, New York (2003), 63--120
\bibitem[BGR]{bib:Inter} J. A. Ball, I. Gohberg, L. Rodman, \textit{Interpolation of Rational Matrix Functions},
Operator theory: advances and applications, Birkhauser, 1990.
\bibitem[BFKD]{bib:BGKD} H.~Bart, I.~Gohberg, M.A.~Kaashok, P.~Van~Dooren. \textit{Factorizations of transfer functions}, SIAM J. Control Optim.,
18(6): 675--696, 1980.
\bibitem[CoLe]{bib:CoddLev} E.A Coddington, N. Levinson, \textit{Theory of ordinary differential equation},
Mc-Graw Hill, 1955.
\bibitem[B]{bib:JordanBrodskii}M.S. Brodskii, \textit{Triangular and Jordan representation of linear operators},
Moscow, Nauka, 1969 (Russian); English trans.: Amer. Math. Soc., Providence, 1974.
\bibitem[Ga]{bib:Gauchman} H. Gauchman, \textit{Curvature colligations}. Integral Equations Operator Theory 7 (1984), no. 1, 45--59
\bibitem[Li]{bib:Vortices} M.S. Liv\v sic, \textit{Vortices of 2D systems}, Operator Theory: Advances and
Applications, Vol. 123, Birkhauser-Verlag, Basel 2001.
\bibitem[LiB]{bib:SpectrAnal} M.S. Brodskii, M. Liv\v sic,  \textit{Spectral analysis of non-selfadjoint operators and intermediate systems}. Amer. Math. Soc. Transl. (2) 13 1960 265--346. 46.00  
\bibitem[LKMV]{bib:TheoryNonComm} M.S. Liv\v sic, N.Kravitsky, A.S. Markus and V. Vinnikov, \textit{Theory of Commuting Nonselfadjoint Operators},
Kluwer Acad. Press, 1995.
\bibitem[La]{bib:Latushkin} C. Chicone, Y. Latushkin,
\textit{Evolution Semigroups in Dynamic Systems and Differential Equations},
Math. Surveys and Monographs, Vol. 70, Amer. Math. Soc. 1999.
\bibitem[MV]{bib:NonConsTheory} A.P. Melnikov, V. Vinnikov, \textit{Overdetermined $2D$ Systems Invariant in One Direction and Their
 Transfer Functions}, ArXiv file: arXiv:0812.3779v1, see http://arxiv.org/abs/0812.3779.
\bibitem[NF]{bib:HarmAnal} B. Sz.-Nagy, C. Foia's, \textit{Harmonic analysis of operators on Hilbert space}, Translated
from the French and revised, North-Holland, Amsterdam, 1970. 
\bibitem[N]{bib:LaxEquation} A. C. Newell, \textit{Solitons in mathematics and physics}, Society for Industr. and
Appl. Math., 1985.
\bibitem[P]{bib:Potapov} V.P. Potapov, \textit{On the multiplicative structure of $J$-nonexpanding matrix functions},
Trudy Moskov. Mat. Obschestva, 4 (1955), 125-236 (Russian); English transl.: Amer. Math. Soc. Transl. (2),
15 (1960), 131-243.
\bibitem[V]{bib:VinnMSRI}  V. Vinnikov, \textit{Commuting operators and function theory on a Riemann surface},
Holomorphic spaces (Berkeley, CA, 1995), 445--476, Math. Sci. Res. Inst. Publ., 33,
Cambridge Univ. Press, Cambridge, 1998.
\end{thebibliography}
\end{document}